\documentclass[11pt]{article}

	\newcommand{\mainTitle}{Identities involving cyclic sums of regularized multiple zeta values each of depth less than $5$}

	\newcommand{\authorName}{Machide, Tomoya}
	\newcommand{\organizationNameFst}{National Institute of Informatics}
	\newcommand{\placeAddressFst}{2-1-2 Hitotsubashi, Chiyoda-ku, Tokyo 101-8430, Japan}
	\newcommand{\emailAddressFst}{machide@nii.ac.jp}
	\newcommand{\organizationNameScd}{JST, ERATO, Kawarabayashi Large Graph Project}
	\newcommand{\departmentNameScd}{c/o Global Research Center for Big Data Mathematics}
	\newcommand{\placeAddressScd}{NII, 2-1-2 Hitotsubashi, Chiyoda-ku, Tokyo 101-8430, Japan}
	
	\newcommand{\MSCname}{11M32(Primary), 16S34,20C05(Secondary)} %
	\newcommand{\keyWord}{regularized multiple zeta value, cyclic sum, symmetric sum, group ring of symmetric group}

\usepackage{graphicx}
\usepackage{geometry}               
\usepackage{ascmac}
\usepackage{amsmath}
\usepackage{amssymb}
\usepackage{amsthm}
\usepackage[hypertex, colorlinks=true, bookmarks=true]{hyperref}
\usepackage{arydshln}
\usepackage{mathabx}
\usepackage{multicol}

	\DeclareMathOperator*{\OPlus}{\bigoplus}

	\newcommand{\RCTI}{\!{\ooalign{\hfill$\lceil$\hfill\crcr\hfill$\>\>\rfloor$\hfill}}}

	\newcommand{\RCTbTmp}{{\ooalign{\hfill${\ooalign{\hfill${\ooalign{\hfill$\!\mid$\hfill\crcr\hfill$\!\mid$\hfill}}$\hfill\crcr\hfill$\!\!\mid$\hfill}}$\hfill\crcr\hfill\hspace{-5pt}$\mid$\hfill}}}
	\newcommand{\RCTbI}{{\ooalign{\hfill$\RCTbTmp$\hfill\crcr\hfill$\hspace{1.5pt}\RCTbTmp$\hfill}}}
	\newcommand{\nbk}[3]{#1#3#2}		
	\newcommand{\bgbk}[3]{\bigl{#1}#3\bigr{#2}}	
	\newcommand{\Bgbk}[3]{\Bigl{#1}#3\Bigr{#2}}			
	\newcommand{\bggbk}[3]{\biggl{#1}#3\biggr{#2}}			
	\newcommand{\Bggbk}[3]{\Biggl{#1}#3\Biggr{#2}}
	\newcommand{\autobk}[3]{\left#1#3\right#2}
	\newcommand{\nbkD}[5]{#1#2#5#3#4}		
	\newcommand{\bgbkD}[5]{\bigl{#1}\bigl{#2}#5\bigr{#3}\bigr{#4}}	
	\newcommand{\BgbkD}[5]{\Bigl{#1}\Bigl{#2}#5\Bigr{#3}\Bigr{#4}}	
	\newcommand{\bggbkD}[5]{\biggl{#1}\biggl{#2}#5\biggr{#3}\biggr{#4}}	
	\newcommand{\BggbkD}[5]{\Biggl{#1}\Biggl{#2}#5\Biggr{#3}\Biggr{#4}}	
	\newcommand{\autobkD}[5]{\left#1\left#2#5\right#3\right#4}	
	\newcommand{\mcbk}[4][?]{\ifx n#1\nbk{#2}{#3}{#4}\else\ifx b#1\bgbk{#2}{#3}{#4}\else\ifx B#1\Bgbk{#2}{#3}{#4}\else\ifx g#1\bggbk{#2}{#3}{#4}\else\ifx G#1\Bggbk{#2}{#3}{#4}\else\ifx a#1\autobk{#2}{#3}{#4}\else\ifx !#1{#4}\else#4\fi\fi\fi\fi\fi\fi\fi}
	\newcommand{\mcbkD}[4][?]{\ifx n#1\nbkD{#2}{#2}{#3}{#3}{#4}\else\ifx b#1\bgbkD{#2}{#2}{#3}{#3}{#4}\else\ifx B#1\BgbkD{#2}{#2}{#3}{#3}{#4}\else\ifx g#1\bggbkD{#2}{#2}{#3}{#3}{#4}\else\ifx G#1\BggbkD{#2}{#2}{#3}{#3}{#4}\else\ifx a#1\autobkD{#2}{#2}{#3}{#3}{#4}\else\ifx !#1{#4}\else#4\fi\fi\fi\fi\fi\fi\fi}
	\newcommand{\nsgsb}[1]{#1}		
	\newcommand{\bgsgsb}[1]{\big{#1}}	
	\newcommand{\Bgsgsb}[1]{\Big{#1}}			
	\newcommand{\bggsgsb}[1]{\bigg{#1}}			
	\newcommand{\Bggsgsb}[1]{\Bigg{#1}}
	\newcommand{\mcsgsb}[2][?]{\ifx n#1\nsgsb{#2}\else\ifx b#1\bgsgsb{#2}\else\ifx B#1\Bgsgsb{#2}\else\ifx g#1\bggsgsb{#2}\else\ifx G#1\Bggsgsb{#2}\else#2\fi\fi\fi\fi\fi}
	\newcommand{\myEqSpace}{\,} 	\newlength{\myEqSpaceLen} 	
	\newcommand{\mLt}[1]{\widetilde{#1}}
	\newcommand{\oLt}[1]{\widecheck{#1}}

	\newcommand{\bkR}[2][n]{\mcbk[#1]{(}{)}{#2}}						
	\newcommand{\bkS}[2][n]{\mcbk[#1]{[}{]}{#2}}						
	\newcommand{\bkB}[2][n]{\mcbk[#1]{\{}{\}}{#2}}						
	\newcommand{\bkA}[2][n]{\mcbk[#1]{\langle}{\rangle}{#2}}				
	\newcommand{\bkAll}[4][n]{\mcbk[#1]{#2}{#3}{#4}}
		
	\newcommand{\nFc}[3][n]{#2\bkR[#1]{#3}}
	
	\newcommand{\idFc}[4][n]{\id{#2}{#3}\bkR[#1]{#4}}
	\newcommand{\pwFc}[4][n]{\pw{#2}{#3}\bkR[#1]{#4}}
	\newcommand{\ipFc}[5][n]{\ip{#2}{#3}{#4}\bkR[#1]{#5}}
	\newcommand{\nFcR}[3][n]{(#2)\bkR[#1]{#3}}
		\newcommand{\Fc}{\nFc}
			\newcommand{\Fcr}{\nFcR}

	\newcommand{\alp}{\alpha}
	\newcommand{\gam}{\gamma} 
	\newcommand{\Gam}{\Gamma}

	\newcommand{\sig}{\sigma}
	\newcommand{\Sig}{\Sigma}

	\newcommand{\modd}[1][\ ]{\mathrm{mod}#1}
	\newcommand{\mo}{(-1)}
	\newcommand{\mVert}[1][n]{{\,\mcsgsb[#1]{\vert}\,}}

	\newcommand{\SetO}[2][n]{\bkB[#1]{#2}}
	\newcommand{\SetT}[3][n]{\bkB[#1]{#2\mVert#3}}
		\newcommand{\Set}{\SetO}

	\newcommand{\GpO}[2][n]{\bkA[#1]{#2}}
		\newcommand{\Gp}{\GpO}

	\newcommand{\setN}{\mathbb{N}}
	\newcommand{\setZ}{\mathbb{Z}} 
	\newcommand{\setQ}{\mathbb{Q}}	
		
	\newcommand{\setR}{\mathbb{R}}

	\newcommand{\matI}[1][?]{\ifx #1?I\else I^{(#1)}\fi}	
	
	\newcommand{\gpSym}[2][?]{S^{(#2)}}			
	\newcommand{\gpAlt}[2][?]{A^{(#2)}}				
	\newcommand{\gpKleinF}[1][?]{V}

	\newcommand{\gpu}[1][?]{\ifx?#1e\else e^{(#1)}\fi}	
	
	\newcommand{\vPack}[1][10]{\vspace{-#1pt}}
	\newcommand{\vWiden}[1][10]{\vspace{#1pt}}
	\newcommand{\lnA}[1][]{&  &}
	\newcommand{\lnP}[1]{\myEqSpace#1\myEqSpace}
	\newcommand{\lnAP}[2][]{& #2 &}
	\newcommand{\lnAH}[1][\nonumber]{#1\\ & &}

		\newcommand{\slnAH}[1][?]{\\}
		
	\newcommand{\refEq}[1]{(\ref{#1})}	
	\newcommand{\refEqA}[1]{(#1)}	
	\newcommand{\pcstSpForRefThm}{\;}		
	\newcommand{\refHL}[2]{#1\pcstSpForRefThm\ref{#2}}		
	\newcommand{\refHLm}[3][?]{\ifx?#1#2\pcstSpForRefThm#3\else#2#3\fi}
	\newcommand{\refThm}[2][?]{\ifx?#1\refHL{Theorem}{#2}\else\ifx s#1\refHL{Theorems}{#2}\else{[argument error]}\fi\fi}
	\newcommand{\refProp}[2][?]{\ifx?#1\refHL{Proposition}{#2}\else\ifx s#1\refHL{Propositions}{#2}\else{[argument error]}\fi\fi}
	\newcommand{\refLem}[2][?]{\ifx?#1\refHL{Lemma}{#2}\else\ifx s#1\refHL{Lemmas}{#2}\else{[argument error]}\fi\fi}
	\newcommand{\refCor}[2][?]{\ifx?#1\refHL{Corollary}{#2}\else\ifx s#1\refHL{Corollaries}{#2}\else{[argument error]}\fi\fi}
	\newcommand{\refDef}[2][?]{\ifx?#1\refHL{Definition}{#2}\else\ifx s#1\refHL{Definitions}{#2}\else{[argument error]}\fi\fi}
	\newcommand{\refRem}[2][?]{\ifx?#1\refHL{Remark}{#2}\else\ifx s#1\refHL{Remarks}{#2}\else{[argument error]}\fi\fi}
	\newcommand{\refTab}[2][?]{\ifx?#1\refHL{Table}{#2}\else\ifx s#1\refHL{Tables}{#2}\else{[argument error]}\fi\fi}
	\newcommand{\refSec}[2][?]{\ifx?#1\refHL{Section}{#2}\else\ifx s#1\refHL{Sections}{#2}\else{[argument error]}\fi\fi}
	\newcommand{\refApp}[2][?]{\ifx?#1\refHL{Appendix}{#2}\else\ifx s#1\refHL{Appendices}{#2}\else{[argument error]}\fi\fi}
		\newcommand{\refSect}{\refSec}
	\newcommand{\refThmA}[2][?]{\refHLm[#1]{Theorem}{#2}}
	\newcommand{\refPropA}[2][?]{\refHLm[#1]{Proposition}{#2}}

	\newcommand{\refSecA}[2][?]{\refHLm[#1]{Section}{#2}}
	
		\newcommand{\refSectA}{\refSecA}
		
	\newcommand{\Vc}[2][?]{\ifx ?#1\vec{#2}\else\ifx l#1\overrightarrow{#2}\else\ifx b#1{\bf#2}\else[error]\fi\fi\fi}
	
	\newcommand{\nopF}[3][?]{\ifx s#1#2/#3\else\ifx b#1(#2)/#3\else\ifx d#1\dfrac{#2}{#3}\else\ifx t#1\tfrac{#2}{#3}\else\frac{#2}{#3}\fi\fi\fi\fi}

		\newcommand{\opF}{\nopF}

	\newcommand{\myVopLetter}{$\cdot$} \newlength{\myVopLetterHeight}   

	\newcommand{\racFT}[3][?]{#2\vert\bkR[#1]{#3}}		\newcommand{\racFrT}[3][?]{(#2)\vert\bkR[#1]{#3}}				\newcommand{\racFrrT}[3][?]{(#2\vert\bkR[#1]{#3})}	
												\newcommand{\racF}{\racFT}	\newcommand{\racFr}{\racFrT}		\newcommand{\racFrr}{\racFrrT}
	\newcommand{\racFATh}[4][n]{\bkR[#1]{#2\vert#3}(#4)}				
												\newcommand{\racFA}{\racFATh}

	\newcommand{\pw}[3][?]{\ifx!#3{#2}^{#3}\else#2^{#3}\fi}
	\newcommand{\id}[3][?]{#2_{#3}}
	\newcommand{\ip}[4][?]{{#2}_{#3}^{#4}}
	\newcommand{\pwR}[3][a]{\ifx!#1{\bkR[#1]{#2}}^{#3}\else\bkR[#1]{#2}^{#3}\fi}
	\newcommand{\pwB}[3][a]{\ifx!#1{\bkB[#1]{#2}}^{#3}\else\bkB[#1]{#2}^{#3}\fi}
	\newcommand{\pwS}[3][a]{\ifx!#1{\bkS[#1]{#2}}^{#3}\else\bkS[#1]{#2}^{#3}\fi}
					
	\newcommand{\mpIlett}{id}		\newcommand{\mpI}[1][?]{\ifx?#1\mpIlett\else \mpIlett_{#1}\fi}	
	
	\newcommand{\tpT}[3][a]{ {#2}\atop \bkR[#1]{#3} }

	\newcommand{\nSmO}[2][?]{\ifx l#1\sum\limits_{#2}\else\ifx t#1{\textstyle\sum\limits_{#2}}\else\sum_{#2}\fi\fi}
	\newcommand{\nSmT}[3][?]{\ifx l#1\sum\limits_{#2}^{#3}\else\if t#1{\textstyle\sum\limits_{#2}^{#3}}\else\sum_{#2}^{#3}\fi\fi}	
	\newcommand{\nSmN}[1][?]{\ifx l#1\sum\limits\else\ifx t#1{\textstyle\sum\limits}\else\sum\fi\fi}
	\newcommand{\pSm}[2][?]{\ifx t#1 \sum_{#2}^{\prime} \else \sideset{}{^\prime}\sum_{#2} \fi}
	\newcommand{\pSmT}[3][?]{\ifx t#1 \sum_{#2}^{\prime#3} \else \sideset{}{^\prime}\sum_{#2}^{#3} \fi}	
	\newcommand{\pSmN}[1][?]{\ifx t#1 \sum^{\prime} \else \sideset{}{^\prime}\sum \fi}
	\newcommand{\dSm}[2][?]{\ifx t#1 \sum_{#2}^{\dagger} \else \sideset{}{^\dagger}\sum_{#2} \fi}
	\newcommand{\dSmT}[3][?]{\ifx t#1 \sum_{#2}^{\dagger#3} \else \sideset{}{^\dagger}\sum_{#2}^{#3} \fi}	
	\newcommand{\dSmN}[1][?]{\ifx t#1 \sum^{\dagger} \else \sideset{}{^\dagger}\sum \fi}
	\newcommand{\tpTSm}[3][?]{\nSmO[#1]{\tpT{#2}{#3}}}

		\newcommand{\Sm}{\nSmO}			\newcommand{\SmT}{\nSmT}			\newcommand{\SmN}{\nSmN}
		\newcommand{\tpSm}{\tpTSm}
		
	\newcommand{\nPd}[2][?]{\ifx l#1 \prod\limits_{#2} \else \prod_{#2} \fi}
	\newcommand{\nPdT}[3][?]{\ifx l#1 \prod\limits_{#2}^{#3} \else \prod_{#2}^{#3} \fi}

		\newcommand{\PdT}{\nPdT}

	\newcommand{\nOPs}[2][?]{\ifx l#1 \OPlus\limits_{#2} \else \OPlus_{#2} \fi}
	\newcommand{\nOPsT}[3][?]{\ifx l#1 \OPlus\limits_{#2}^{#3} \else \OPlus_{#2}^{#3} \fi}	
	\newcommand{\pOPs}[2][?]{\ifx t#1 \OPlus_{#2}^{\prime} \else \sideset{}{^\prime}\OPlus_{#2} \fi}
	\newcommand{\pOPsT}[3][?]{\ifx t#1 \OPlus_{#2}^{\prime#3} \else \sideset{}{^\prime}\OPlus_{#2}^{#3} \fi}

	\newcommand{\nIs}[2][?]{\ifx l#1 \bigcap\limits_{#2}\else\ifx b#1 \bigcap_{#2}\else{\textstyle\bigcap\limits_{#2}}\fi\fi}
	\newcommand{\nIsT}[3][?]{\ifx l#1 \bigcap\limits_{#2}^{#3}\else\ifx b#1 \bigcap_{#2}^{#3}\else{\textstyle\bigcap\limits_{#2}^{#3}}\fi\fi}	
	\newcommand{\pIs}[2][?]{\ifx t#1 \bigcap_{#2}^{\prime} \else \sideset{}{^\prime}\bigcap_{#2} \fi}
	\newcommand{\pIsT}[3][?]{\ifx t#1 \bigcap_{#2}^{\prime#3} \else \sideset{}{^\prime}\bigcap_{#2}^{#3} \fi}

	\newcommand{\nUn}[2][?]{\ifx L#1 \bigcup\limits_{#2}\else\ifx b#1 \bigcup_{#2}\else\ifx t#1{\textstyle\bigcup_{#2}}\else{\textstyle\bigcup\limits_{#2}}\fi\fi\fi}
	\newcommand{\nUnT}[3][?]{\ifx L#1 \bigcup\limits_{#2}^{#3}\else\ifx b#1 \bigcup_{#2}^{#3}\else\ifx t#1{\textstyle\bigcup_{#2}^{#3}}\else{\textstyle\bigcup\limits_{#2}^{#3}}\fi\fi\fi}	
	\newcommand{\pUn}[2][?]{\ifx t#1 \bigcup_{#2}^{\prime} \else \sideset{}{^\prime}\bigcup_{#2} \fi}
	\newcommand{\pUnT}[3][?]{\ifx t#1 \bigcup_{#2}^{\prime#3} \else \sideset{}{^\prime}\bigcup_{#2}^{#3} \fi}

		\newcommand{\Un}{\nUn}
		\newcommand{\UnT}{\nUnT}

	\newcommand{\nLm}[2][?]{\ifx l#1 \lim\limits_{#2} \else \lim_{#2} \fi}
		
	\newcommand{\glcondEnvLineHead}[1]{ \ifx*#1 \begin{eqnarray*} \else \begin{eqnarray}  \label{#1} \fi }
	\newcommand{\glcondEnvLineTail}[1]{ \ifx*#1 \end{eqnarray*} \else \end{eqnarray} \fi }
	\newcommand{\glcondDis}[1]{\ifx d#1 \displaystyle \fi}
	\newcommand{\glcmdEqShift}{\hspace{-20pt}}
	\newcommand{\glcmdHLineCWiden}{\rule{0cm}{15pt}}	\newcommand{\glcdH}{\glcmdHLineCWiden}
	\newcommand{\lccondPar}[1]{\ifx#1p \\ \fi}

		\newcommand{\envMO}[2][*]{$\ifx d#1 \displaystyle \fi#2$}
		\newcommand{\envMT}[3][*]{$\ifx d#1 \displaystyle \fi#2=#3$}
		\newcommand{\envMTDef}[3][*]{$\ifx d#1 \displaystyle \fi#2:=#3$}
		\newcommand{\envMTPt}[4][*]{$\ifx d#1 \displaystyle \fi#3#2#4$}
			\newcommand{\envM}{\envMT}
			
			\newcommand{\envMPt}{\envMTPt}
		\newcommand{\envMTh}[4][*]{$\ifx d#1 \displaystyle \fi#2=#3=#4$}
		\newcommand{\envMThDef}[4][*]{$\ifx d#1 \displaystyle \fi#2:=#3=#4$}
		\newcommand{\envMThPt}[5][*]{$\ifx d#1 \displaystyle \fi#3#2#4#2#5$}
		\newcommand{\envMThPte}[6][*]{$\ifx d#1 \displaystyle \fi#2#3#4#5#6$}
		\newcommand{\envMF}[5][*]{$\ifx d#1 \displaystyle \fi#2=#3=#4=#5$}
		\newcommand{\envMFPt}[6][*]{$\ifx d#1 \displaystyle \fi#3#2#4#2#5#2#6$}
		
	\newcommand{\envMLineT}[3][*]{ \ifx*#1 \begin{multline*} #2\lnP{=}#3\end{multline*} \else \begin{multline} \label{#1} #2\lnP{=}#3\end{multline} \fi }
	\newcommand{\envMLineTDef}[3][*]{ \ifx*#1 \begin{multline*} #2\lnP{:=}#3\end{multline*} \else \begin{multline} \label{#1} #2\lnP{:=}#3\end{multline} \fi }
	\newcommand{\envMLineTPt}[4][*]{ \ifx*#1 \begin{multline*} #3\lnP{#2}#4\end{multline*} \else \begin{multline} \label{#1} #3\lnP{#2}#4\end{multline} \fi }

		\newcommand{\envHLineT}[3][*]{ \glcondEnvLineHead{#1} #2&=&#3\glcondEnvLineTail{#1} }
		\newcommand{\envHLineTDef}[3][*]{ \glcondEnvLineHead{#1} #2&:=&#3\glcondEnvLineTail{#1} }
		\newcommand{\envHLineTPt}[4][*]{\glcondEnvLineHead{#1} #3&#2&#4\glcondEnvLineTail{#1}}
		
			\newcommand{\envHLine}{\envHLineT}
			\newcommand{\envHLineDef}{\envHLineTDef}
			\newcommand{\envHLinePt}{\envHLineTPt}
			
		\newcommand{\envHLineThPt}[5][*]{\glcondEnvLineHead{#1} #3&#2&#4\nonumber\\&#2&#5\glcondEnvLineTail{#1}}
		
		\newcommand{\envHLineF}[5][*]{ \glcondEnvLineHead{#1} #2&=&#3\nonumber\\&=&#4\\&=&#5\nonumber \glcondEnvLineTail{#1}}
		\newcommand{\envHLineFPt}[6][*]{\glcondEnvLineHead{#1} #3&#2&#4\nonumber\\&#2&#5\\&#2&#6\nonumber\glcondEnvLineTail{#1}}
		
		\newcommand{\envHLineCFLaa}[5][*]{\envPLine[#1]{#2\lnP{=}#3\qquad\text{and}\qquad#4\lnP{=}#5}}
				
		\newcommand{\envHLineCSNme}[7][*]{\begin{eqnarray} #2&=&#3\\\glcdH#4&=&#5\\\glcdH#6&=&#7\end{eqnarray}}
		
		\newcommand{\envHLineCSNmePt}[8][*]{\glcondEnvLineHead{#1} #3&#2&#4\\\glcdH#5&#2&#6\\\glcdH#7&#2&#8\glcondEnvLineTail{#1} }
		
		\newcommand{\envPLine}[2][*]{\glcondEnvLineHead{#1} #2\glcondEnvLineTail{#1}}
		\newcommand{\envPLineNme}[2][*]{\begin{eqnarray}#2\end{eqnarray}}
		
		\newcommand{\envPLineTNme}[3][*]{\begin{eqnarray} #2\\#3\end{eqnarray}}

		\newcommand{\envOTLine}[4][*]{\glcondEnvLineHead{#1} #2\lnAP{=}#3\lnP{=}#4\glcondEnvLineTail{#1}}
		
		\newcommand{\envOTLinePt}[5][*]{\glcondEnvLineHead{#1} #3\lnP{#2}#4\lnP{#2}#5\glcondEnvLineTail{#1}}

	\newcommand{\lcparaCase}{\vspace{3pt}}

	\newcommand{\envSMatT}[3][a]{\autobk{(}{)}{\begin{smallmatrix}#2\\#3\end{smallmatrix}}}

		\newcommand{\SMatT}{\envSMatT}

	\newcommand{\matu}[1][?]{\ifx#1?I\else I_{#1}\fi}

	\newcommand{\vlA}[2][n]{\bkAll[#1]{|}{|}{#2}}
				
	\newcommand{\nmSet}[2][n]{\vlA[#1]{#2}}	
		\newcommand{\vlNS}{\nmSet}

	\newcommand{\degLett}{\deg}
	\newcommand{\degg}[2][?]{\Fc[#1]{\degLett}{#2}}

	\newcommand{\sbLandau}[2][n]{\Fc[#1]{O}{#2}}
		
		\newcommand{\sbL}{\sbLandau}

	\newcommand{\exx}[2][n]{ \Fc[#1]{\exp}{#2} }
	\newcommand{\exN}[1]{ e^{#1}}

	\newcommand{\fcGam}[2][n]{\Fc[#1]{\Gam}{#2}}

		\newcommand{\fcG}{\fcGam}
	
	\newcommand{\ctG}[1][?]{\gam}

	\newcommand{\cTxT}[2]{\textcolor{#1}{#2}}
	
		\newcommand{\cTx}{\cTxT}
	\newcommand{\alTx}[1]{\cTx{red}{#1}}
	\newcommand{\rmTx}[1]{\cTx{blue}{#1}}
	\newcommand{\ntTx}[1]{\cTx{Magenta}{#1}}
	
	\newcommand{\cmoTx}[1]{\cTx{Gray}{#1}}
	\newcommand{\sTx}[2][?]{ \ifx t#1{\tiny #2} \else \ifx s#1{\scriptsize #2} \else \ifx f#1{\footnotesize #2} \else \ifx S#1{\small #2} \else \ifx n#1{\normalsize #2} \else \ifx l#1{\large #2} \else \ifx L#1{\Large #2} \else \ifx R#1{\LARGE #2} \else \ifx h#1{\huge #2} \else \ifx H#1{\Huge #2} \else \ifx ?#1 #2 \else #2 \fi\fi\fi\fi\fi\fi\fi\fi\fi\fi\fi }
	\newcommand{\bfTx}[1]{{\bf#1}}
	\newcommand{\usbMTx}[3][?]{\underset{#2}{\underbrace{#3}}}

		\newcommand{\usbTx}{\usbMTx}
		
	\newcommand{\raTx}[3][?]{\raisebox{#2pt}[0pt][0pt]{\ifx d#1\displaystyle\fi#3}}
	\newcommand{\raMTx}[3][?]{\raisebox{#2pt}[0pt][0pt]{$\ifx d#1\displaystyle\fi#3$}}

	\newcommand{\roMTx}[3][?]{\rotatebox[origin=c]{#2}{$#3$}}
	\newcommand{\orMTx}[3][?]{{\ooalign{\hfill$#2$\hfill\crcr\hfill$#3$\hfill}}}

		\newcommand{\envHLineTCl}[3][a]{ \ifx a#1\alTx{\envHLineT{#2}{#3}} \else \ifx r#1 \rmTx{\envHLineT{#2}{#3}} \else\ifx n#1 \ntTx{\envHLineT{#2}{#3}}\else\ifx c#1 \cmoTx{\envHLineT{#2}{#3}} \else \text{[argument error]} \fi\fi\fi\fi \vPack[18] }

		\newcommand{\envHLineCSClPart}[8][?]{\ifx*#1 \begin{eqnarray*} \else \begin{eqnarray}  \label{#1}  \fi \alTx{#3}&#2&\alTx{#4}\\\glcdH#5&#2&#6\nonumber\\\glcdH\alTx{#7}&#2&\alTx{#8}\nonumber\glcondEnvLineTail{*}}
		
			\newcommand{\HLineCTCl}[3][?]{\alTx{#2}&=&\alTx{#3}\nonumber \ifx#1p \\\glcdH \fi}
			\newcommand{\HLineCTClDef}[3][?]{\alTx{#2}&:=&\alTx{#3}\nonumber \ifx#1p \\\glcdH \fi}
			\newcommand{\HLineCFCl}[5][?]{\alTx{#2}&=&\alTx{#3}\nonumber\\\glcdH#4&=&#5\nonumber \ifx#1p \\\glcdH \fi}
			\newcommand{\HLineCFClDef}[5][?]{\alTx{#2}&:=&\alTx{#3}\nonumber\\\glcdH#4&:=&#5\nonumber \ifx#1p \\\glcdH \fi}
			\newcommand{\HLineCSCl}[7][?]{\alTx{#2}&=&\alTx{#3}\nonumber\\\glcdH#4&=&#5\nonumber\\\glcdH\alTx{#6}&=&\alTx{#7}\nonumber\ifx#1p \\\glcdH \fi}
			\newcommand{\HLineCSClDef}[7][?]{\alTx{#2}&:=&\alTx{#3}\nonumber\\\glcdH#4&:=&#5\nonumber\\\glcdH\alTx{#6}&:=&\alTx{#7}\nonumber\ifx#1p \\\glcdH \fi}									
			\newcommand{\HLineCECl}[9][?]{\alTx{#2}&=&\alTx{#3}\nonumber\\\glcdH#4&=&#5\nonumber\\\glcdH\alTx{#6}&=&\alTx{#7}\nonumber\\\glcdH#8&=&#9\nonumber\ifx#1p \\\glcdH \fi}
			\newcommand{\HLineCEClDef}[9][?]{\alTx{#2}&:=&\alTx{#3}\nonumber\\\glcdH#4&:=&#5\nonumber\\\glcdH\alTx{#6}&:=&\alTx{#7}\nonumber\\\glcdH#8&:=&#9\nonumber\ifx#1p \\\glcdH \fi}
		
	\newcommand{\envCenter}[2][*]{\ifx*#1\begin{center}\else\begin{center}[#1]\fi #2\end{center}}
	\newcommand{\envFlushleft}[2][*]{\ifx*#1\begin{flushleft}\else\begin{flushleft}[#1]\fi #2\end{flushleft}}
	\newcommand{\envFlushright}[2][*]{\ifx*#1\begin{flushright}\else\begin{flushright}[#1]\fi #2\end{flushright}}
		
	\newcommand{\envItemIm}[2][*]{\ifx*#1\begin{itemize}\else\begin{itemize}[#1]\fi #2\end{itemize}}
	\newcommand{\envItemDp}[2][*]{\ifx*#1\begin{description}\else\begin{description}[#1]\fi #2\end{description}}
	\newcommand{\envItemEm}[2][*]{\ifx*#1\begin{enumerate}\else\begin{enumerate}[#1]\fi #2\end{enumerate}}

	\newcommand{\envMultCol}[3][*]{\ifx1#2#3\else\begin{multicols}{#2}\ifx*#1\else\mbox{}\vspace{-#1pt}\fi#3\end{multicols}\fi}
	
\theoremstyle{plain}
\newtheorem{theorem}{THEOREM}[section]
\newtheorem{proposition}[theorem]{PROPOSITION}
\newtheorem{lemma}[theorem]{LEMMA}
\newtheorem{corollary}[theorem]{COROLLARY}
\theoremstyle{definition}

\theoremstyle{remark}
\newtheorem{remark}[theorem]{REMARK}
\theoremstyle{plain}

\theoremstyle{definition}

\theoremstyle{remark}

\theoremstyle{plain}
\newtheorem{theoremN}{THEOREM}[section]

\theoremstyle{definition}

\theoremstyle{remark}
\newtheorem{remarkN}[theoremN]{REMARK}
\allowdisplaybreaks[4]
\numberwithin{equation}{section}

	\newcommand{\lccondBibitem}[3][]{ \if ?#2 \bibitem{#3} \else \bibitem[#2]{#3} \fi}
	\newcommand{\refPaper}[8][?]{
			\lccondBibitem{#1}{#2}
				#3,			
				\emph{#4}, 	
				#5\ 			
				{\bf #6}		
				(#7),			
				#8.			
		}
	\newcommand{\refPreprint}[6][?]{
			\lccondBibitem{#1}{#2}
				#3,			
				\emph{#4}, 	
				preprint; #5,	
				#6.			
		}
	\newcommand{\refPaperRep}[9][?]{
			\lccondBibitem{#1}{#2}
				#3,			
				\emph{#4}, 	
				#5\ 			
				{\bf #6}		
				(#7),			
				#8			
				; reprinted in #9	
		}

	\newcommand{\refPaperAlm}[5][?]{
			\lccondBibitem{#1}{#2}
				#3,	 		
				\emph{#4}, 	
				#5		
		}

	\newcommand{\etalTx}[2][?]{#2 \emph{et al.}\!}
	\newcommand{\respTx}[1]{(resp. #1)}

	\newcommand{\glcondEnvLineTailPd}[1]{.\ifx*#1 \end{eqnarray*} \else \end{eqnarray} \fi  }
	\newcommand{\glcondEnvLineTailCm}[1]{,\ifx*#1 \end{eqnarray*} \else \end{eqnarray} \fi }
	\newcommand{\prcondEnvEqSpHead}[1]{ \ifx*#1 \begin{equation*}[ERROR] \else \begin{equation}  \label{#1} \fi  }
	\newcommand{\prcondEnvEqSpTail}[1]{\ifx*#1 [ERROR]\end{equation*} \else \end{equation} \fi }

	\newcommand{\envProof}[2][?]{ \par\mbox{}\vspace{-5pt}\\ \ifx?#1\emph{Proof.}\else\emph{Proof of #1.}\fi \ #2 \hfill $\Box$\\ \par}
	
		\newcommand{\envLineTPd}[3][*]{ \glcondEnvLineHead{#1} & &\glcmdEqShift#2\nonumber\\&=&#3 \glcondEnvLineTailPd{#1} }
		\newcommand{\envLineTDefPd}[3][*]{ \glcondEnvLineHead{#1} & &\glcmdEqShift#2\nonumber\\&:=&#3 \glcondEnvLineTailPd{#1} }
		\newcommand{\envLineTCm}[3][*]{ \glcondEnvLineHead{#1} & &\glcmdEqShift#2\nonumber\\&=&#3 \glcondEnvLineTailCm{#1} }
		
			\newcommand{\envLinePd}{\envLineTPd}
			\newcommand{\envLineDefPd}{\envLineTDefPd}
			\newcommand{\envLineCm}{\envLineTCm}
			
		\newcommand{\envLineTCmPt}[4][*]{\glcondEnvLineHead{#1} & &\glcmdEqShift#3\nonumber\\&#2&#4 \glcondEnvLineTailCm{#1}}

			\newcommand{\envLineCmPt}{\envLineTCmPt}

		\newcommand{\envLineThPd}[4][*]{ \glcondEnvLineHead{#1} & &\glcmdEqShift#2\nonumber\\&=&#3 \\&=&#4\nonumber \glcondEnvLineTailPd{#1} }
		\newcommand{\envLineThCm}[4][*]{ \glcondEnvLineHead{#1} & &\glcmdEqShift#2\nonumber\\&=&#3 \\&=&#4\nonumber \glcondEnvLineTailCm{#1} }

		\newcommand{\envLineFPd}[5][*]{ \glcondEnvLineHead{#1} & &\glcmdEqShift#2\nonumber\\&=&#3\nonumber \\&=&#4 \\&=&#5\nonumber \glcondEnvLineTailPd{#1} }
		
		\newcommand{\envHLineTPd}[3][*]{ \glcondEnvLineHead{#1} #2&=&#3\glcondEnvLineTailPd{#1} }
		\newcommand{\envHLineTDefPd}[3][*]{ \glcondEnvLineHead{#1} #2&:=&#3\glcondEnvLineTailPd{#1} }
		\newcommand{\envHLineTCm}[3][*]{ \glcondEnvLineHead{#1} #2&=&#3\glcondEnvLineTailCm{#1} }
		\newcommand{\envHLineTCmDef}[3][*]{ \glcondEnvLineHead{#1} #2&:=&#3\glcondEnvLineTailCm{#1} }
		\newcommand{\envHLineTCmPt}[4][*]{\glcondEnvLineHead{#1} #3&#2&#4\glcondEnvLineTailCm{#1}}					
		\newcommand{\envHLineTPdPt}[4][*]{\glcondEnvLineHead{#1} #3&#2&#4\glcondEnvLineTailPd{#1}}

			\newcommand{\envHLinePd}{\envHLineTPd}
			\newcommand{\envHLineDefPd}{\envHLineTDefPd}
			\newcommand{\envHLineCm}{\envHLineTCm}
			\newcommand{\envHLineCmDef}{\envHLineTCmDef}
			\newcommand{\envHLineCmPt}{\envHLineTCmPt}
			\newcommand{\envHLinePdPt}{\envHLineTPdPt}

		\newcommand{\envHLineThPd}[4][*]{ \glcondEnvLineHead{#1} #2&=&#3\nonumber\\&=&#4\glcondEnvLineTailPd{#1}	 }
		
		\newcommand{\envHLineThCm}[4][*]{ \glcondEnvLineHead{#1} #2&=&#3\nonumber\\&=&#4\glcondEnvLineTailCm{#1}}

		\newcommand{\envHLineFPd}[5][*]{ \glcondEnvLineHead{#1} #2&=&#3\nonumber\\&=&#4 \\&=&#5\nonumber \glcondEnvLineTailPd{#1} }
		\newcommand{\envHLineFCm}[5][*]{ \glcondEnvLineHead{#1} #2&=&#3\nonumber\\&=&#4 \\&=&#5\nonumber \glcondEnvLineTailCm{#1} }
		\newcommand{\envHLineFCmPt}[6][*]{\glcondEnvLineHead{#1} #3&#2&#4\nonumber\\&#2&#5 \\&#2&#6\nonumber\glcondEnvLineTailCm{#1}}

		\newcommand{\envHLineFiPd}[6][*]{ \glcondEnvLineHead{#1} #2&=&#3\nonumber\\&=&#4\nonumber \\&=&#5  \\&=&#6\nonumber\glcondEnvLineTailPd{#1}	}
		\newcommand{\envHLineFiCm}[6][*]{ \glcondEnvLineHead{#1} #2&=&#3\nonumber\\&=&#4\nonumber \\&=&#5  \\&=&#6\nonumber\glcondEnvLineTailCm{#1}}

		\newcommand{\envHLineCFCm}[5][*]{\glcondEnvLineHead{#1} #2&=&#3,\nonumber\\\glcdH#4&=&#5\glcondEnvLineTailCm{#1}}

		\newcommand{\envHLineCFCmNme}[5][*]{\begin{eqnarray} #2&=&#3,\\\glcdH#4&=&#5 \glcondEnvLineTailCm{?} }
		\newcommand{\envHLineCFNmePd}[5][*]{\begin{eqnarray} #2&=&#3,\\\glcdH#4&=&#5 \glcondEnvLineTailPd{?} }
		\newcommand{\envHLineCFCmDefNme}[5][*]{\begin{eqnarray} #2&:=&#3,\\\glcdH#4&:=&#5 \glcondEnvLineTailCm{?} }
		\newcommand{\envHLineCFDefNmePd}[5][*]{\begin{eqnarray} #2&:=&#3,\\\glcdH#4&:=&#5 \glcondEnvLineTailPd{?} }
		\newcommand{\envHLineCFCmPt}[6][*]{\glcondEnvLineHead{#1} #3&#2&#4,\nonumber\\\glcdH#5&#2&#6\glcondEnvLineTailCm{#1}}

		\newcommand{\envHLineCFNmePdPt}[6][*]{\begin{eqnarray}#3&#2&#4,\\\glcdH#5&#2&#6\glcondEnvLineTailPd{?}}
		\newcommand{\envHLineCFCmNmePt}[6][*]{\begin{eqnarray}#3&#2&#4,\\\glcdH#5&#2&#6\glcondEnvLineTailCm{?}}
		\newcommand{\envHLineCFNmePdPte}[7][*]{\begin{eqnarray}#2&#3&#4,\\\glcdH#5&#6&#7\glcondEnvLineTailPd{?}}
		\newcommand{\envHLineCFCmNmePte}[7][*]{\begin{eqnarray}#2&#3&#4,\\\glcdH#5&#6&#7\glcondEnvLineTailCm{?}}

		\newcommand{\envHLineCFLqqPd}[5][*]{\envPLinePd[#1]{#2\lnP{=}#3,\qquad#4\lnP{=}#5}}

		\newcommand{\envHLineCFLaaCm}[5][*]{\envPLineCm[#1]{#2\lnP{=}#3\qquad\text{and}\qquad#4\lnP{=}#5}}		
		\newcommand{\envHLineCFLaaPd}[5][*]{\envPLinePd[#1]{#2\lnP{=}#3\qquad\text{and}\qquad#4\lnP{=}#5}}			
		\newcommand{\envHLineCFLaaCmDef}[5][*]{\envPLineCm[#1]{#2\lnP{:=}#3\qquad\text{and}\qquad#4\lnP{:=}#5}}

		\newcommand{\envHLineCSCm}[7][*]{\glcondEnvLineHead{#1} #2&=&#3,\nonumber\\\glcdH#4&=&#5,\\\glcdH#6&=&#7\nonumber\glcondEnvLineTailCm{#1}}
			
		\newcommand{\envHLineCSNmePd}[7][*]{\begin{eqnarray} #2&=&#3,\\\glcdH#4&=&#5,\\\glcdH#6&=&#7\glcondEnvLineTailPd{?}}
		\newcommand{\envHLineCSDefNmePd}[7][*]{\begin{eqnarray} #2&:=&#3,\\\glcdH#4&:=&#5,\\\glcdH#6&:=&#7\glcondEnvLineTailPd{?}}
		\newcommand{\envHLineCSCmNme}[7][*]{\begin{eqnarray} #2&=&#3,\\\glcdH#4&=&#5,\\\glcdH#6&=&#7\glcondEnvLineTailCm{?}}
		\newcommand{\envHLineCSCmDefNme}[7][*]{\begin{eqnarray} #2&:=&#3,\\\glcdH#4&:=&#5,\\\glcdH#6&:=&#7\glcondEnvLineTailCm{?}}

		\newcommand{\envHLineCSNmePdPt}[8][*]{\begin{eqnarray}#3&#2&#4,\\\glcdH#5&#2&#6,\\\glcdH#7&#2&#8\glcondEnvLineTailPd{?}}
		\newcommand{\envHLineCSCmNmePt}[8][*]{\begin{eqnarray}#3&#2&#4,\\\glcdH#5&#2&#6,\\\glcdH#7&#2&#8\glcondEnvLineTailCm{?}}
		\newcommand{\envHLineCSNmePdPte}[9][*]{\begin{eqnarray}#2&#3&#4,\\\glcdH#5&#6&#7,\\\glcdH#8&#2&#9\glcondEnvLineTailPd{?}}
		\newcommand{\envHLineCSCmNmePte}[9][*]{\begin{eqnarray}#2&#3&#4,\\\glcdH#5&#6&#7,\\\glcdH#8&#2&#9\glcondEnvLineTailCm{?}}
		
		\newcommand{\envHLineCECm}[9][*]{\glcondEnvLineHead{#1} #2&=&#3,\nonumber\\\glcdH#4&=&#5,\nonumber\\\glcdH#6&=&#7,\\\glcdH#8&=&#9\nonumber\glcondEnvLineTailCm{#1}}
		
		\newcommand{\envHLineCENmePd}[9][*]{\begin{eqnarray} #2&=&#3,\\\glcdH#4&=&#5,\\\glcdH#6&=&#7,\\\glcdH#8&=&#9\glcondEnvLineTailPd{?}}
		\newcommand{\envHLineCEDefNmePd}[9][*]{\begin{eqnarray} #2&:=&#3,\\\glcdH#4&:=&#5,\\\glcdH#6&:=&#7,\\\glcdH#8&:=&#9\glcondEnvLineTailPd{?}}
		\newcommand{\envHLineCECmNme}[9][*]{\begin{eqnarray} #2&=&#3,\\\glcdH#4&=&#5,\\\glcdH#6&=&#7,\\\glcdH#8&=&#9\glcondEnvLineTailCm{?}}
		\newcommand{\envHLineCECmDefNme}[9][*]{\begin{eqnarray} #2&:=&#3,\\\glcdH#4&:=&#5,\\\glcdH#6&:=&#7,\\\glcdH#8&:=&#9\glcondEnvLineTailCm{?}}

		\newcommand{\envHLineCECmNmePt}[9]{\begin{eqnarray}#2&#1&#3,\\\glcdH#4&#1&#5,\\\glcdH#6&#1&#7,\\\glcdH#8&#1&#9,\end{eqnarray}}
		
			\newcommand{\pccondPaForPar}[1]{\ifx#1p \\\glcdH \fi}
			\newcommand{\pccondPaForNonnum}[1]{\ifx#1* \nonumber \fi}
			
			\newcommand{\HLineCTPd}[3][?]{#2&=&#3. \nonumber\pccondPaForPar{#1}}

			\newcommand{\HLineCTCmNme}[3][?]{#2&=&#3, \pccondPaForPar{#1}}
			\newcommand{\HLineCTNmePd}[3][?]{#2&=&#3. \pccondPaForPar{#1}}

			\newcommand{\HLineCTCmNmePt}[4][?]{#3&#2&#4, \pccondPaForPar{#1}}
			
			\newcommand{\HLineCFCm}[5][?]{#2&=&#3,\nonumber\\\glcdH#4&=&#5, \pccondPaForNonnum{#1} \pccondPaForPar{#1}} 
			\newcommand{\HLineCFPd}[5][?]{#2&=&#3,\nonumber\\\glcdH#4&=&#5. \pccondPaForNonnum{#1} \pccondPaForPar{#1}} 
			\newcommand{\HLineCFCmNme}[5][?]{#2&=&#3,\\\glcdH#4&=&#5, \pccondPaForPar{#1}}
			\newcommand{\HLineCFNmePd}[5][?]{#2&=&#3,\\\glcdH#4&=&#5. \pccondPaForPar{#1}}

			\newcommand{\HLineCSCm}[7][?]{#2&=&#3,\nonumber\\\glcdH#4&=&#5,\nonumber\\\glcdH#6&=&#7,\pccondPaForPar{#1}}

			\newcommand{\HLineCSCmNme}[7][?]{#2&=&#3,\\\glcdH#4&=&#5,\\\glcdH#6&=&#7,\pccondPaForPar{#1}}	
			\newcommand{\HLineCSNmePd}[7][?]{#2&=&#3,\\\glcdH#4&=&#5,\\\glcdH#6&=&#7.\pccondPaForPar{#1}}

			\newcommand{\HLineCSCmNmePt}[8][?]{#3&#2&#4,\\\glcdH#5&#2&#6,\\\glcdH#7&#2&#8,\pccondPaForPar{#1}}	
			\newcommand{\HLineCSNmePdPt}[8][?]{#3&#2&#4,\\\glcdH#5&#2&#6,\\\glcdH#7&#2&#8.\pccondPaForPar{#1}}			
			\newcommand{\HLineCECm}[9][?]{#2&=&#3,\nonumber\\\glcdH#4&=&#5,\nonumber\\\glcdH#6&=&#7,\nonumber\\\glcdH#8&=&#9,\pccondPaForPar{#1}}

			\newcommand{\HLineCECmNme}[9][?]{#2&=&#3,\\\glcdH#4&=&#5,\\\glcdH#6&=&#7,\\\glcdH#8&=&#9,\pccondPaForPar{#1}}
			
		\newcommand{\envPLinePd}[2][*]{\glcondEnvLineHead{#1} #2\glcondEnvLineTailPd{#1}}
		\newcommand{\envPLineCm}[2][*]{\glcondEnvLineHead{#1} #2\glcondEnvLineTailCm{#1}}
		\newcommand{\envOTLineCm}[4][*]{\glcondEnvLineHead{#1} #2\lnAP{=}#3\lnP{=}#4,\glcondEnvLineTail{#1}}
		\newcommand{\envOTLinePd}[4][*]{\glcondEnvLineHead{#1} #2\lnAP{=}#3\lnP{=}#4.\glcondEnvLineTail{#1}}

		\newcommand{\envOTLineCmPt}[5][*]{\glcondEnvLineHead{#1} #3\lnAP{#2}#4\lnP{#2}#5\glcondEnvLineTailCm{#1}}

			\newcommand{\envOTLineThCm}{\envOTLineCm}
			\newcommand{\envOTLineThPd}{\envOTLinePd}

			\newcommand{\envOTLineThCmPt}{\envOTLineCmPt}

		\newcommand{\envOFLineCm}[5][*]{\glcondEnvLineHead{#1} #2\lnAP{=}#3\lnP{=}#4\lnP{=}#5,\glcondEnvLineTail{#1}}

		\newcommand{\envMOCm}[2][*]{$\ifx d#1 \displaystyle \fi#2$,}
		\newcommand{\envMOPd}[2][*]{$\ifx d#1 \displaystyle \fi#2$.}
		
		\newcommand{\envMTCm}[3][*]{$\ifx d#1 \displaystyle \fi#2=#3$,}
		\newcommand{\envMTPd}[3][*]{$\ifx d#1 \displaystyle \fi#2=#3$.}
		\newcommand{\envMTCmDef}[3][*]{$\ifx d#1 \displaystyle \fi#2:=#3$,}
		\newcommand{\envMTDefPd}[3][*]{$\ifx d#1 \displaystyle \fi#2:=#3$.}
		\newcommand{\envMTCmPt}[4][*]{$\ifx d#1 \displaystyle \fi#3#2#4$,}
		\newcommand{\envMTPdPt}[4][*]{$\ifx d#1 \displaystyle \fi#3#2#4$.}
			\newcommand{\envMCm}{\envMTCm}
			\newcommand{\envMPd}{\envMTPd}

			\newcommand{\envMCmPt}{\envMTCmPt}
			
		\newcommand{\envMThCm}[4][*]{$\ifx d#1 \displaystyle \fi#2=#3=#4$,}
		\newcommand{\envMThPd}[4][*]{$\ifx d#1 \displaystyle \fi#2=#3=#4$.}
		\newcommand{\envMThCmPt}[5][*]{$\ifx d#1 \displaystyle \fi#3#2#4#2#5$,}
		\newcommand{\envMThPdPt}[5][*]{$\ifx d#1 \displaystyle \fi#3#2#4#2#5$.}
		\newcommand{\envMFCm}[5][*]{$\ifx d#1 \displaystyle \fi#2=#3=#4=#5$,}
		\newcommand{\envMFPd}[5][*]{$\ifx d#1 \displaystyle \fi#2=#3=#4=#5$.}
		\newcommand{\envMFCmPt}[6][*]{$\ifx d#1 \displaystyle \fi#3#2#4#2#5#2#6$,}
		\newcommand{\envMFPdPt}[6][*]{$\ifx d#1 \displaystyle \fi#3#2#4#2#5#2#6$.}

		\newcommand{\prcondHLCPNm}{\hspace{-1pt}}
		\newcommand{\envHLineCFCmNm}[5][*]{ \begin{equation}\begin{split} \ifx*#1 \text{[ERROR;need label name]} \else \label{#1} \fi #2&\prcondHLCPNm\lnP{=}\prcondHLCPNm#3,\\#4&\prcondHLCPNm\lnP{=}\prcondHLCPNm#5, \end{split}\end{equation} }
		\newcommand{\envHLineCFNm}[5][*]{ \begin{equation}\begin{split} \ifx*#1 \text{[ERROR;need label name]} \else \label{#1} \fi #2&\prcondHLCPNm\lnP{=}\prcondHLCPNm#3\\#4&\prcondHLCPNm\lnP{=}\prcondHLCPNm#5, \end{split}\end{equation} }
		\newcommand{\envHLineCFNmPd}[5][*]{ \begin{equation}\begin{split} \ifx*#1 \text{[ERROR;need label name]} \else \label{#1} \fi #2&\prcondHLCPNm\lnP{=}\prcondHLCPNm#3,\\#4&\prcondHLCPNm\lnP{=}\prcondHLCPNm#5. \end{split}\end{equation} }
		\newcommand{\envHLineCFCmDefNm}[5][*]{ \begin{equation}\begin{split} \ifx*#1 \text{[ERROR;need label name]} \else \label{#1} \fi #2&\prcondHLCPNm\lnP{:=}\prcondHLCPNm#3,\\#4&\prcondHLCPNm\lnP{:=}\prcondHLCPNm#5, \end{split}\end{equation} }
		\newcommand{\envHLineCFDefNm}[5][*]{ \begin{equation}\begin{split} \ifx*#1 \text{[ERROR;need label name]} \else \label{#1} \fi #2&\prcondHLCPNm\lnP{:=}\prcondHLCPNm#3\\#4&\prcondHLCPNm\lnP{:=}\prcondHLCPNm#5, \end{split}\end{equation} }
		\newcommand{\envHLineCFDefNmPd}[5][*]{ \begin{equation}\begin{split} \ifx*#1 \text{[ERROR;need label name]} \else \label{#1} \fi #2&\prcondHLCPNm\lnP{:=}\prcondHLCPNm#3,\\#4&\prcondHLCPNm\lnP{:=}\prcondHLCPNm#5. \end{split}\end{equation} }
		\newcommand{\envHLineCSCmNm}[7][*]{ \begin{equation}\begin{split} \ifx*#1 \text{[ERROR;need label name]} \else \label{#1} \fi #2&\prcondHLCPNm\lnP{=}\prcondHLCPNm#3,\\#4&\prcondHLCPNm\lnP{=}\prcondHLCPNm#5,\\#6&\prcondHLCPNm\lnP{=}\prcondHLCPNm#7 \end{split}\end{equation} }
		\newcommand{\envHLineCSNm}[7][*]{ \begin{equation}\begin{split} \ifx*#1 \text{[ERROR;need label name]} \else \label{#1} \fi #2&\prcondHLCPNm\lnP{=}\prcondHLCPNm#3\\#4&\prcondHLCPNm\lnP{=}\prcondHLCPNm#5\\#6&\prcondHLCPNm\lnP{=}\prcondHLCPNm#7 \end{split}\end{equation} }
		\newcommand{\envHLineCSNmPd}[7][*]{ \begin{equation}\begin{split} \ifx*#1 \text{[ERROR;need label name]} \else \label{#1} \fi #2&\prcondHLCPNm\lnP{=}\prcondHLCPNm#3,\\#4&\prcondHLCPNm\lnP{=}\prcondHLCPNm#5,\\#6&\prcondHLCPNm\lnP{=}\prcondHLCPNm#7. \end{split}\end{equation} }
		\newcommand{\envHLineCSCmDefNm}[7][*]{ \begin{equation}\begin{split} \ifx*#1 \text{[ERROR;need label name]} \else \label{#1} \fi #2&\prcondHLCPNm\lnP{:=}\prcondHLCPNm#3,\\#4&\prcondHLCPNm\lnP{:=}\prcondHLCPNm#5,\\#6&\prcondHLCPNm\lnP{:=}\prcondHLCPNm#7 \end{split}\end{equation} }
		\newcommand{\envHLineCSDefNm}[7][*]{ \begin{equation}\begin{split} \ifx*#1 \text{[ERROR;need label name]} \else \label{#1} \fi #2&\prcondHLCPNm\lnP{:=}\prcondHLCPNm#3\\#4&\prcondHLCPNm\lnP{:=}\prcondHLCPNm#5\\#6&\prcondHLCPNm\lnP{:=}\prcondHLCPNm#7 \end{split}\end{equation} }
		\newcommand{\envHLineCSDefNmPd}[7][*]{ \begin{equation}\begin{split} \ifx*#1 \text{[ERROR;need label name]} \else \label{#1} \fi #2&\prcondHLCPNm\lnP{:=}\prcondHLCPNm#3,\\#4&\prcondHLCPNm\lnP{:=}\prcondHLCPNm#5,\\#6&\prcondHLCPNm\lnP{:=}\prcondHLCPNm#7. \end{split}\end{equation} }
		\newcommand{\envHLineCECmNm}[9][*]{ \begin{equation}\begin{split} \ifx*#1 \text{[ERROR;need label name]} \else \label{#1} \fi #2&\prcondHLCPNm\lnP{=}\prcondHLCPNm#3,\\#4&\prcondHLCPNm\lnP{=}\prcondHLCPNm#5,\\#6&\prcondHLCPNm\lnP{=}\prcondHLCPNm#7,\\#8&\prcondHLCPNm\lnP{=}\prcondHLCPNm#9,  \end{split}\end{equation} }
		\newcommand{\envHLineCENm}[9][*]{ \begin{equation}\begin{split} \ifx*#1 \text{[ERROR;need label name]} \else \label{#1} \fi #2&\prcondHLCPNm\lnP{=}\prcondHLCPNm#3\\#4&\prcondHLCPNm\lnP{=}\prcondHLCPNm#5\\#6&\prcondHLCPNm\lnP{=}\prcondHLCPNm#7\\#8&\prcondHLCPNm\lnP{=}\prcondHLCPNm#9  \end{split}\end{equation} }
		\newcommand{\envHLineCENmPd}[9][*]{ \begin{equation}\begin{split} \ifx*#1 \text{[ERROR;need label name]} \else \label{#1} \fi #2&\prcondHLCPNm\lnP{=}\prcondHLCPNm#3,\\#4&\prcondHLCPNm\lnP{=}\prcondHLCPNm#5,\\#6&\prcondHLCPNm\lnP{=}\prcondHLCPNm#7,\\#8&\prcondHLCPNm\lnP{=}\prcondHLCPNm#9.  \end{split}\end{equation} }
		\newcommand{\envHLineCECmDefNm}[9][*]{ \begin{equation}\begin{split} \ifx*#1 \text{[ERROR;need label name]} \else \label{#1} \fi #2&\prcondHLCPNm\lnP{:=}\prcondHLCPNm#3,\\#4&\prcondHLCPNm\lnP{:=}\prcondHLCPNm#5,\\#6&\prcondHLCPNm\lnP{:=}\prcondHLCPNm#7,\\#8&\prcondHLCPNm\lnP{:=}\prcondHLCPNm#9,  \end{split}\end{equation} }
		\newcommand{\envHLineCEDefNm}[9][*]{ \begin{equation}\begin{split} \ifx*#1 \text{[ERROR;need label name]} \else \label{#1} \fi #2&\prcondHLCPNm\lnP{:=}\prcondHLCPNm#3\\#4&\prcondHLCPNm\lnP{:=}\prcondHLCPNm#5\\#6&\prcondHLCPNm\lnP{:=}\prcondHLCPNm#7\\#8&\prcondHLCPNm\lnP{:=}\prcondHLCPNm#9  \end{split}\end{equation} }
		\newcommand{\envHLineCEDefNmPd}[9][*]{ \begin{equation}\begin{split} \ifx*#1 \text{[ERROR;need label name]} \else \label{#1} \fi #2&\prcondHLCPNm\lnP{:=}\prcondHLCPNm#3,\\#4&\prcondHLCPNm\lnP{:=}\prcondHLCPNm#5,\\#6&\prcondHLCPNm\lnP{:=}\prcondHLCPNm#7,\\#8&\prcondHLCPNm\lnP{:=}\prcondHLCPNm#9.  \end{split}\end{equation} }
	\newcommand{\envMLineTPd}[3][*]{ \ifx*#1 \begin{multline*} #2\lnP{=}#3.\end{multline*} \else \begin{multline} \label{#1} #2\lnP{=}#3.\end{multline} \fi }
	\newcommand{\envMLineTCm}[3][*]{ \ifx*#1 \begin{multline*} #2\lnP{=}#3,\end{multline*} \else \begin{multline} \label{#1} #2\lnP{=}#3,\end{multline} \fi }
	\newcommand{\envMLineTDefPd}[3][*]{ \ifx*#1 \begin{multline*} #2\lnP{:=}#3.\end{multline*} \else \begin{multline} \label{#1} #2\lnP{:=}#3.\end{multline} \fi }
	\newcommand{\envMLineTCmDef}[3][*]{ \ifx*#1 \begin{multline*} #2\lnP{:=}#3,\end{multline*} \else \begin{multline} \label{#1} #2\lnP{:=}#3,\end{multline} \fi }

	\newcommand{\envCaseTCm}[3][?]{\begin{cases} \glcondDis{#1}#2,\lcparaCase\\\glcondDis{#1}#3,\end{cases}}
	\newcommand{\envCaseTPd}[3][?]{\begin{cases} \glcondDis{#1}#2,\lcparaCase\\\glcondDis{#1}#3.\end{cases}}
	\newcommand{\envCaseThCm}[4][?]{\begin{cases} \glcondDis{#1}#2,\lcparaCase\\\glcondDis{#1}#3,\lcparaCase\\\glcondDis{#1}#4,\end{cases}}
	\newcommand{\envCaseThPd}[4][?]{\begin{cases} \glcondDis{#1}#2,\lcparaCase\\\glcondDis{#1}#3,\lcparaCase\\\glcondDis{#1}#4.\end{cases}}
	
	\newcommand{\envCaseFPd}[5][?]{\begin{cases} \glcondDis{#1}#2,\lcparaCase\\\glcondDis{#1}#3,\lcparaCase\\\glcondDis{#1}#4,\lcparaCase\\\glcondDis{#1}#5.\end{cases}}
	\newcommand{\envCaseFiCm}[6][?]{\begin{cases} \glcondDis{#1}#2,\lcparaCase\\\glcondDis{#1}#3,\lcparaCase\\\glcondDis{#1}#4,\lcparaCase\\\glcondDis{#1}#5,\lcparaCase\\\glcondDis{#1}#6,\end{cases}}

	\newcommand{\envItemEmPa}[2][*]{\vspace{-5pt}\ifx*#1\begin{enumerate}\else\begin{enumerate}[#1]\fi \renewcommand{\itemsep}{0pt}\renewcommand{\itemindent}{0pt}\renewcommand{\theenumi}{\bfTx{(\roman{enumi})}}#2\end{enumerate}}

			\newcommand{\OTLineCThCm}[4][?]{#2&=&#3=#4,\nonumber \ifx#1p \\ \fi}
			\newcommand{\OTLineCThPd}[4][?]{#2&=&#3=#4.\nonumber \ifx#1p \\ \fi}
			\newcommand{\OTLineCThCmDef}[4][?]{#2&:=&#3=#4,\nonumber \ifx#1p \\ \fi}
			\newcommand{\OTLineCThDefPd}[4][?]{#2&:=&#3=#4.\nonumber \ifx#1p \\ \fi}
			\newcommand{\OTLineCSCm}[7][?]{#2&=&#3=#4\nonumber#5&=&#6=#7,\nonumber \ifx#1p \\ \fi}
			\newcommand{\OTLineCSPd}[7][?]{#2&=&#3=#4\nonumber#5&=&#6=#7.\nonumber \ifx#1p \\ \fi}

	\newcommand{\envMLineTCmPt}[4][*]{ \ifx*#1 \begin{multline*} #3\lnP{#2}#4,\end{multline*} \else \begin{multline} \label{#1} #3\lnP{#2}#4,\end{multline} \fi }
	\newcommand{\envMLineTPdPt}[4][*]{ \ifx*#1 \begin{multline*} #3\lnP{#2}#4.\end{multline*} \else \begin{multline} \label{#1} #3\lnP{#2}#4.\end{multline} \fi }

	\DeclareFontEncoding{OT2}{}{}
	\DeclareFontSubstitution{OT2}{cmr}{m}{n}
	\DeclareFontFamily{OT2}{cmr}{\hyphenchar\font45}
	\DeclareFontShape{OT2}{cmr}{m}{n}{<5><6><7><8><9>gen*wncyr <10><10.95><12><14.4><17.28><20.74><24.88>wncyr10}{}
	\DeclareFontShape{OT2}{cmr}{b}{n}{<5><6><7><8><9>gen*wncyb<10><10.95><12><14.4><17.28><20.74><24.88>wncyb10}{}
	\DeclareMathAlphabet{\mathcyr}{OT2}{cmr}{m}{n}
	\DeclareMathAlphabet{\mathcyb}{OT2}{cmr}{b}{n}
	\SetMathAlphabet{\mathcyr}{bold}{OT2}{cmr}{b}{n}

	\newcommand{\sh}{\mathcyr{sh}}
	\newcommand{\shS}[1][\;]{#1\mathcyr{sh}#1}

	\newcommand{\vlWT}[2][n]{\Fc[#1]{\mathrm{wt}}{#2}}
	\newcommand{\vlDP}[2][n]{\Fc[#1]{\mathrm{dep}}{#2}}
	
	\newcommand{\alHfnd}{\mathfrak{H}}		
	\newcommand{\alHR}{\alHfnd^1}		\newcommand{\alHRh}{\alHfnd_*^1}			\newcommand{\alHRs}{\alHfnd_\sh^1}		
	\newcommand{\alHN}{\alHfnd^0}		\newcommand{\alHNh}{\alHfnd_*^0}			\newcommand{\alHNs}{\alHfnd_\sh^0}	
	\newcommand{\alHF}{\alHfnd}							
		
	\newcommand{\fcZBlettHelp}{\bullet}
	\newcommand{\fcZBN}[1][?]{\zeta^\fcZBlettHelp}			
	\newcommand{\fcZHN}[1][?]{\zeta^*}				\newcommand{\fcZHO}[2][n]{\pwFc[#1]{\zeta}{*}{#2}}					\newcommand{\fcZH}{\fcZHO}
	\newcommand{\fcZRN}[1][?]{\zeta^\ddagger}								
	\newcommand{\fcZSN}[1][?]{\zeta^\sh}			\newcommand{\fcZSO}[2][n]{\pwFc[#1]{\zeta}{\sh}{#2}}				\newcommand{\fcZS}{\fcZSO}
	\newcommand{\fcZMN}[1][?]{\zeta_p}								
	\newcommand{\fcMZO}[2][?]{\zeta_{#2}}				\newcommand{\fcMZT}[3][n]{\idFc[#1]{\zeta}{#2}{#3}}					\newcommand{\fcMZ}{\fcMZT}
	\newcommand{\fcMZBO}[2][?]{\zeta_{#2}^\fcZBlettHelp}	\newcommand{\fcMZBT}[3][n]{\ipFc[#1]{\zeta}{#2}{\fcZBlettHelp}{#3}}		\newcommand{\fcMZB}{\fcMZBT}
	\newcommand{\fcMZHO}[2][?]{\zeta_{#2}^*}			\newcommand{\fcMZHT}[3][n]{\ipFc[#1]{\zeta}{#2}{*}{#3}}					\newcommand{\fcMZH}{\fcMZHT}
									
	\newcommand{\fcMZSO}[2][?]{\zeta_{#2}^\sh}			\newcommand{\fcMZST}[3][n]{\ipFc[#1]{\zeta}{#2}{\sh}{#3}}				\newcommand{\fcMZS}{\fcMZST}
								
	\newcommand{\LetHelpFcCh}{\chi}
	\newcommand{\fcChBN}[1][?]{\LetHelpFcCh^\fcZBlettHelp}			
	\newcommand{\fcChHN}[1][?]{\LetHelpFcCh^*}									
	\newcommand{\fcChRN}[1][?]{\LetHelpFcCh^\ddagger}								
	\newcommand{\fcChSN}[1][?]{\LetHelpFcCh^\sh}							

	\newcommand{\fcChBmN}[1][?]{\mLt{\LetHelpFcCh}^\fcZBlettHelp}			
	\newcommand{\fcChHmN}[1][?]{\mLt{\LetHelpFcCh}^*}							
	\newcommand{\fcChRmN}[1][?]{\mLt{\LetHelpFcCh}^\ddagger}								
	\newcommand{\fcChSmN}[1][?]{\mLt{\LetHelpFcCh}^\sh}							
								
	\newcommand{\fcMChBO}[2][?]{\LetHelpFcCh_{#2}^\fcZBlettHelp}	\newcommand{\fcMChBT}[3][n]{\ipFc[#1]{\LetHelpFcCh}{#2}{\fcZBlettHelp}{#3}}		\newcommand{\fcMChB}{\fcMChBT}
				\newcommand{\fcMChHT}[3][n]{\ipFc[#1]{\LetHelpFcCh}{#2}{*}{#3}}					\newcommand{\fcMChH}{\fcMChHT}
									
	\newcommand{\fcMChSO}[2][?]{\LetHelpFcCh_{#2}^\sh}			\newcommand{\fcMChST}[3][n]{\ipFc[#1]{\LetHelpFcCh}{#2}{\sh}{#3}}				\newcommand{\fcMChS}{\fcMChST}

	\newcommand{\fcMChSoO}[2][?]{\oLt{\LetHelpFcCh}_{#2}^\sh}			\newcommand{\fcMChSoT}[3][n]{\ipFc[#1]{\oLt{\LetHelpFcCh}}{#2}{\sh}{#3}}				\newcommand{\fcMChSo}{\fcMChSoT}

	\newcommand{\mpEVlett}{Z}
	\newcommand{\mpEVN}[2][n]{\Fc[#1]{\mpEVlett}{#2}}			\newcommand{\mpEVNi}[3][n]{\idFc[#1]{\mpEVlett}{#2}{#3}}	
	\newcommand{\mpEVH}[2][n]{\pwFc[#1]{\mpEVlett}{*}{#2}}		\newcommand{\mpEVHi}[3][n]{\ipFc[#1]{\mpEVlett}{#2}{*}{#3}}			
	\newcommand{\mpEVS}[2][n]{\pwFc[#1]{\mpEVlett}{\sh}{#2}}		\newcommand{\mpEVSi}[3][n]{\ipFc[#1]{\mpEVlett}{#2}{\sh}{#3}}

	\newcommand{\lbfL}{{\bf l}}
	
	\newcommand{\lbfS}{{\bf s}}
	\newcommand{\lfgLSN}{[\orMTx{\raMTx{6}{\roMTx{90}{\RCTbI}} }{\raMTx{-2}{\roMTx{45}{\;\RCTI}}}]}						\newcommand{\lfgLSO}[2][?]{{\ifx#1?\lfgLSN\else\reflectbox{\lfgLSN}\fi}^{(#2)}}			
	\newcommand{\lfgLSsN}{[\orMTx{\raMTx{6}{\roMTx{90}{\RCTbI}}}{\raMTx{-1}{\!\roMTx{45}{\RCTI}\!\!\raMTx{-3}{\circ}}}]}		\newcommand{\lfgLSsO}[2][?]{{\ifx#1?\lfgLSsN\else\reflectbox{\lfgLSsN}\fi}^{(#2)}}			
	\newcommand{\lfgLSSN}{[\orMTx{\raMTx{-2}{\roMTx{90}{\RCTbI}} }{\raMTx{2}{\roMTx{-45}{\RCTI}}}]}						\newcommand{\lfgLSSO}[2][?]{{\ifx#1?\lfgLSSN\else\reflectbox{\lfgLSSN}\fi}^{(#2)}}			
	\newcommand{\lfgLSSsN}{[\orMTx{\raMTx{-2}{\roMTx{90}{\RCTbI}}}{\raMTx{2}{\roMTx{-45}{\!\RCTI}\!\!\raMTx{3}{\circ}}}]}		\newcommand{\lfgLSSsO}[2][?]{{\ifx#1?\lfgLSSsN\else\reflectbox{\lfgLSSsN}\fi}^{(#2)}}			
	\newcommand{\lfgSHiN}{sh}		\newcommand{\lfgSHiT}[3][?]{\lfgSHiN_{#2}^{(#3)}}

	\newcommand{\lsstSh}[1][?]{\ifx#1?{\mathsf{S}}_{\sh}\else{\mathsf{S}}_{\sh}^{(#1)}\fi}
	\newcommand{\lsstShT}[2][?]{\ifx#1?{\mathsf{S}}_{\sh,#2}\else{\mathsf{S}}_{\sh,#2}^{(#1)}\fi}

	\newcommand{\pgpC}[2][?]{\mathfrak{C}_{#2}}			\newcommand{\gpC}{\pgpC}
	\newcommand{\pgpCs}[2][?]{\mathfrak{C}_{#2}^0}		\newcommand{\gpCs}{\pgpCs}
	\renewcommand{\vlWT}[2][n]{\Fc[#1]{w}{#2}}
	\renewcommand{\vlDP}[2][n]{\Fc[#1]{d}{#2}}
	\newcommand{\pmpDPN}{\vlDP[]{}}						\newcommand{\pmpDPT}[3][n]{\idFc[#1]{\pmpDPN}{#2}{#3}}
																						
																					\newcommand{\mpDPT}{\pmpDPT}
																					
	\newcommand{\pmpWTN}{\vlWT[]{}}			\newcommand{\pmpWTO}[2][?]{\pmpWTN_#2}			\newcommand{\pmpWTT}[3][n]{\idFc[#1]{\pmpWTN}{#2}{#3}}
																					\newcommand{\mpWTO}{\pmpWTO}	\newcommand{\mpWT}{\mpWTO}
																					\newcommand{\mpWTT}{\pmpWTT}
	\renewcommand{\fcZBlettHelp}{\dagger}
	\newcommand{\letB}{\fcZBlettHelp}

	\newcommand{\pempWTN}{\vlWT[]{}}	\newcommand{\pempWTT}[3][?]{\pempWTN_#2^{(#3)}}		\newcommand{\pempWTTh}[4][n]{\ipFc[#1]{\pempWTN}{#2}{(#3)}{#4}}
																					\newcommand{\empWTT}{\pempWTT}	\newcommand{\empWT}{\empWTT}
																					\newcommand{\empWTTh}{\pempWTTh}
	\newcommand{\pmpRH}[2][n]{\Fc[#1]{\rho}{#2}}
	\newcommand{\pmpSN}{\Sig}	\newcommand{\pmpSO}[2][?]{\pmpSN_{#2}}
		\newcommand{\pmpS}{\pmpSO}		\newcommand{\mpS}{\pmpS}
	\newcommand{\pstUT}[1][?]{U_3}		\newcommand{\pstUF}[1][?]{U_4}				
	\newcommand{\pstVF}[1][?]{\ifx?#1V_4\else V_4^#1\fi}		
	\newcommand{\pstWF}[1][?]{\ifx?#1W_4\else W_4^#1\fi}		
	\newcommand{\pstWFO}[2][?]{\ifx?#1W_{4,#2}\else W_{4,#2}^#1\fi}		
	\newcommand{\pstXF}[1][?]{\ifx?#1X_4\else X_4^#1\fi}		
	\newcommand{\pstPO}[2][?]{\mathcal{P}_{#2}}			\newcommand{\stP}{\pstPO}
	\newcommand{\pstPhT}[3][?]{\mathcal{P}_{#2;#3}}		\newcommand{\pstPhTh}[4][?]{\mathcal{P}_{#2;#3}^{(#4)}}	
		\newcommand{\stPh}{\pstPhT}		\newcommand{\stPhTh}{\pstPhTh}
	
			\newcommand{\EqvT}[3][n]{\bkS[#1]{#2}_{#3}}		\newcommand{\Eqv}{\EqvT}
	
	\newcommand{\pfcCP}[3][n]{\idFc[#1]{\tilde{c}}{#2}{#3}}
	\newcommand{\pfcZP}[3][n]{\idFc[#1]{\zeta}{1}{#2;#3}}			\newcommand{\fcZP}{\pfcZP}
	\newcommand{\pfcZPB}[3][n]{\ipFc[#1]{\zeta}{1}{\letB}{#2;#3}}		\newcommand{\fcZPB}{\pfcZPB}
	\newcommand{\pfcZPH}[3][n]{\ipFc[#1]{\zeta}{1}{*}{#2;#3}}		\newcommand{\fcZPH}{\pfcZPH}
	\newcommand{\pfcZPS}[3][n]{\ipFc[#1]{\zeta}{1}{\sh}{#2;#3}}		\newcommand{\fcZPS}{\pfcZPS}

	\newcommand{\pfcPChS}[3][n]{\Fc[n]{\LetHelpFcCh^\sh}{#2;#3}}			\newcommand{\fcPChS}{\pfcPChS}

	\newcommand{\pfcA}[2][n]{\Fc[#1]{A}{#2}}
	\newcommand{\pfcCVNLett}{ct}				\newcommand{\pfcCVN}[1][?]{\ifx ?#1\pfcCVNLett\else\pfcCVNLett_{#1}\fi}	
	\newcommand{\pfcCVoNLett}{\mathit{1}}		\newcommand{\pfcCVoN}[1][?]{\ifx ?#1\pfcCVoNLett\else\pfcCVoNLett_{#1}\fi}	
				
	\newcommand{\pfcCVoAO}[2][?]{\Fc[n]{\pfcCVoN[#1]}{#2}}	
			\newcommand{\pfcCVo}{\pfcCVoN}				\newcommand{\pfcCVoA}{\pfcCVoAO}		
				\newcommand{\fcCVo}{\pfcCVo}				\newcommand{\fcCVoA}{\pfcCVoA}	

	\newlength{\tmpSpaceLen} 
	\newcommand{\tmpSpace}{\hspace{\the\tmpSpaceLen}}
	
	\renewcommand{\gpSym}[2][?]{\mathfrak{S}_{#2}}
	
	\renewcommand{\gpAlt}[2][?]{\mathfrak{A}_{#2}}
	\newcommand{\gpAlts}[2][?]{\mathfrak{A}_{#2}^0}
	\renewcommand{\gpu}[1][?]{\ifx?#1e\else e_{#1}\fi}

\allowdisplaybreaks[4]
	\title{\mainTitle}
	\author{\authorName
			\thanks{\organizationNameFst, \placeAddressFst}
			\mbox{}
			\thanks{\organizationNameScd, \departmentNameScd, \placeAddressScd}
		}
	\date{}

\begin{document}
\maketitle
\renewcommand{\thefootnote}{\fnsymbol{footnote}}
\footnote[0]{e-mail:\,\emailAddressFst} 
\footnote[0]{MSC-class: \MSCname}
\renewcommand{\thefootnote}{\arabic{footnote}}\setcounter{footnote}{0}
\vPack[30]

\begin{abstract}
In the present paper,
	we give identities involving cyclic sums of regularized multiple zeta values of depth less than $5$.
As a corollary,
	we present two kinds of extensions of Hoffman's theorem for symmetric sums of multiple zeta values for this case. 
\end{abstract}
\section{Introduction and statement of results} \label{sectOne}
Multiple zeta values (MZVs) and their generalizations, which are called regularized multiple zeta values (RMZVs), are real numbers that 
	are variations of special values of the Riemann zeta function 
	\envM{
		\fcMZ{1}{s}
	}{
		\SmT{m=1}{\infty} \opF[s]{1}{m^s}
	}
	with integer arguments.
It is known that these values satisfy a great many relations over $\setQ$, including, for example,
	extended harmonic and shuffle relations, Drinfel'd associator relations, and Kawashima's relations (e.g., see \cite{Drinfeld90,IKZ06,Kawashima09}).
New classes of relations are being studied,
	but their exact structure is not yet fully understood.
In the present paper,
	we give new identities involving cyclic sums of RMZVs of depth less than $5$. 
As a corollary,
	we offer two kinds of extensions of Hoffman's theorem \cite[\refThmA{2.2}]{Hoffman92} for symmetric sums of MZVs for this case.

We will begin by introducing the notation and terminology that will be used to state our results.
An MZV is a convergent series defined by
	\envHLineCmDef
	{
		\fcMZ{n}{\lbfL_n}
	}
	{
		\Sm{m_1>\cdots>m_n>0} \opF{1}{ \pw{m_1}{l_1}\cdots\pw{m_n}{l_n} }
	} 
	where $\lbfL_n=(l_1,\ldots,l_n)$ is an (ordered) index set of positive integers with $l_1\geq2$.
In other words,
	MZVs are images under the real-valued function $\fcMZO{n}$ with the domain $\SetT{(l_1,\ldots,l_n)\in\setN^n}{l_1\geq2}$.
We call  
	\envM{\mpWTT{n}{\lbfL_n}}{l_1+\cdots+l_n} 
	the weight,
	and 
	\envM{\mpDPT{n}{\lbfL_n}}{n} 
	the depth.
Ihara, Kaneko, and Zagier \cite{IKZ06} extended MZVs to two types of RMZV 
	(harmonic and shuffle)
	with two different renormalization procedures for divergent series $\fcMZ{n}{\lbfL_n}$ of $l_1=1$.
The former and latter  types
	are denoted by $\fcMZH{n}{\lbfL_n}$ and $\fcMZS{n}{\lbfL_n}$, 
	and they
	inherit the harmonic and shuffle relation structures,
	respectively.
The following are a few examples of these values:
	$\fcMZH{1}{1}=\fcMZS{1}{1}=\fcMZS{2}{1,1}=0$ 
	and
	$\fcMZH{2}{1,1}=-\opF[s]{\fcMZ{1}{2}}{2}\neq0$.
In other words,
	RMZVs $\fcMZH{n}{\lbfL_n}$ and $\fcMZS{n}{\lbfL_n}$ are images under two different extension functions  of $\fcMZO{n}$ to the domain $\setN^n$.

Let  $\gpSym{n}$ denote the symmetric group of degree $n$,
	and let $\gpu=\gpu[n]$ denote its unit element.
Let $\gpC{3}$ and $\gpC{4}$ 
	be the cyclic subgroups in $\gpSym{3}$ and $\gpSym{4}$ given by
	\envMTh{
		\gpC{3}
	}{
		\Gp{(123)}
	}{
		\SetO{\gpu, (123), (132)}
	}
	and
	\envMThCm{
		\gpC{4}
	}{
		\Gp{(1234)}
	}{
		\SetO{\gpu, (1234), (13)(24), (1432)}
	}
	respectively.
We set $\gpC{2}=\Gp{(12)}$ (or $\gpC{2}=\gpSym{2}$) for convenience.
The group ring $\setZ[\gpSym{n}]$ of $\gpSym{n}$ over $\setZ$ acts on a function $f$ of $n$ variables in a natural way by
	\envHLineCmDef
	{
		\racFA{f}{\Gam}{x_1,\ldots,x_n}
	}
	{
		\SmN[t] a_i 
		\Fc{f}{x_{\sig^{-1}(1)},\ldots,x_{\sig^{-1}(n)}}
	}
	where $\Gam=\SmN[t] a_i\sig_i\in\setZ[\gpSym{n}]$.
This is a right action,
	that is,
	$\racF{f}{(\Gam_1\Gam_2)} =	\racF{\racFrr{f}{\Gam_1}}{\Gam_2}$.
For a subset $H$ in $\gpSym{n}$,
	we define the sum of all elements in $H$ by
	\envHLineDefPd
	{
		\mpS{H}
	}
	{
		\Sm[t]{\sig\in H} \sig
	\lnP{\in}
		\setZ[\gpSym{n}]
	}
That is,
	\envMO
	{
		\racFA{f}{ \mpS{H} }{x_1,\ldots,x_n}
	}
	is 
	\envMOPd
	{
		\Sm{\sig\in H} \Fc{f}{x_{\sig^{-1}(1)},\ldots,x_{\sig^{-1}(n)}}
	}
In particular,
	if $H$ is a group,
	it is 
	\envMO{
		\Sm{\sig\in H} \Fc{f}{x_{\sig(1)},\ldots,x_{\sig(n)}}
	} 
	because $H=H^{-1}$.
For real-valued functions $f_{n_1},\ldots,f_{n_j}$ such that each $f_{n_i}$ has  $n_i$ variables,
	we define the function $f_{n_1}\otimes\cdots\otimes f_{n_j}$ of $n=n_1+\cdots+n_j$ variables by
	\envLineDefPd
	{
		\Fc{f_{n_1}\otimes f_{n_2}\otimes\cdots\otimes f_{n_j}}{x_1,\ldots,x_n} 
	}
	{
		\Fc{f_{n_1}}{x_1,\ldots,x_{n_1}}
		\Fc{f_{n_2}}{x_{n_1+1},\ldots,x_{n_1+n_2}}
		\cdots
		\Fc{f_{n_j}}{x_{n_1+n_2+\cdots+n_{j-1}+1},\ldots,x_n}
	}
For example,
	\envM
	{
		\Fc{\fcMZO{1}\otimes\fcMZO{1}}{\lbfL_2}
	}{
		\fcMZ{1}{l_1}\fcMZ{1}{l_2}
	}
	and
	\envMPd{
		\Fc{\fcMZO{2}\otimes\fcMZO{1}}{\lbfL_3}
	}{
		\fcMZ{2}{l_1,l_2}\fcMZ{1}{l_3}
	}
Note that the operation $\otimes$ is not commutative.
We define the characteristic functions $\fcMChH[]{n}{}$ and $\fcMChS[]{n}{}$ of the set $\setN^n$ by
	\envHLineCFLaa
	{
		\fcMChH{n}{\lbfL_n}
	}
	{
		\fcCVoA[n]{\lbfL_n}
	\lnP{=}
		1
	}
	{
		\fcMChS{n}{\lbfL_n}
	}
	{
		\envCaseTCm{
			0			&	(l_1=\cdots=l_n=1)
		}{
			1			&	(otherwise)
		}
	}
	respectively.

Our results are stated as follows. 

\begin{theorem}\label{1_Thm1}
Let $\lbfL_n=(l_1,\ldots,l_n)$ be an index set in $\setN^n$,
	and let $L_n=\mpWTT{n}{\lbfL_n}$ be its weight. 
Then we have the following identities for RMZVs $\fcMZH{n}{\lbfL_n}$ and $\fcMZS{n}{\lbfL_n}$ of $n=2,3$, and $4$:
	\envHLineCSCmNme
	{\label{1_Thm1_EqDpt2}
		\racFA{\fcMZBO{2}}{ \mpS{\gpC{2}} }{\lbfL_2}
	}
	{
		\Fc{\fcMZBO{1}\otimes\fcMZBO{1}}{\lbfL_2} 
		-
		\fcMChB{2}{\lbfL_2}\fcMZ{1}{L_2}	
	}
	{\label{1_Thm1_EqDpt3}\rule{0cm}{20pt}
		\racFA{\fcMZBO{3}}{\mpS{\gpC{3}}}{\lbfL_3}
	}
	{
		-
		\Fc{\fcMZBO{1}\otimes\fcMZBO{1}\otimes\fcMZBO{1}}{\lbfL_3} 
		+
		\racFA{\fcMZBO{2}\otimes\fcMZBO{1}}{\mpS{\gpC{3}}}{\lbfL_3} 
		+
		\fcMChB{3}{\lbfL_3}\fcMZ{1}{L_3}	
	}
	{\label{1_Thm1_EqDpt4}\rule{0cm}{20pt}
		\racFA{\fcMZBO{4}}{\mpS{\gpC{4}}}{\lbfL_4}
	}
	{
		\Fc{\fcMZBO{1}\otimes\fcMZBO{1}\otimes\fcMZBO{1}\otimes\fcMZBO{1}}{\lbfL_4} 
		-
		\racFA{\fcMZBO{2}\otimes\fcMZBO{1}\otimes\fcMZBO{1}}{\mpS{\gpC{4}}}{\lbfL_4}
		\lnAH
		+
		\racFA{\fcMZBO{2}\otimes\fcMZBO{2}}{\mpS{\gpCs{4}}}{\lbfL_4}
		+
		\racFA{\fcMZBO{3}\otimes\fcMZBO{1}}{\mpS{\gpC{4}}}{\lbfL_4}
		-
		\fcMChB{4}{\lbfL_4}\fcMZ{1}{L_4}	
	}
	where $\letB\in\Set{*,\sh}$,
	and $\gpCs{4}$ in \refEq{1_Thm1_EqDpt4} is the subset $\Set{\gpu, (1234)}$ of $\gpC{4}$.
\end{theorem}

We note that \refEq{1_Thm1_EqDpt2} can be easily obtained from the harmonic relations 
	$\fcMZH{1}{l_1}\fcMZH{1}{l_2}=\fcMZH{2}{l_1,l_2}+\fcMZH{2}{l_2,l_1}+\fcMZH{1}{l_1+l_2}$ 
	for RMZVs of harmonic type of depth $2$;
	thus our main results are \refEq{1_Thm1_EqDpt3} and \refEq{1_Thm1_EqDpt4}
	(see \refSec{sectFive} for their straightforward expressions).

We now recall Hoffman's theorem.
For any partition $\Pi=\bkB{P_1,\ldots,P_j}$ of $\SetO{1,\ldots,n}$, 
	we define an integer $\pfcCP{n}{\Pi}$ and a real number $\fcZP{\lbfL_n}{\Pi}$ by
	\envHLineCFLaaCmDef[1_PL_DefNmPartition]
	{
		\pfcCP{n}{\Pi}
	}
	{
		\mo^{n-j} \PdT{i=1}{j} (\vlNS{P_i}-1)! 
	}
	{
		\fcZP{\lbfL_n}{\Pi}
	}
	{
		\PdT{i=1}{j} \fcMZ[G]{1}{ \Sm{p\in P_i} l_p }
	}
	respectively,
	where $\vlNS{P}$ is the number of elements of the set $P$,
	and
	$\fcMZ{1}{1}$ stands for $0\;(=\fcMZB{1}{1})$. %
(The value $\fcMZ{1}{1}$ is not necessary for Hoffman's theorem,
	but we will use it in \refCor{1_Cor1}, 
	below.)
Hoffman \cite[\refThmA{2.2}]{Hoffman92} proved 
	the following identity involving symmetric sums of MZVs of depth $n$:
	\envHLineCm[1_PL_EqSymMZV]
	{
		\racFA{\fcMZO{n}}{ \mpS{\gpSym{n}} }{\lbfL_n}
	}
	{
		\Sm{\mathrm{partitions}\;\Pi\;\mathrm{of}\;\SetO{1,\ldots,n}} \pfcCP{n}{\Pi} \fcZP{\lbfL_n}{\Pi}
	}
	where every component of $\lbfL_n$ is at least $2$.
Note 
	that \refEq{1_PL_EqSymMZV} was proved under the weaker condition 
	that 
	$\lbfL_n$ is an index set $\lbfS_n$ of real numbers greater than $1$.

We require some notation to state \refCor{1_Cor1} of \refThm{1_Thm1}.
For any subset $P=\Set{p_1,\ldots, p_m}$ of $\Set{1,\ldots,n}$,
	let $\fcPChS{\lbfL_n}{P}$ be 0 if $l_p=1$ for all $p\in P$, and $1$ otherwise.
For example,
	\envMCm{
		\fcPChS{(2,1,1)}{\Set{2,3}}
	}{
		0
	}
	and
	\envMThPd{
		\fcPChS{(2,1,1)}{\Set{1,2}}
	}{
		\fcPChS{(2,1,1)}{\Set{1,3}}
	}{
		1
	}
In other words,
	\envMPd
	{
		\fcPChS{\lbfL_n}{P}
	}
	{
		\fcMChS{m}{l_{p_1},\ldots,l_{p_m}}
	}
We then define the real numbers $\fcZPH{\lbfL_n}{\Pi}$ and $\fcZPS{\lbfL_n}{\Pi}$
	by
	\envHLineCFLaaCmDef
	{
		\fcZPH{\lbfL_n}{\Pi}
	}
	{
		\fcZP{\lbfL_n}{\Pi}
	}
	{
		\fcZPS{\lbfL_n}{\Pi}
	}
	{
		\PdT{i=1}{j} \fcPChS{\lbfL_n}{P_i} \fcMZ[G]{1}{ \Sm{p\in P_i} l_p }
	}
	respectively.
	
\refCor{1_Cor1} is stated as follows.

\begin{corollary}\label{1_Cor1}
Let $\lbfL_n$ and $\letB$ be as in \refThm{1_Thm1}.
For any $n\in\Set{2,3,4}$,
	\envHLinePd[1_Cor_EqSymRMZV]
	{
		\racFA{\fcMZBO{n}}{ \mpS{\gpSym{n}} }{\lbfL_n}
	}
	{
		\Sm{\mathrm{partitions}\;\Pi\;\mathrm{of}\;\SetO{1,\ldots,n}} \pfcCP{n}{\Pi} \fcZPB{\lbfL_n}{\Pi}
	}
	
\end{corollary}

For a depth less than $5$, we can easily deduce \refEq{1_PL_EqSymMZV} from \refEq{1_Cor_EqSymRMZV},
	since 
	$\fcMZB{n}{\lbfL_n}=\fcMZ{n}{\lbfL_n}$
	and
	$\fcZPB{\lbfL_n}{\Pi}=\fcZP{\lbfL_n}{\Pi}$ when $l_1,\ldots,l_n>1$.
This means that we obtain two kinds of extensions of \refEq{1_PL_EqSymMZV} for such a depth.

We now briefly explain how \refThm{1_Thm1} and \refCor{1_Cor1} can be proved.
We first prove the identities in \refThm{1_Thm1} for $\letB=*$ by using harmonic relations of RMZVs $\fcMZH{n}{\lbfL_n}$.
We then derive the identities for $\letB=\sh$ from those for $\letB=*$ and
	relations over $\setQ$ between RMZVs $\fcMZH{n}{\lbfL_n}$ and $\fcMZS{n}{\lbfL_n}$ given in \cite[\refThmA{1}]{IKZ06} (which we call renormalization relations).
We prove \refCor{1_Cor1} by using the identities in \refThm{1_Thm1} and the fact that $\gpC{n}$ is a subgroup of $\gpSym{n}$, 
	i.e.,
	$\racFA{\fcMZBO{n}}{\mpS{\gpC{n}}}{\lbfL_n}$ is a partial sum of $\racFA{\fcMZBO{n}}{\mpS{\gpSym{n}}}{\lbfL_n}$.

It is worth noting that \refThm{1_Thm1} gives the following property,
	which is an analog of the parity property \cite{BG96,Euler1776,IKZ06,Tsumura04};
	any cyclic sum of RMZVs of depth less than $5$,
	or
	\envHLine[1_PL_EqCycSumDetail]
	{
		\racFA{\fcMZBO{n}}{ \mpS{\gpC{n}} }{\lbfL_n}
	}
	{
		\SmT{j=1}{n} \fcMZB{n}{l_j,\dots,l_n,l_1,\ldots,l_{j-1}}
	}
	for $n=2,3,4$ and $\letB\in\Set{*,\sh}$,
	is a linear combination of the Riemann zeta value $\fcMZ{1}{l_1+\cdots+l_n}$ and products of RMZVs of smaller depth and weight.
It appears that the existence of
	such a property for depths greater than $4$ is an open problem.
(The case of symmetric sums of general depth easily follows from Hoffman's identity \refEq{1_PL_EqSymMZV};
there is a stronger property from \refEq{1_PL_EqSymMZV}, such that 
	any symmetric sum can be written in terms of only Riemann zeta values.)
It is also worth noting that 
	Hoffman and Ohno \cite{HO03} studied a class of relations involving
	\envPLineCm
	{
		\SmT{j=1}{n} \fcMZ{n}{l_j+1,l_{j+1},\dots,l_n,l_1,\ldots,l_{j-1}}
	}
	whose form is quite similar to \refEq{1_PL_EqCycSumDetail},
	but
	the first indices differ.

The paper is organized as follows.
In \refSec{sectTwo}, we review some facts of RMZVs by referring to \cite{Hoffman97,IKZ06}.
\refSec[s]{sectThree} and \ref{sectFour} each have two subsections.
\refSec[s]{sectThreeOne} and \ref{sectThreeTwo} are devoted to calculating 
	harmonic relations for RMZVs $\fcMZH{n}{\lbfL_n}$ and renormalization relations between RMZVs $\fcMZH{n}{\lbfL_n}$ and $\fcMZS{n}{\lbfL_n}$,
	respectively,
	for depths less than $5$.
We then prove \refThm{1_Thm1} in \refSec{sectFourOne}
	and \refCor{1_Cor1} in  \refSec{sectFourTwo}.
We give some examples of \refThm{1_Thm1} and \refCor{1_Cor1} in \refSec{sectFive}.

\begin{remarkN}\label{1_Rem1}
\bfTx{(i)}
Although the ideas of the proofs are the same,
	the computational complexity of proving \refEq{1_Thm1_EqDpt4} is much greater than that required to prove \refEq{1_Thm1_EqDpt2} and \refEq{1_Thm1_EqDpt3}. 
We recommend that, on first reading, those readers who are interested only in the ideas 
	skip over the statements relating to the proof of \refEq{1_Thm1_EqDpt4} (or statements in the case of depth $4$).
\mbox{}\bfTx{(ii)}
The present paper is an expansion of \refSectA{2.1} in \cite{Machide12Ax}.
The remainder of the results of \cite{Machide12Ax} will be amplified in a forthcoming paper \cite{Machide14Ax}.
\end{remarkN}

\section{Preparation}\label{sectTwo}
Let $\alHF=\setQ\bkA{x,y}$ be the noncommutative polynomial algebra over $\setQ$ in two indeterminates $x$ and $y$,
	and let $\alHN$ and $\alHR$ be its subalgebras $\setQ+x\alHF y$ and $\setQ+\alHF y$, respectively.
These satisfy the inclusion relations $\alHN\subset\alHR\subset\alHF$.
Let $z_l$ denote $x^{l-1}y$ for any integer $l\geq1$.
Every word $w=w_0y$ in the set $\Set{x,y}$ with terminal letter $y$
	is expressed as $w=z_{l_1}\cdots z_{l_n}$ uniquely,
	and so $\alHR$ is the free algebra generated by $z_l\;(l=1,2,3,\ldots)$.
We define the harmonic product $*$ on $\alHR$ inductively by
	\envPLineTNme
	{\label{2_PL_DefHarProdUnit}&
		1 * w \lnP{=} w * 1 \lnP{=} w
	,&}
	{\label{2_PL_DefHarProdMain}&
		z_kw_1 * z_lw_2 \lnP{=} z_k(w_1*z_lw_2) + z_l(z_kw_1*w_2) + z_{k+l}(w_1*w_2)
	,&}
	for any integers $k,l\geq1$ and words $w, w_1, w_2\in\alHR$,
	and then extend it by $\setQ$-bilinearity.
This product gives the subalgebras $\alHN$ and $\alHR$ structures of commutative $\setQ$-algebras \cite{Hoffman97},	
	which we denote by $\alHNh$ and $\alHRh$, respectively;
	note that $\alHNh$ is a subalgebra of $\alHRh$. 
In a similar way,
	we can define the shuffle product $\sh$ on $\alHR$ 
	and the commutative $\setQ$-algebras $\alHNs$ and $\alHRs$ (see  \cite{IKZ06,Reutenauer93} for details).

Let $\mpEVN[]{}:\alHN\to\setR$ be the $\setQ$-linear map (evaluation map) given by
	\envHLinePd[2_PL_DefMapEV]
	{
		\mpEVN{z_{l_1}\cdots z_{l_n}}
	}
	{
		\fcMZ{n}{\lbfL_n}
		\qquad
		(z_{l_1}\cdots z_{l_n}\in\alHN)
	}
We know from \cite{Hoffman97} that 
	$\mpEVN[]{}$ is homomorphic on both products $*$ and $\sh$, 
	i.e.,
	\envOTLine
	{
		\mpEVN{w_1*w_2}
	}
	{
		\mpEVN{w_1 \shS w_2}
	}
	{
		\mpEVN{w_1} \mpEVN{w_2}
	}
	for $w_1,w_2\in\alHN$. 
Let $\setR[T]$ be the polynomial ring in a single indeterminate with real coefficients.
Through the isomorphisms $\alHRh\simeq\alHNh[y]$ and $\alHRs\simeq\alHNs[y]$, which were proved in \cite{Hoffman97} and \cite{Reutenauer93}, respectively,
	\etalTx{Ihara} \cite[\refPropA{1}]{IKZ06} considered 
	the algebra homomorphisms $\mpEVH[]{} : \alHRh \to \setR[T]$ and $\mpEVS[]{} : \alHRs \to \setR[T]$,
	respectively,
	which are uniquely characterized by the property that they extend the evaluation map $\mpEVN[]{}$ and send $y$ to $T$.		
For any word $w=z_{l_1}\cdots z_{l_n}\in\alHR$,
	we denote by $\mpEVHi{\lbfL_n}{T}$ and $\mpEVSi{\lbfL_n}{T}$ 
	the images under the maps $\mpEVH[]{}$ and $\mpEVS[]{}$, respectively, of the word $w$,
	that is,
	\envHLineCFLaaPd[2_PL_DefEVimage]
	{
		\mpEVHi{l_1,\ldots,l_n}{T}
	}
	{
		\mpEVH{z_{l_1}\cdots z_{l_n}}
	}
	{
		\mpEVSi{l_1,\ldots,l_n}{T}
	}
	{
		\mpEVS{z_{l_1}\cdots z_{l_n}}
	}
(The notation $\mpEVHi{\lbfL_n}{T}$ and $\mpEVSi{\lbfL_n}{T}$ will be used 
	when we focus on the variable $T$ and the corresponding index set $\lbfL_n$ of the word $z_{l_1}\cdots z_{l_n}$.)
Then the RMZVs $\fcZH{\lbfL_n}$ and $\fcZS{\lbfL_n}$ of the harmonic and shuffle types are defined as
	\envHLineCFLaaCmDef[2_PL_DefRMZV]
	{
		\fcZH{l_1,\ldots,l_n}
	}
	{
		\mpEVHi{l_1,\ldots,l_n}{0}
	}
	{
		\fcZS{l_1,\ldots,l_n}
	}
	{
		\mpEVSi{l_1,\ldots,l_n}{0}
	}
	respectively.
Obviously,
	$\fcMZH{n}{\lbfL_n}=\fcMZS{n}{\lbfL_n}=\fcMZ{n}{\lbfL_n}$ if $l_1>1$. 
We have
	\envHLine
	{
		\mpEVH{z_{k_1}\cdots z_{k_m}*z_{l_1}\cdots z_{l_n}}
	}
	{
		\mpEVH{z_{k_1}\cdots z_{k_m}} \mpEVH{z_{l_1}\cdots z_{l_n}}
	}
	for index sets $(k_1,\ldots,k_m)$ and $(l_1,\ldots,l_n)$,
	since $\mpEVH[]{}$ is homomorphic,
	and so
	we see from the first equations of \refEq{2_PL_DefEVimage} and \refEq{2_PL_DefRMZV} that the
	RMZVs $\fcMZH{n}{\lbfL_n}$ satisfy the harmonic relations.
In \refSect{sectThreeOne}, we will calculate these relations in detail for depths less than $5$.
(We can also see that the RMZVs $\fcMZS{n}{\lbfL_n}$ satisfy the shuffle relations since $\mpEVS[]{}$ is homomorphic,
	but 
	we will not discuss this in the present paper.)

Let $\pfcA{u}=\SmT{k=0}{\infty} \gam_k u^k$ be the Taylor expansion of $\exN{\gam u}\fcG{1+u}$ near $u=0$,
	where $\gam$ is Euler's constant and $\fcG{x}$ is the gamma function.
The renormalization map $\pmpRH[]{}:\setR[T]\to\setR[T]$
	is an $\setR$-linear map defined by
	\envHLinePd[2_PL_DefMapRho]
	{
		\pmpRH{\exN{Tu}}
	}
	{
		\pfcA{u} \exN{Tu}
	}
That is,
	images $\pmpRH{T^m}$ are determined by comparing the coefficients of $u^m$ on both sides of \refEq{2_PL_DefMapRho},
	and
	expressed as
	\envHLinePd[2_PL_DefMapRhoDetail]
	{\qquad
		\pmpRH{T^m}
	}
	{
		m! \SmT{i=0}{m} \gam_i \opF{T^{m-i}}{(m-i)!}
		\qquad
		(m=0,1,2,\ldots)
	}
Then the renormalization formula proved by \etalTx{Ihara} \cite[\refThmA{1}]{IKZ06} is
	\envHLinePd[3.2_PL_EqFundRelRegIKZ]
	{
		\pmpRH{\mpEVHi{\lbfL_n}{T}}
	}
	{
		\mpEVSi{\lbfL_n}{T}
	}
Combining \refEq{2_PL_DefRMZV}  and \refEq{3.2_PL_EqFundRelRegIKZ} with $T=0$,
	we can obtain 
	relations between RMZVs $\fcMZH{n}{\lbfL_n}$ and $\fcMZS{n}{\lbfL_n}$,
	or
	renormalization relations.
In \refSect{sectThreeTwo}, we will calculate these relations in detail for depths less than $5$.

\section{Relations}\label{sectThree}	
\subsection{Harmonic relations}\label{sectThreeOne}	
We begin by defining the notation that we will use
	to state the harmonic relations of RMZVs $\fcMZH{n}{\lbfL_n}$ of depth less than $5$
	in terms of real-valued functions.
	
We first define analogs of the weight map $\mpWT{n}:\setN^{n}\to\setN$ of depth $n$.
For positive integers $n_1,\ldots,n_j$, and $n$ with $n_1+\cdots+n_j=n$,
	we define the map $\empWT{n}{n_1,\ldots,n_j}$ from $\setN^{n}$ to $\setN^j$ by
	\envHLineDefPd[3.1_PL_DefFuncW]
		{
			\empWT{n}{n_1,\ldots,n_j}
		}
		{
			\mpWT{{n_1}}\otimes\cdots\otimes\mpWT{{n_j}}
		}
For example,
	$\empWTTh{3}{2,1}{l_1,l_2,l_3}=(l_1+l_2,l_3)$ and $\empWTTh{4}{1,2,1}{l_1,l_2,l_3,l_4}=(l_1,l_2+l_3,l_4)$.
Obviously,
	\envMPd{
		\empWT{n}{n}
	}{
		\mpWT{n}
	}
We define a subset $\pstUT$ in $\gpSym{3}$ as
	\envHLineCm[3.1_PL_DefGpDpt3]
	{
		\pstUT
	}
	{
		\Set{\gpu[3],(23),(123)}
	}
	and subsets $\pstUF$, $\pstVF[0]$, $\pstVF$, $\pstWF[0]$, $\pstWF[1]$, $\pstWF$, and $\pstXF$ in $\gpSym{4}$ as
	\envPLine[????]{\HLineCTCmNme[p]
	{\label{3.1_PL_DefGpDpt4U}\rule{0cm}{12pt}
		\pstUF
	}
	{
		\Set{\gpu[4],(34),(234),(1234)}
	}\HLineCFCm[p]
	{\label{3.1_PL_DefGpDpt4V}\rule{0cm}{20pt}
		\pstVF[0]
	}
	{
		\Set{(23),(1243)}
	}
	{
		\pstVF
	}
	{
		\Set{\gpu[4],(13)(24),(123),(243)}
		\cup
		\pstVF[0]
	}\HLineCFCm[p]
	{\label{3.1_PL_DefGpDpt4W}\rule{0cm}{20pt}
		\pstWF[0]
	}
	{
		\Set{(23),(24)}
	}
	{
		\pstWF[1]
	}
	{
		\Set{(34),(1234),(1243),(1324)}
		\cup
		\pstWF[0]
	}\HLineCFPd
	{
		\pstWF
	}
	{
		\Set{\gpu[4],(13)(24),(123),(124),(234),(243)}
		\cup
		\pstWF[1]
	}
	{\label{3.1_PL_DefGpDpt4X}\rule{0cm}{20pt}
		\pstXF
	}
	{
		\Set{(14), (23)} \cup \gpC{4}
	}
	}
We have the inclusion relations
	$\pstVF[0]\subset\pstVF$ and $\pstWF[0]\subset\pstWF[1]\subset\pstWF$.
We denote by $\gpAlt{3}$ and $\gpAlt{4}$ the alternating groups of degree $3$ and $4$,
	respectively. 
Note that 
	$\gpAlt{3}=\gpC{3}$. 
The operator $\circ$ of the function composition satisfies the distributive law,
	i.e.,
	\envMCm
	{
		(\Sm{i} f_i)\circ(\Sm{j} g_j)
	}
	{
		\Sm{i,j} f_i \circ g_j
	}
	for real-valued functions $f_i$ with domain $\setN^n$ and vector-valued functions $g_j$ whose images are included in $\setN^n$.
We use the symbol $\racF{f\circ g}{\sig}$ for both $\racFr{f\circ g}{\sig}$ and $f\circ(\racF{g}{\sig})$ for any $\sig\in\gpSym{n}$
	since, by definition, they are equal.

\begin{remark}\label{3.1_Rem1}
For integers $l,n$ with $1\leq j \leq n-1$,
	let $\lfgSHiT{j}{n}$ be the shuffle elements given in \cite{IKZ06},
	which 
	are elements in $\setQ[\gpSym{n}]$ 
	and defined as
	\envHLineDefPd
	{
		\lfgSHiT{j}{n}
	}
	{
		\tpSm{\sig\in\gpSym{n}}{\sig(1)<\cdots<\sig(j) \atop \sig(j+1)<\cdots<\sig(n)} \sig
	}
The elements $\mpS{\pstUT}$,  $\mpS{\pstUF}$, and $\mpS{\pstVF}$ are equal to $\lfgSHiT{2}{3}$, $\lfgSHiT{3}{4}$, and $\lfgSHiT{2}{4}$,
	respectively.
The element $\mpS{\pstWF}$ cannot be written in terms of only a shuffle element,
	but it is equal to 
	\envM
	{
		\mpS{\pstVF} \mpS{\Gp{(34)}}
	}
	{
		\lfgSHiT{2}{4} \mpS{\Gp{(34)}},
	}
	as we will see in \refEq{3.1_Lem5i_Eq3}, below.
\end{remark}

The harmonic relations we desire are listed below. For brevity,  
	we will denote by
	${\fcMZHO{n}}^{\otimes m}$ the function $\usbTx{m}{\fcMZHO{n}\otimes\cdots\otimes \fcMZHO{n}}$ of $mn$ variables.
\begin{proposition}[Case of depth $2$]\label{3.1_Prop1}
We have  
	\envHLinePd[3.1_Prop1_Eq]
	{
		{\fcMZHO{1}}^{\otimes2}
	}
	{
		\racF{\fcMZHO{2}}{ \mpS{\gpC{2}} } + \fcMZO{1} \circ \mpWT{2}
	}
\end{proposition}
\begin{proposition}[Case of depth $3$]\label{3.1_Prop2}
We have 
	\envHLineCFNmePd
	{\label{3.1_Prop2_Eq1}
		\fcMZHO{2}\otimes\fcMZHO{1}
	}
	{
		\racF{\fcMZHO{3}}{ \mpS{\pstUT} } + \fcMZHO{2} \circ  ( \racF{ \empWT{3}{2,1} }{(123)}  + \empWT{3}{1,2} )
	}
	{\label{3.1_Prop2_Eq2}\rule{0cm}{20pt}
		{\fcMZHO{1}}^{\otimes3}
	}
	{
		\racF{\fcMZHO{3}}{ \mpS{\gpSym{3}} } 
		+ 
		\racF{ \fcMZHO{2} \circ (\empWT{3}{2,1}+\empWT{3}{1,2}) }{ \mpS{\gpC{3}} } 
		+
		\fcMZO{1} \circ \mpWT{3}
	}
\end{proposition}
\begin{proposition}[Case of depth $4$]\label{3.1_Prop3}
We have 
	\envHLineCECmNme
	{\label{3.1_Prop3_Eq1}\hspace{-20pt}\rule{0cm}{15pt}
		\fcMZHO{3}\otimes\fcMZHO{1}
	}
	{
		\racF{\fcMZHO{4}}{ \mpS{\pstUF} } + \fcMZHO{3} \circ \bkR[b]{ \racF{(\empWT{4}{2,1,1}+ \empWT{4}{1,2,1}) }{ (234) } + \empWT{4}{1,1,2} }
	}
	{\label{3.1_Prop3_Eq2}\hspace{-20pt}\rule{0cm}{20pt}
		{\fcMZHO{2}}^{\otimes2}
	}
	{
		\racF{\fcMZHO{4}}{ \mpS{\pstVF} } 
		+
		\racF{ \fcMZHO{3} \circ (\empWT{4}{2,1,1} +  \empWT{4}{1,2,1} + \empWT{4}{1,1,2}) }{ \mpS{\pstVF[0]} }
		+
		\racF{ \fcMZHO{2} \circ \empWT{4}{2,2} }{ (23) }
	}
	{\label{3.1_Prop3_Eq3}\hspace{-20pt}\rule{0cm}{20pt}
		\fcMZHO{2}\otimes{\fcMZHO{1}}^{\otimes2}
	}
	{
		\racF{\fcMZHO{4}}{ \mpS{\pstWF} } 
		\lnAH
		+
		\fcMZHO{3} \circ \bkR[b]{ 
			\racF{ \empWT{4}{2,1,1} }{ \mpS{\pstWFO[1]{(34)}} } 
			+ 
			\racF{ \empWT{4}{1,2,1} }{ \mpS{\pstWFO[1]{(1234)}} } 
			+ 
			\racF{ \empWT{4}{1,1,2} }{ \mpS{\pstWFO[1]{(1324)}} } 
		}
		\lnAH[]
		+
		\fcMZHO{2} \circ  \bkR[b]{ \racF{ \empWT{4}{2,2} }{ \mpS{\pstWF[0]} } + \racF{ \empWT{4}{3,1} }{ (24) } + \empWT{4}{1,3} }
		\nonumber
	}
	{\label{3.1_Prop3_Eq4}\hspace{-20pt}\rule{0cm}{20pt}
		{\fcMZHO{1}}^{\otimes4}
	}
	{
		\racF{\fcMZHO{4}}{ \mpS{\gpSym{4}} } 
		+
		\racF{ \fcMZHO{3} \circ (\empWT{4}{2,1,1} + \empWT{4}{1,2,1} + \empWT{4}{1,1,2}) }{ \mpS{\gpAlt{4}} }
		\lnAH
		+
		\fcMZHO{2} \circ \bkR[b]{ 
			\racF{ \empWT{4}{2,2} }{ \mpS{\pstXF} } + \racF{ (\empWT{4}{3,1} + \empWT{4}{1,3}) }{ \mpS{\gpC{4}} }
		}
		+
		\fcMZO{1} \circ \mpWT{4}
	}	
	where 
	$\pstWFO[1]{\sig}$ $(\sig\in\Set{(34),(1234),(1324)})$ in \refEq{3.1_Prop3_Eq3} mean the subsets $\pstWF[1]\setminus\Set{\sig}$. 
\end{proposition}

We will show \refLem[s]{3.1_Lem1}, \ref{3.1_Lem2}, and \ref{3.1_Lem3} 
	to prove \refProp[s]{3.1_Prop1}, \ref{3.1_Prop2}, and \ref{3.1_Prop3},
	respectively.
	These calculate the harmonic products of the generators $z_l$ of $\alHRh$ for the corresponding depths.
\begin{lemma}[Case of depth $2$]\label{3.1_Lem1}
For positive integers $l_1,l_2$,
	we have
	\envHLinePd[3.1_Lem1_Eq]
	{
		z_{l_1}*z_{l_2}
	}
	{
		z_{l_1}z_{l_2} + z_{l_2}z_{l_1} + z_{l_1+l_2}
	}
\end{lemma}
\begin{lemma}[Case of depth $3$]\label{3.1_Lem2}
For positive integers $l_1,l_2,l_3$,
	we have
	\envHLineCFNmePd
	{\label{3.1_Lem2_Eq1}
		z_{l_1}z_{l_2}*z_{l_3}
	}
	{
		z_{l_1}z_{l_2}z_{l_3} + z_{l_1}z_{l_3}z_{l_2} + z_{l_3}z_{l_1}z_{l_2} + z_{l_1+l_3}z_{l_2} + z_{l_1}z_{l_2+l_3}
	}
	{\label{3.1_Lem2_Eq2}\rule{0cm}{15pt}
		z_{l_1}*z_{l_2}*z_{l_3}
	}
	{
		z_{l_1}z_{l_2}z_{l_3} + z_{l_1}z_{l_3}z_{l_2} + z_{l_2}z_{l_1}z_{l_3} 
		+ 
		z_{l_2}z_{l_3}z_{l_1} + z_{l_3}z_{l_1}z_{l_2} + z_{l_3}z_{l_2}z_{l_1} 
		\lnAH
		+ 
		z_{l_1+l_2}z_{l_3} + z_{l_1+l_3}z_{l_2} + z_{l_2+l_3}z_{l_1}
		\lnAH[]
		+
		z_{l_1}z_{l_2+l_3} + z_{l_2}z_{l_1+l_3} + z_{l_3}z_{l_1+l_2} 
		+ 
		z_{l_1+l_2+l_3}
		\nonumber
	}
\end{lemma}
\begin{lemma}[Case of depth $4$]\label{3.1_Lem3}
For positive integers $l_1,l_2,l_3,l_4$,
	we have
	\envHLineCENmePd
	{\label{3.1_Lem3_Eq1}
		z_{l_1}z_{l_2}z_{l_3}*z_{l_4}
	}
	{
		z_{l_1}z_{l_2}z_{l_3}z_{l_4} + z_{l_1}z_{l_2}z_{l_4}z_{l_3} + z_{l_1}z_{l_4}z_{l_2}z_{l_3} + z_{l_4}z_{l_1}z_{l_2}z_{l_3}
		\lnAH
		+
		z_{l_1+l_4}z_{l_2}z_{l_3} + z_{l_1}z_{l_2+l_4}z_{l_3} + z_{l_1}z_{l_2}z_{l_3+l_4}
	}
	{\label{3.1_Lem3_Eq2}\rule{0cm}{15pt}
		z_{l_1}z_{l_2}*z_{l_3}z_{l_4}
	}
	{
		z_{l_1}z_{l_2}z_{l_3}z_{l_4} + z_{l_1}z_{l_3}z_{l_2}z_{l_4} + z_{l_1}z_{l_3}z_{l_4}z_{l_2} 
		\lnAH
		+ 
		z_{l_3}z_{l_1}z_{l_2}z_{l_4} + z_{l_3}z_{l_1}z_{l_4}z_{l_2} + z_{l_3}z_{l_4}z_{l_1}z_{l_2}
		\lnAH
		+
		z_{l_1+l_3}z_{l_2}z_{l_4} + z_{l_1+l_3}z_{l_4}z_{l_2} + z_{l_1}z_{l_2+l_3}z_{l_4}
		\lnAH[]
		+
		z_{l_3}z_{l_1+l_4}z_{l_2} + z_{l_1}z_{l_3}z_{l_2+l_4} + z_{l_3}z_{l_1}z_{l_2+l_4}
		+
		z_{l_1+l_3}z_{l_2+l_4}
		\nonumber
	}
	{\label{3.1_Lem3_Eq3}\rule{0cm}{15pt}
		z_{l_1}z_{l_2}*z_{l_3}*z_{l_4}
	}
	{
		z_{l_1}z_{l_2}*z_{l_3}z_{l_4} + z_{l_1}z_{l_2}*z_{l_4}z_{l_3}
		\lnAH
		+
		z_{l_3+l_4}z_{l_1}z_{l_2}  + z_{l_1}z_{l_3+l_4}z_{l_2} + z_{l_1}z_{l_2}z_{l_3+l_4}
		\lnAH[]
		+ 
		z_{l_1+l_3+l_4}z_{l_2} + z_{l_1}z_{l_2+l_3+l_4}
		\nonumber
	}
	{\label{3.1_Lem3_Eq4}\rule{0cm}{15pt}
		z_{l_1}*z_{l_2}*z_{l_3}*z_{l_4}
	}
	{
		z_{l_1}z_{l_2}*z_{l_3}*z_{l_4} + z_{l_2}z_{l_1}*z_{l_3}*z_{l_4}
		\lnAH
		+
		z_{l_1+l_2}z_{l_3}z_{l_4} + z_{l_1+l_2}z_{l_4}z_{l_3} + z_{l_3}z_{l_1+l_2}z_{l_4} + z_{l_4}z_{l_1+l_2}z_{l_3} 
		\lnAH
		+ 
		z_{l_3}z_{l_4}z_{l_1+l_2} + z_{l_4}z_{l_3}z_{l_1+l_2} + z_{l_1+l_2}z_{l_3+l_4} + z_{l_3+l_4}z_{l_1+l_2} 
		\lnAH[]
		+ 
		z_{l_1+l_2+l_3}z_{l_4} + z_{l_1+l_2+l_4}z_{l_3} + z_{l_3}z_{l_1+l_2+l_4} + z_{l_4}z_{l_1+l_2+l_3} 
		\lnAH
		+ 
		z_{l_1+l_2+l_3+l_4}
		\nonumber
	}
\end{lemma}
\envProof[\refLem{3.1_Lem1}]{
Identity \refEq{3.1_Lem1_Eq} follows from \refEq{2_PL_DefHarProdUnit} and \refEq{2_PL_DefHarProdMain} with $w_1=w_2=1$.
}
\envProof[\refLem{3.1_Lem2}]{
We see from \refEq{2_PL_DefHarProdUnit}, \refEq{2_PL_DefHarProdMain}, and \refEq{3.1_Lem1_Eq} that
	\envHLineFCm
	{
		z_{l_1}z_{l_2}*z_{l_3}
	}
	{
		z_{l_1} (z_{l_2}*z_{l_3}) + z_{l_3} (z_{l_1}z_{l_2}*1) + z_{l_1+l_3}(z_{l_2}*1)
	}
	{
		z_{l_1} ( z_{l_2}z_{l_3} + z_{l_3}z_{l_2} + z_{l_2+l_3} ) + z_{l_3}z_{l_1}z_{l_2} + z_{l_1+l_3}z_{l_2}
	}
	{
		z_{l_1}z_{l_2}z_{l_3} + z_{l_1}z_{l_3}z_{l_2} + z_{l_3}z_{l_1}z_{l_2} + z_{l_1+l_3}z_{l_2} + z_{l_1}z_{l_2+l_3}
	}
	which proves \refEq{3.1_Lem2_Eq1}.
We see from \refEq{3.1_Lem1_Eq} and \refEq{3.1_Lem2_Eq1} that
	\envHLineFCm
	{
		z_{l_1}*z_{l_2}*z_{l_3}
	}
	{
		(z_{l_1}z_{l_2}+z_{l_2}z_{l_1}+z_{l_1+l_2})*z_{l_3}
	}
	{
		z_{l_1}z_{l_2}*z_{l_3} + z_{l_2}z_{l_1}*z_{l_3} + z_{l_1+l_2}*z_{l_3}
	}
	{
		z_{l_1}z_{l_2}z_{l_3} + z_{l_1}z_{l_3}z_{l_2} + z_{l_3}z_{l_1}z_{l_2} + z_{l_1+l_3}z_{l_2} + z_{l_1}z_{l_2+l_3}
		\lnAH
		+
		z_{l_2}z_{l_1}z_{l_3} + z_{l_2}z_{l_3}z_{l_1} + z_{l_3}z_{l_2}z_{l_1} + z_{l_2+l_3}z_{l_1} + _{l_2}z_{l_1+l_3}
		\lnAH
		+
		z_{l_1+l_2}z_{l_3} + z_{l_3}z_{l_1+l_2} + z_{l_1+l_2+l_3}
	}
	which proves \refEq{3.1_Lem2_Eq2},
	and completes the proof.
}
\envProof[\refLem{3.1_Lem3}]{
We see from \refEq{2_PL_DefHarProdUnit}, \refEq{2_PL_DefHarProdMain}, and \refEq{3.1_Lem2_Eq1} that
	\envHLineFCm
	{
		z_{l_1}z_{l_2}z_{l_3}*z_{l_4}
	}
	{
		z_{l_1}(z_{l_2}z_{l_3}*z_{l_4}) + z_{l_4}(z_{l_1}z_{l_2}z_{l_3}*1) + z_{l_1+l_4}(z_{l_2}z_{l_3}*1)
	}
	{
		z_1(z_{l_2}z_{l_3}z_{l_4} + z_{l_2}z_{l_4}z_{l_3} + z_{l_4}z_{l_2}z_{l_3} + z_{l_2+l_4}z_{l_3} + z_{l_2}z_{l_3+l_4})
		\lnAH
		+
		z_{l_4}z_{l_1}z_{l_2}z_{l_3} + z_{l_1+l_4}z_{l_2}z_{l_3}
	}
	{
		z_{l_1}z_{l_2}z_{l_3}z_{l_4} + z_{l_1}z_{l_2}z_{l_4}z_{l_3} + z_{l_1}z_{l_4}z_{l_2}z_{l_3} + z_{l_4}z_{l_1}z_{l_2}z_{l_3}
		\lnAH
		+
		z_{l_1+l_4}z_{l_2}z_{l_3} + z_{l_1}z_{l_2+l_4}z_{l_3} + z_{l_1}z_{l_2}z_{l_3+l_4}
	}
	which proves \refEq{3.1_Lem3_Eq1}.
We see from \refEq{2_PL_DefHarProdMain}, \refEq{3.1_Lem1_Eq}, and \refEq{3.1_Lem2_Eq1} that
	\envHLineFCm
	{
		z_{l_1}z_{l_2}*z_{l_3}z_{l_4}
	}
	{
		z_{l_1}(z_{l_2}*z_{l_3}z_{l_4}) + z_{l_3}(z_{l_1}z_{l_2}*z_{l_4}) + z_{l_1+l_3}(z_{l_2}*z_{l_4})
	}
	{
		z_{l_1} (z_{l_3}z_{l_4}z_{l_2} + z_{l_3}z_{l_2}z_{l_4} + z_{l_2}z_{l_3}z_{l_4} + z_{l_3+l_2}z_{l_4} + z_{l_3}z_{l_4+l_2})
		\lnAH
		+
		z_{l_3} (z_{l_1}z_{l_2}z_{l_4} + z_{l_1}z_{l_4}z_{l_2} + z_{l_4}z_{l_1}z_{l_2} + z_{l_1+l_4}z_{l_2} + z_{l_1}z_{l_2+l_4})
		\lnAH
		+
		z_{l_1+l_3} (z_{l_2}z_{l_4}+z_{l_4}z_{l_2}+z_{l_2+l_4})
	}
	{
		z_{l_1}z_{l_2}z_{l_3}z_{l_4} + z_{l_1}z_{l_3}z_{l_2}z_{l_4} + z_{l_1}z_{l_3}z_{l_4}z_{l_2} + z_{l_3}z_{l_1}z_{l_2}z_{l_4} 
		\lnAH
		+ 
		z_{l_3}z_{l_1}z_{l_4}z_{l_2} + z_{l_3}z_{l_4}z_{l_1}z_{l_2} + z_{l_1+l_3}z_{l_2}z_{l_4} + z_{l_1+l_3}z_{l_4}z_{l_2}
		\lnAH
		+
		z_{l_1}z_{l_2+l_3}z_{l_4} + z_{l_3}z_{l_1+l_4}z_{l_2} + z_{l_1}z_{l_3}z_{l_2+l_4} + z_{l_3}z_{l_1}z_{l_2+l_4} + z_{l_1+l_3}z_{l_2+l_4}
	}
	which proves \refEq{3.1_Lem3_Eq2}.
We see from \refEq{3.1_Lem1_Eq} and \refEq{3.1_Lem2_Eq1} that
	\envHLineFCm
	{
		z_{l_1}z_{l_2}*z_{l_3}*z_{l_4}
	}
	{
		z_{l_1}z_{l_2}*(z_{l_3}z_{l_4}+z_{l_4}z_{l_3}+z_{l_3+l_4})
	}
	{
		z_{l_1}z_{l_2}*z_{l_3}z_{l_4} + z_{l_1}z_{l_2}*z_{l_4}z_{l_3} + z_{l_1}z_{l_2}*z_{l_3+l_4}
	}
	{
		z_{l_1}z_{l_2}*z_{l_3}z_{l_4} + z_{l_1}z_{l_2}*z_{l_4}z_{l_3}
		\lnAH
		+
		z_{l_1}z_{l_2}z_{l_3+l_4} + z_{l_1}z_{l_3+l_4}z_{l_2} + z_{l_3+l_4}z_{l_1}z_{l_2} 
		+ 
		z_{l_1+l_3+l_4}z_{l_2} + z_{l_1}z_{l_2+l_3+l_4}
	}
	which proves \refEq{3.1_Lem3_Eq3}.
We see from \refEq{3.1_Lem1_Eq} and \refEq{3.1_Lem2_Eq2} that
	\envHLineFCm
	{
		z_{l_1}*z_{l_2}*z_{l_3}*z_{l_4}
	}
	{
		(z_{l_1}z_{l_2}+z_{l_2}z_{l_1}+z_{l_1+l_2})*z_{l_3}*z_{l_4}
	}
	{
		z_{l_1}z_{l_2}*z_{l_3}*z_{l_4} + z_{l_2}z_{l_1}*z_{l_3}*z_{l_4} + z_{l_1+l_2}*z_{l_3}*z_{l_4}
	}
	{
		z_{l_1}z_{l_2}*z_{l_3}*z_{l_4} + z_{l_2}z_{l_1}*z_{l_3}*z_{l_4}
		\lnAH
		+
		z_{l_1+l_2}z_{l_3}z_{l_4} + z_{l_1+l_2}z_{l_4}z_{l_3} + z_{l_3}z_{l_1+l_2}z_{l_4} + z_{l_4}z_{l_1+l_2}z_{l_3} 
		\lnAH
		+ 
		z_{l_3}z_{l_4}z_{l_1+l_2} + z_{l_4}z_{l_3}z_{l_1+l_2} 
		+ 
		z_{l_1+l_2+l_3}z_{l_4} + z_{l_1+l_2+l_4}z_{l_3} + z_{l_3+l_4}z_{l_1+l_2} 
		\lnAH
		+
		z_{l_1+l_2}z_{l_3+l_4} + z_{l_3}z_{l_1+l_2+l_4} + z_{l_4}z_{l_1+l_2+l_3} 
		+ 
		z_{l_1+l_2+l_3+l_4}
	}
	which proves \refEq{3.1_Lem3_Eq4}, 
	and completes the proof.
}

We are now in a position to prove \refProp[s]{3.1_Prop1} and \ref{3.1_Prop2}.

\envProof[\refProp{3.1_Prop1}]{
Let $\lbfL_2=(l_1,l_2)$ be an index set in $\setN^2$. 
Applying the map $\mpEVH[]{}$ to both sides of \refEq{3.1_Lem1_Eq} and substituting $T=0$,
	we obtain
	\envOTLineThPd
	{
		\fcMZH{1}{l_1}\fcMZH{1}{l_2}
	}
	{
		\fcMZH{2}{l_1,l_2} + \fcMZH{2}{l_2,l_1} + \fcMZH{1}{l_1+l_2}
	}
	{
		\Sm{\sig\in\gpC{2}} \fcMZH{2}{l_{\sig^{-1}(1)},l_{\sig^{-1}(2)}} + \fcMZH{1}{l_1+l_2}
	}
We thus have
	\envHLineCm
	{
		{\fcMZHO{1}}^{\otimes2}(\lbfL_2)
	}
	{
		\racFA{\fcMZHO{2}}{ \mpS{\gpC{2}} }{\lbfL_2} + \fcMZHO{1}\circ\mpWTT{2}{\lbfL_2}
	}
	which proves \refEq{3.1_Prop1_Eq}
	because $\lbfL_2$ is arbitrary and $\fcMZHO{1}\circ\mpWTT{2}{\lbfL_2}=\fcMZO{1}\circ\mpWTT{2}{\lbfL_2}$ by $\mpWTT{2}{\lbfL_2}\geq2$.
}
\envProof[\refProp{3.1_Prop2}]{
Let $\lbfL_3=(l_1,l_2,l_3)$ be an index set in $\setN^3$. 
Applying the map $\mpEVH[]{}$ to both sides of \refEq{3.1_Lem2_Eq1} and substituting $T=0$,
	we obtain
	\envLineThCm
	{
		\fcMZH{2}{l_1,l_2}\fcMZH{1}{l_3}
	}
	{
		\fcMZH{3}{l_1,l_2,l_3} + \fcMZH{3}{l_1,l_3,l_2} + \fcMZH{3}{l_3,l_1,l_2} 
		+ 
		\fcMZH{2}{l_1+l_3,l_2} + \fcMZH{2}{l_1,l_2+l_3}
	}
	{
		\Sm{\sig\in\pstUT} \fcMZH{3}{l_{\sig^{-1}(1)},l_{\sig^{-1}(2)},l_{\sig^{-1}(3)}}
		+
		\fcMZH{2}{l_{\tau^{-1}(1)}+l_{\tau^{-1}(2)},l_{\tau^{-1}(3)}} + \fcMZH{2}{l_1,l_2+l_3}
	}
	where $\tau=(123)$.
We thus have
	\envHLineCm
	{
		\fcMZHO{2}\otimes\fcMZHO{1}(\lbfL_3)
	}
	{
		\racFA{\fcMZHO{3}}{ \mpS{\pstUT} }{\lbfL_3} 
		+ 
		\racFA{ \fcMZHO{2} \circ \empWT{3}{2,1} }{(123)}{\lbfL_3} 
		+ 
		\fcMZHO{2} \circ \empWTTh{3}{1,2}{\lbfL_3} 
	}
	which proves \refEq{3.1_Prop2_Eq1}.
In a similar way,
	we obtain from \refEq{3.1_Lem2_Eq2} that
	\envHLineCm
	{
		{\fcMZHO{1}}^{\otimes3}(\lbfL_3)
	}
	{
		\racFA{\fcMZHO{3}}{ \mpS{\gpSym{3}} }{\lbfL_3}
		+ 
		\racFA{ \fcMZHO{2} \circ (\empWT{3}{2,1} + \empWT{3}{1,2}) }{ \mpS{\gpC{3}} }{\lbfL_3}
		+
		\fcMZHO{1} \circ \mpWTT{3}{\lbfL_3}
	}
	which proves \refEq{3.1_Prop2_Eq2}.
}

We require another lemma for the proof of \refProp{3.1_Prop3}, 
	since the proof is more complicated than those of \refProp[s]{3.1_Prop1} and \ref{3.1_Prop2}.

\begin{lemma}\label{3.1_Lem4}
Let $\mpI=\mpI[4]$ mean the identity map on $\setN^4$.
We have the following equations in maps with the domain $\setN^4$:
\mbox{}\\\bfTx{(i)}
	\envHLineCSNme
	{\label{3.1_Lem4i_Eq1}\hspace{-20pt}
		\racF{ \empWT{4}{2,2} }{ (23)\mpS{\Gp{(34)}} }
	}
	{
		\racF{ \empWT{4}{2,2} }{ \mpS{ \pstWF[0] } }
	,}
	{\label{3.1_Lem4i_Eq2}\hspace{-20pt}\rule{0cm}{30pt}
		\racF{ \empWT{4}{i,j,k} }{ \mpS{\pstVF[0]} \mpS{\Gp{(34)}} }
	}
	{
		\envCaseTCm{
			\racF[n]{ \empWT{4}{i,j,k} }{ \mpS{\pstWF[1]} - (34) - (1324) }
			&
			((i,j,k)\in I)
		}{
			\racF[n]{ \empWT{4}{i,j,k} }{ \mpS{\pstWF[1]} - (24) - (1234) }
			&
			((i,j,k)=(1,2,1))
		}
	}
	{\label{3.1_Lem4i_Eq3}\hspace{-20pt}\rule{0cm}{20pt}
		\racF{\mpI}{ \mpS{\pstVF} \mpS{\Gp{(34)}} }
	}
	{
		\racF{ \mpI }{ \mpS{\pstWF} }
	,}
	where $I$ in \refEq{3.1_Lem4i_Eq2} means the set $\Set{(2,1,1),(1,1,2)}$.
\mbox{}\vWiden[5]\\
\bfTx{(ii)}
	\envPLine[????]{\HLineCECmNme[p]
	{\label{3.1_Lem4ii_Eq1}
		\racF{ \empWT{4}{3,1} }{ (24)\mpS{\Gp{(12)}} }
	}
	{
		\racF[n]{ \empWT{4}{3,1} }{ \mpS{\gpC{4}} - \gpu - (1234) }
	}
	{\label{3.1_Lem4ii_Eq2}\rule{0cm}{20pt}
		\racF{ \empWT{4}{1,3} }{ \mpS{\Gp{(12)}} }
	}
	{
		\racF[n]{ \empWT{4}{1,3} }{ \mpS{\gpC{4}} - (13)(24) - (1234) }
	}
	{\label{3.1_Lem4ii_Eq3}\rule{0cm}{20pt}
		\racF{ \empWT{4}{2,2} }{ \mpS{\pstWF[0]} \mpS{\Gp{(12)}} }
	}
	{
		\racF[n]{ \empWT{4}{2,2} }{ \mpS{\pstXF} - \gpu - (13)(24) }
	}
	{\label{3.1_Lem4ii_Eq4}\rule{0cm}{20pt}
		\racF{ \empWT{4}{2,1,1} }{ \mpS{\pstWFO[1]{(34)}} \mpS{\Gp{(12)}} } 
	}
	{
		\racF[n]{ \empWT{4}{2,1,1} }{ \mpS{\gpAlt{4}} - \gpu - (12)(34) } 
	}\HLineCSNmePd
	{\label{3.1_Lem4ii_Eq5}\rule{0cm}{20pt}
		\racF{ \empWT{4}{1,2,1} }{ \mpS{\pstWFO[1]{(1234)}} \mpS{\Gp{(12)}} } 
	}
	{
		\racF[n]{ \empWT{4}{1,2,1} }{ \mpS{\gpAlt{4}} - (123) - (134) } 
	}
	{\label{3.1_Lem4ii_Eq6}\rule{0cm}{20pt}
		\racF{ \empWT{4}{1,1,2} }{ \mpS{\pstWFO[1]{(1324)}} \mpS{\Gp{(12)}} } 
	}
	{
		\racF[n]{ \empWT{4}{1,1,2} }{ \mpS{\gpAlt{4}} - (13)(24) - (14)(23) } 
	}
	{\label{3.1_Lem4ii_Eq7}\rule{0cm}{20pt}
		\racF{\mpI}{ \mpS{\pstWF} \mpS{\Gp{(12)}} }
	}
	{
		\racF{\mpI}{ \mpS{\gpSym{4}} }
	}
	}
\end{lemma}

We now prove \refProp{3.1_Prop3}.
We will then discuss a proof of \refLem{3.1_Lem4}.

\envProof[\refProp{3.1_Prop3}]{
Let $\lbfL_4=(l_1,l_2,l_3,l_4)$ be an index set in $\setN^4$. 
Applying the map $\mpEVH[]{}$ to both sides of \refEq{3.1_Lem3_Eq1} and substituting $T=0$,
	we obtain 
	\envLineThCm
	{
		\fcMZH{3}{l_1,l_2,l_3}\fcMZH{1}{l_4}
	}
	{
		\fcMZH{4}{l_1,l_2,l_3,l_4} + \fcMZH{4}{l_1,l_2,l_4,l_3} + \fcMZH{4}{l_1,l_4,l_2,l_3} + \fcMZH{4}{l_4,l_1,l_2,l_3}
		\lnAH[]
		+
		\fcMZH{3}{l_1+l_4,l_2,l_3} + \fcMZH{3}{l_1,l_2+l_4,l_3} + \fcMZH{3}{l_1,l_2,l_3+l_4}
	}
	{
		\Sm{\sig\in\pstUF} \fcMZH{4}{l_{\sig^{-1}(1)},l_{\sig^{-1}(2)},l_{\sig^{-1}(3)},l_{\sig^{-1}(4)}}
		\lnAH
		+
		\fcMZH{3}{l_{\tau^{-1}(1)}+l_{\tau^{-1}(2)},l_{\tau^{-1}(3)},l_{\tau^{-1}(4)}} + \fcMZH{3}{l_{\tau^{-1}(1)},l_{\tau^{-1}(2)}+l_{\tau^{-1}(3)},l_{\tau^{-1}(4)}} 
		\lnAH
		+ 
		\fcMZH{3}{l_1,l_2,l_3+l_4}
	}
	where $\tau=(234)$.
We thus have
	\envHLineCm
	{\hspace{-10pt}
		\fcMZHO{3}\otimes\fcMZHO{1}(\lbfL_4)
	}
	{
		\racFA{\fcMZHO{4}}{ \mpS{\pstUF} }{\lbfL_4}
		+
		\fcMZHO{3} \circ \Fcr{ \racFr{ \empWT{4}{2,1,1} + \empWT{4}{1,2,1} }{ (234) } + \empWT{4}{1,1,2} }{\lbfL_4}
	}
	which proves \refEq{3.1_Prop3_Eq1}.
Similarly,
	we have by \refEq{3.1_Lem3_Eq2} 
	\envLineCm
	{
		{\fcMZHO{2}}^{\otimes2}(\lbfL_4)
	}
	{
		\racFA{\fcMZHO{4}}{ \mpS{\pstVF} }{\lbfL_4}
		+
		\racFA{ \fcMZHO{3} \circ (\empWT{4}{2,1,1} + \empWT{4}{1,2,1} + \empWT{4}{1,1,2}) }{ \mpS{\pstVF[0]} }{\lbfL_4}
		+
		\racFA{ \fcMZHO{2} \circ \empWT{4}{2,2} }{ (23) }{\lbfL_4}
	}
	which proves \refEq{3.1_Prop3_Eq2}.

As we have proved \refEq{3.1_Prop3_Eq1} and \refEq{3.1_Prop3_Eq2}, 
	we can deduce from \refEq{3.1_Lem3_Eq3} and \refEq{3.1_Lem3_Eq4} that
	\envHLine[3.1_Prop3Pr_Eq3]
	{
		\fcMZHO{2}\otimes{\fcMZHO{1}}^{\otimes2}
	}
	{
		\racF{{\fcMZHO{2}}^{\otimes2}}{ \mpS{\Gp{(34)}} }
		+
		\fcMZHO{3} \circ \bkR[b]{ \racF{ \empWT{4}{2,1,1} }{ (1324) } + \racF{ \empWT{4}{1,2,1} }{ (24) } + \racF{ \empWT{4}{1,1,2} }{ (34) } }
		\lnAH
		+
		\fcMZHO{2} \circ \bkR[b]{ \racF{ \empWT{4}{3,1} }{ (24) } + \empWT{4}{1,3} }
	}
	and
	\envHLineCm[3.1_Prop3Pr_Eq4]
	{
		{\fcMZHO{1}}^{\otimes4}
	}
	{
		\racF{\fcMZHO{2}\otimes{\fcMZHO{1}}^{\otimes2}}{ \mpS{\Gp{(12)}} }
		+
		\fcMZHO{3} \circ \bkR[b]{ 
			\racF[n]{ \empWT{4}{2,1,1} }{ \gpu + (12)(34)} 
			+ 
			\racF[n]{ \empWT{4}{1,2,1} }{ (123)+(134) } 
			\lnAH
			+ 
			\racF[n]{ \empWT{4}{1,1,2} }{ (13)(24)+(14)(23) } 
		}
		+
		\fcMZHO{2} \circ \bkR[b]{ 
			\racF[n]{ \empWT{4}{2,2} }{ \gpu + (13)(24) } 
			\lnAH[]
			+ 
			\racF[n]{ \empWT{4}{3,1} }{ \gpu+(1234) } 
			+ 
			\racF[n]{ \empWT{4}{1,3} }{ (13)(24)+(1234) } 
			\nonumber
		}
		+
		\fcMZO{1}  \circ \mpWT{4}
	}
	respectively.
Combining \refEq{3.1_Prop3_Eq2} and the equations in (i) of \refLem{3.1_Lem4},
	we obtain
	\envLineThPd
	{
		\racF{{\fcMZHO{2}}^{\otimes2}}{ \mpS{\Gp{(34)}} }
	}
	{
		\racF{\fcMZHO{4}}{ \mpS{\pstVF}\mpS{\Gp{(34)}} } 
		+
		\racF{ \fcMZHO{3} \circ (\empWT{4}{2,1,1} +  \empWT{4}{1,2,1} + \empWT{4}{1,1,2}) }{ \mpS{\pstVF[0]}\mpS{\Gp{(34)}} }
		+
		\racF{ \fcMZHO{2} \circ \empWT{4}{2,2} }{ (23)\mpS{\Gp{(34)}} }
	}
	{
		\racF{ \fcMZHO{4} }{ \mpS{\pstWF} }
		+
		\fcMZHO{3} \circ \bkR[b]{
			\racF[n]{ (\empWT{4}{2,1,1} + \empWT{4}{1,1,2}) }{ \mpS{\pstWF[1]} - (34) - (1324) }
			\lnAH
			+
			\racF[n]{ \empWT{4}{1,2,1} }{ \mpS{\pstWF[1]} - (24) - (1234) }
		}
		+
		\racF{ \fcMZHO{2} \circ \empWT{4}{2,2} }{ \mpS{\pstWF[0]} }
	}
Substituting this into the right-hand side of \refEq{3.1_Prop3Pr_Eq3} gives
	\envHLineCm
	{
		\fcMZHO{2}\otimes{\fcMZHO{1}}^{\otimes2}
	}
	{
		\racF{\fcMZHO{4}}{ \mpS{\pstWF} } 
		+
		\racF{ \fcMZHO{3} \circ (\empWT{4}{2,1,1} + \empWT{4}{1,2,1} + \empWT{4}{1,1,2}) }{ \mpS{\pstWF[1]} }
		\lnAH
		-
		\fcMZHO{3} \circ \bkR[b]{ \racF{ \empWT{4}{2,1,1} }{ (34) } + \racF{ \empWT{4}{1,2,1} }{ (1234) } + \racF{ \empWT{4}{1,1,2} }{ (1324) } }
		\lnAH[]
		+
		\fcMZHO{2} \circ  \bkR[b]{ \racF{ \empWT{4}{2,2} }{ \mpS{\pstWF[0]} } + \racF{ \empWT{4}{3,1} }{ (24) } + \empWT{4}{1,3} }
	}
	which proves \refEq{3.1_Prop3_Eq3},
	since
	\envM
	{
		\mpS{\pstWF[1]} - \sig
	}
	{
		\mpS{\pstWF[1]\setminus\Set{\sig}}
	}
	for $\sig\in\pstWF[1]$.
Similarly,
	combining \refEq{3.1_Prop3_Eq3} and the equations in (ii) of \refLem{3.1_Lem4},
	we obtain
	\envLinePd
	{
		\racF{\fcMZHO{2}\otimes{\fcMZHO{1}}^{\otimes2}}{ \mpS{\Gp{(12)}} }
	}
	{
		\racF{\fcMZHO{4}}{ \mpS{\gpSym{4}} } 
		+
		\fcMZHO{3} \circ \bkR[b]{
			\racF[n]{ \empWT{4}{2,1,1} }{ \mpS{\gpAlt{4}} - \gpu - (12)(34) }
			\lnAH
			+
			\racF[n]{ \empWT{4}{1,2,1} }{ \mpS{\gpAlt{4}} - (123) - (134) }
			+
			\racF[n]{ \empWT{4}{1,1,2} }{ \mpS{\gpAlt{4}} - (13)(24) - (14)(23) }
		}
		\lnAH
		+
		\fcMZHO{2} \circ  \bkR[b]{ 
			\racF[n]{ \empWT{4}{2,2} }{ \mpS{\pstXF} - \gpu - (13)(24) } 
			+ 
			\racF[n]{ \empWT{4}{3,1} }{ \mpS{\gpC{4}} - \gpu - (1234) } 
			\lnAH
			+ 
			\racF[n]{ \empWT{4}{1,3} }{ \mpS{\gpC{4}} - (13)(24) - (1234) } 
		}
	}
Substituting this into the right-hand side of \refEq{3.1_Prop3Pr_Eq4} proves \refEq{3.1_Prop3_Eq4}.
}

We will show \refLem{3.1_Lem4} for the completeness of the proof of \refProp{3.1_Prop3}.

For a subgroup $H$ in $\gpSym{4}$,
	we define an equivalence relation $\equiv$ on $\gpSym{4}$ such that $\sig \equiv\tau\ \modd H$ if and only if $\sig\tau^{-1}\in H$,
	and we denote by $\Eqv{\sig}{H}$ the equivalence class of $\sig$.
Note that $\Eqv{\sig}{H}$ is the right coset $H\sig$ of $\gpSym{4}$.
We extend the relation $\equiv$ on $\gpSym{4}$ to the congruence relation on its group ring $\setZ[\gpSym{4}]$ in the standard way such that
	 \envMPt{\equiv}{
		\SmN[t] a_i\sig_i
	}{
		\SmN[t] a_i \tau_i
	 	\ \modd H
	}
	if and only if $\sig_i \equiv\tau_i\ \modd H$ for every $i$.
\refTab{3.1_PL_Table1} below gives all the equivalence classes in $\gpSym{4}$ modulo certain subgroups,
	where
	we denote by $\Gp{\sig_1,\ldots,\sig_i}$ the subgroup generated by permutations $\sig_1,\ldots,\sig_i$.
(We have already used  $\Gp{\sig}$ to denote a cyclic subgroup.)
The equivalence classes in the table will be necessary when we prove some congruence equations in $\setZ[\gpSym{4}]$.

The following congruence equations in $\setZ[\gpSym{4}]$ are useful  for proving \refLem{3.1_Lem4}.

\begin{lemma}\label{3.1_Lem5}\newcommand{\tmpNm}{-25pt}
\mbox{}\\\bfTx{(i)}
	\envHLineCSNmePt[???]{\equiv}
	{\hspace{\tmpNm}\label{3.1_Lem5i_Eq1}
		(23)\mpS{\Gp{(34)}}
	}
	{
		\mpS{ \pstWF[0] }
		\qquad\qquad\qquad\qquad\;\;\,
		\modd \Gp{(12),(34)}
		,
	}
	{\hspace{\tmpNm}\label{3.1_Lem5i_Eq2}\rule{0cm}{30pt}
		\mpS{\pstVF[0]} \mpS{\Gp{(34)}}
	}
	{
		\envCaseTCm{
			\mpS{\pstWF[1]} - (34) - (1324)
			&
			\modd \Gp{(12)} \;\;\mathrm{or}\;\; \modd \Gp{(34)}
		}{
			\mpS{\pstWF[1]} - (24) - (1234) 
			&
			\modd \Gp{(23)}
		}
	}
	{\hspace{\tmpNm}\label{3.1_Lem5i_Eq3}\rule{0cm}{20pt}
		\mpS{\pstVF} \mpS{\Gp{(34)}}
	}
	{
		\mpS{\pstWF}
		\qquad\qquad\qquad\qquad\;\;\,\,
		\modd \Gp{\gpu}
		.
	}
\bfTx{(ii)}
	\envPLine[????]{\HLineCTCmNmePt[p]{\equiv}
	{\hspace{\tmpNm}\label{3.1_Lem5ii_Eq1}
		(24)\mpS{\Gp{(12)}}
	}
	{
		\mpS{\gpC{4}} - \gpu - (1234)
		\qquad\qquad\qquad
		\modd \Gp{(12), (123)}
	}\HLineCSCmNmePt[p]{\equiv}
	{\hspace{\tmpNm}\label{3.1_Lem5ii_Eq2}\rule{0cm}{20pt}
		\mpS{\Gp{(12)}}
	}
	{
		\mpS{\gpC{4}} - (13)(24) - (1234)
		\qquad\ \;\,\,
		\modd \Gp{(23), (234)}
	}
	{\hspace{\tmpNm}\label{3.1_Lem5ii_Eq3}\rule{0cm}{20pt}
		\mpS{\pstWF[0]} \mpS{\Gp{(12)}}
	}
	{
		\mpS{\pstXF} - \gpu - (13)(24) 
		\qquad\qquad\quad
		\modd \Gp{(12), (34)}
	}
	{\hspace{\tmpNm}\label{3.1_Lem5ii_Eq4}\rule{0cm}{20pt}
		\mpS{\pstWFO[1]{(34)}} \mpS{\Gp{(12)}}
	}
	{
		\mpS{\gpAlt{4}} - \gpu - (12)(34) 
		\qquad\qquad\;\;\;\,\,
		\modd \Gp{(12)}
	}\HLineCSNmePdPt{\equiv}
	{\hspace{\tmpNm}\label{3.1_Lem5ii_Eq5}\rule{0cm}{20pt}
		\mpS{\pstWFO[1]{(1234)}} \mpS{\Gp{(12)}}
	}
	{
		\mpS{\gpAlt{4}} - (123) - (134) 
		\qquad\qquad\;\;\,
		\modd \Gp{(23)}
	}
	{\hspace{\tmpNm}\label{3.1_Lem5ii_Eq6}\rule{0cm}{20pt}
		\mpS{\pstWFO[1]{(1324)}} \mpS{\Gp{(12)}}
	}
	{
		\mpS{\gpAlt{4}} - (13)(24) - (14)(23) 
		\qquad\;
		\modd \Gp{(34)}
	}
	{\hspace{\tmpNm}\label{3.1_Lem5ii_Eq7}\rule{0cm}{20pt}
		\mpS{\pstWF} \mpS{\Gp{(12)}}
	}
	{
		\mpS{\gpSym{4}} 
		\qquad\qquad\qquad\qquad\qquad\quad\;\,\,
		\modd \Gp{\gpu}
	}
	}
\end{lemma}
\begin{table}[!t]\renewcommand{\arraystretch}{1.3}\begin{center}{\footnotesize
\caption{All equivalence classes (or all right cosets $H\sig$) in $\gpSym{4}$ modulo subgroups $H$.}\label{3.1_PL_Table1}
\mbox{}\\
	\begin{tabular}{|c|llll|}  \hline
		\bfTx{mod} 		&\multicolumn{4}{l|}{\bfTx{All equivalence classes}}								\\\hline  
		$\Gp{(12), (123)}$ 	&\multicolumn{4}{l|}{$\Set{\gpu, (12), (13), (23), (123), (132)}$,}		\\
						&\multicolumn{4}{l|}{$\Set{(14), (14)(23), (142), (143), (1423), (1432)}$,}		\\
						&\multicolumn{4}{l|}{$\Set{(24), (13)(24), (124), (243), (1243), (1324)}$,}		\\
						&\multicolumn{4}{l|}{$\Set{(34), (12)(34), (134), (234), (1234), (1342)}$.}		\\\hline
		$\Gp{(23), (234)}$ 	&\multicolumn{4}{l|}{$\Set{\gpu, (23), (24), (34), (234), (243)}$,}				\\
				 		&\multicolumn{4}{l|}{$\Set{(12), (12)(34), (132), (142), (1342), (1432)}$,}		\\
				 		&\multicolumn{4}{l|}{$\Set{(13), (13)(24), (123), (143), (1243), (1423)}$,}		\\
				 		&\multicolumn{4}{l|}{$\Set{(14), (14)(23), (124), (134), (1234), (1324)}$.}		\\\hline
		$\Gp{(12), (34)}$	&\multicolumn{2}{l}{$\Set{\gpu, (12), (34), (12)(34)}$,} 		& \multicolumn{2}{l|}{$\Set{(13), (132), (143), (1432)}$,}					\\
				   		&\multicolumn{2}{l}{$\Set{(14), (134), (142), (1342)}$,}  	& \multicolumn{2}{l|}{$\Set{(23), (123), (243), (1243)}$,}		\\
				   		&\multicolumn{2}{l}{$\Set{(24), (124), (234), (1234)}$,}		& \multicolumn{2}{l|}{$\Set{(13)(24), (14)(23), (1324), (1423)}$.}	\\\hline
		$\Gp{(12)}$    		&\Set{\gpu, (12)},		&\Set{(13),(132)},		&\Set{(14), (142)},		&\Set{(23), (123)},			\\
				   		&\Set{(24),(124)},		&\Set{(34), (12)(34)},		&\Set{(13)(24), (1324)},	&\Set{(14)(23), (1423)},			\\
						&\Set{(134), (1342)},		&\Set{(143),(1432)},		&\Set{(234),(1234)},		&\Set{(243), (1243)}.	\\\hline
		$\Gp{(23)}$    		&\Set{\gpu, (23)},		&\Set{(12), (132)},		&\Set{(13), (123)},		&\Set{(14), (14)(23)},		\\
				   		&\Set{(24), (243)},		&\Set{(34), (234)},		&\Set{(12)(34), (1342)},	&\Set{(13)(24), (1243)},	\\
						&\Set{(124), (1324)},		&\Set{(134), (1234)},		&\Set{(142), (1432)},		&\Set{(143), (1423)}.		\\\hline
		$\Gp{(34)}$    		&\Set{\gpu, (34)},		&\Set{(12), (12)(34)},		&\Set{(13), (143)},		&\Set{(14), (134)},			\\
				   		&\Set{(23), (243)},		&\Set{(24), (234)},		&\Set{(13)(24), (1423)},	&\Set{(14)(23), (1324)},			\\
						&\Set{(123), (1243)},		&\Set{(124), (1234)},		&\Set{(132), (1432)},		&\Set{(142), (1342)}.	\\\hline
		$\Gp{(13)(24)}$    	&\Set{\gpu, (13)(24)},	&\Set{(12), (1423)},		&\Set{(13), (24)},		&\Set{(14), (1243)},			\\
				   		&\Set{(23), (1342)},		&\Set{(34), (1324)},		&\Set{(12)(34), (14)(23)},	&\Set{(123), (142)},			\\
						&\Set{(124), (143)},		&\Set{(132), (234)},		&\Set{(134), (243)},		&\Set{(1234), (1432)}.	\\\hline
	\end{tabular} 
}\end{center}\end{table}

\envProof{
Before proving the congruence equations,
	we introduce an identity in $\setZ[\gpSym{n}]$, 
	which immediately follows from the definition:
	\envHLineCm[3.1_Lem5iPr_Eq00]
	{
		\mpS{H}\mpS{K}
	}
	{
		\mpS{H_1}\mpS{K} + \cdots + \mpS{H_n}\mpS{K}
	}
	where 
	$H$ and $K$ are subsets in $\gpSym{n}$
	such that
	$H_1,\ldots,H_n$ are a partition of $H$
	(i.e., a set of subsets of $H$ satisfying $\UnT[t]{i=1}{n}H_i=H$ and $H_i\cap H_j=\phi$ for $i\neq j$).
	
We first prove the congruence equations stated in (i).
We obtain from 
	\envM
	{
		\mpS{\Gp{(34)}}
	}
	{
		\gpu + (34)
	} 
	that 
	\envHLinePd[3.1_Lem5iPr_Eq1]
	{
		(23)\mpS{\Gp{(34)}}
	}
	{
		(23) + (234) 
	}
Since $\Set{(24), (124), (234), (1234)}$ is an equivalence class modulo $\Gp{(12), (34)}$ as we see in \refTab{3.1_PL_Table1},
	\envHLinePdPt{\equiv}
	{
		(234)
	}
	{
		(24)
		\qquad
		\modd\Gp{(12), (34)}
	}
Thus,
	noting the definition of $\pstWF[0]$ in \refEq{3.1_PL_DefGpDpt4W},
	we have
	\envOTLineThCmPt{\equiv}
	{
		(23)\mpS{\Gp{(34)}}
	}
	{
		(23) + (24) 
	}
	{
		\mpS{\pstWF[0]}
		\qquad
		\modd \Gp{(12),(34)}
	}
	which proves \refEq{3.1_Lem5i_Eq1}.
A calculation shows that
	\envHLineCm[3.1_Lem5iPr_Eq2help]
	{
		(1243)\mpS{\Gp{(34)}} 
	}
	{
		(1243) + (124)
	}
	and so we see from \refEq{3.1_Lem5iPr_Eq00}, \refEq{3.1_Lem5iPr_Eq1}, and \refEq{3.1_Lem5iPr_Eq2help} that 
	\envOTLinePd[3.1_Lem5iPr_Eq2]
	{
		\mpS{\pstVF[0]} \mpS{\Gp{(34)}}
	}
	{
		(23)\mpS{\Gp{(34)}} + (1243)\mpS{\Gp{(34)}} 
	}
	{
		(23) + (124) + (234) + (1243)
	}
Using \refEq{3.1_Lem5iPr_Eq2} and the equivalence classes modulo $\Gp{(12)}$, $\Gp{(23)}$, and $\Gp{(34)}$ in \refTab{3.1_PL_Table1}, 
	we obtain 
	\envHLineThPt{\equiv}
	{
		\mpS{\pstVF[0]} \mpS{\Gp{(34)}}
	}
	{
		\envCaseTCm{
			(23) + (24) + (1234) + (1243)
			\;&
			\modd \Gp{(12)} \text{\;\;\;or\;\;\;} \modd \Gp{(34)}
		}{
			(23) + (34) + (1243) + (1324)
			&
			\modd \Gp{(23)}
		}
	}
	{
		\envCaseTCm{
			\mpS{\pstWF[1]} - (34) - (1324)
			\;\hspace{41pt}&
			\modd \Gp{(12)} \text{\;\;\;or\;\;\;} \modd \Gp{(34)}
		}{
			\mpS{\pstWF[1]} - (24) - (1234)
			&
			\modd \Gp{(23)}
		}
	}
	which proves \refEq{3.1_Lem5i_Eq2}.	
Direct calculations show that
	\envHLineCSCm
	{
		(13)(24)\mpS{\Gp{(34)}}
	}
	{
		(13)(24) + (1324)
	}
	{
		(123)\mpS{\Gp{(34)}}
	}
	{
		(123) + (1234)
	}
	{
		(243)\mpS{\Gp{(34)}}
	}
	{
		(243) + (24)
	}
	which together with \refEq{3.1_Lem5iPr_Eq2} yields
	\envHLineThPd
	{
		\mpS{\pstVF} \mpS{\Gp{(34)}}
	}
	{
		\mpS{\Set{\gpu,(13)(24),(123),(243)}}\mpS{\Gp{(34)}} + \mpS{\pstVF[0]}\mpS{\Gp{(34)}}
	}
	{
		\gpu + (34) + (13)(24) + (1324) + (123) + (1234) + (243) + (24)
		\lnAH
		+ 
		(23) + (124) + (234) + (1243)
	}
We obtain \refEq{3.1_Lem5i_Eq3}
	because 
	the right-hand side of this equation is
	\envMOCm
	{
		\mpS{\pstWF}
	}
	by definition.

We next prove the congruence equations stated in (ii).
We easily see that
	\envHLineCFLaaPd[3.1_Lem5iiPr_Eq1]
	{
		(24)\mpS{\Gp{(12)}}
	}
	{
		(24) + (142)
	}
	{
		\mpS{\Gp{(12)}}
	}
	{
		\gpu + (12)
	}
Using these equations and the equivalence classes modulo $\Gp{(12), (123)}$ and $\Gp{(23), (234)}$ in \refTab{3.1_PL_Table1}, 
	we obtain 
	\envOTLinePt{\equiv}
	{
		(24)\mpS{\Gp{(12)}}
	}
	{
		(13)(24) + (1432)
	}
	{
		\mpS{\gpC{4}} - \gpu - (1234)
		\qquad
		\modd \Gp{(12), (123)}
	}
	and
	\envOTLineCmPt{\equiv}
	{
		\mpS{\Gp{(12)}}
	}
	{
		\gpu + (1432)
	}
	{
		\mpS{\gpC{4}} - (13)(24) - (1234)
		\qquad\hspace{1pt}
		\modd \Gp{(23), (234)}
	}
	which prove \refEq{3.1_Lem5ii_Eq1} and \refEq{3.1_Lem5ii_Eq2},
	respectively.
A direct calculation shows that
	\envHLineCm[3.1_Lem5iiPr_Eq2help]
	{
		(23)\mpS{\Gp{(12)}}
	}
	{
		(23) + (132)
	}
	and so we see from \refEq{3.1_Lem5iiPr_Eq1} and \refEq{3.1_Lem5iiPr_Eq2help} that
	\envOTLinePd[3.1_Lem5iiPr_Eq2]
	{
		\mpS{\pstWF[0]} \mpS{\Gp{(12)}}
	}
	{	
		(23)\mpS{\Gp{(12)}} + (24)\mpS{\Gp{(12)}}
	}
	{
		(23) + (24) + (132) + (142)
	}
Using \refEq{3.1_Lem5iiPr_Eq2} and the equivalence classes modulo  $\Gp{(12), (34)}$ in \refTab{3.1_PL_Table1}, 
	we obtain 
	\envOTLineThCmPt{\equiv}
	{
		\mpS{\pstWF[0]} \mpS{\Gp{(12)}}
	}
	{
		(23) + (1234) + (1432) + (14)
	}
	{
		\mpS{\pstXF} - \gpu - (13)(24)
		\quad
		\modd \Gp{(12), (34)}
	}
	which proves \refEq{3.1_Lem5ii_Eq3}.
Direct calculations show that  
	\envHLineCECm[3.1_Lem5iiPr_Eq3help1]
	{
		(34)\mpS{\Gp{(12)}}
	}
	{
		(34) + (12)(34)
	}
	{
		(1234)\mpS{\Gp{(12)}}
	}
	{
		(1234) + (134)
	}
	{
		(1243)\mpS{\Gp{(12)}}
	}
	{
		(1243) + (143)
	}
	{
		(1324)\mpS{\Gp{(12)}}
	}
	{
		(1324) + (14)(23)
	}
	and so we see from \refEq{3.1_Lem5iiPr_Eq2} and \refEq{3.1_Lem5iiPr_Eq3help1} that
	\envHLineThCm
	{
		\mpS{\pstWF[1]} \mpS{\Gp{(12)}}
	}
	{
		\mpS{\Set{(34),(1234),(1243),(1324)}} \mpS{\Gp{(12)}} + \mpS{\pstWF[0]} \mpS{\Gp{(12)}}
	}
	{
		(34) + (12)(34) + (1234) + (134) + (1243) + (143) + (1324) + (14)(23)
		\lnAH
		+
		(23) + (24) + (132) + (142)
	}
	which can be restated as
	\envHLinePd[3.1_Lem5iiPr_Eq3]
	{
		\mpS{\pstWF[1]} \mpS{\Gp{(12)}}
	}
	{
		(23) + (24) + (34) + (12)(34) + (14)(23)
		\lnAH
		+
		(132) + (134) + (142) + (143) + (1234) + (1243) + (1324)
	}
Equation \refEq{3.1_Lem5iiPr_Eq3} together with the first equation of \refEq{3.1_Lem5iiPr_Eq3help1} gives
	\envHLineFPd
	{
		\mpS{\pstWFO[1]{(34)}} \mpS{\Gp{(12)}}
	}
	{
		\mpS{\pstWF[1]\setminus\Set{(34)}} \mpS{\Gp{(12)}}
	}
	{
		\mpS{\pstWF[1]} \mpS{\Gp{(12)}} - (34)\mpS{\Gp{(12)}}
	}
	{
		(23) + (24) + (14)(23)
		\lnAH[]
		+
		(132) + (134) + (142) + (143) + (1234) + (1243) + (1324)
		\nonumber
	}
Using this equation and the equivalence classes modulo $\Gp{(12)}$ in \refTab{3.1_PL_Table1},
	we obtain 
	\envLineCmPt{\equiv}
	{
		\mpS{\pstWFO[1]{(34)}} \mpS{\Gp{(12)}}
	}
	{
		(123) + (124) + (14)(23)
		\lnAH
		+
		(132) + (134) + (142) + (143) + (234) + (243) + (13)(24)
		\qquad
		\modd \Gp{(12)} 
	}
	which proves \refEq{3.1_Lem5ii_Eq4}
	since
	\envHLinePd[3.1_Lem5iiPr_Eq3help2]
	{
		\mpS{\gpAlt{4}}
	}
	{
		\gpu + (12)(34) + (13)(24) + (14)(23) 
		\lnAH
		+	
		(123) + (124) + (132) + (134) + (142) + (143) + (234) + (243)
	}
Similarly,
	\refEq{3.1_Lem5iiPr_Eq3} together with
	the second and fourth equations of \refEq{3.1_Lem5iiPr_Eq3help1} and the equivalence classes modulo $\Gp{(23)}$ and $\Gp{(34)}$ in \refTab{3.1_PL_Table1} yield
	\envHLineFPt{\equiv}
	{
		\mpS{\pstWFO[1]{(1234)}} \mpS{\Gp{(12)}}
	}
	{
		(23) + (24) + (34) + (12)(34) + (14)(23)
		\lnAH
		+
		(132) + (142) + (143) + (1243) + (1324)
	}
	{
		\gpu + (243) + (234) + (12)(34) + (14)(23)
		\lnAH
		+
		(132) + (142) + (143) + (13)(24) + (124)
	}
	{
		\mpS{\gpAlt{4}} - (123) - (134)
		\qquad\qquad\qquad\qquad\qquad\qquad
		\modd \Gp{(23)}
	}
	and
	\envHLineFCmPt{\equiv}
	{
		\mpS{\pstWFO[1]{(1324)}} \mpS{\Gp{(12)}}
	}
	{
		(23) + (24) + (34) + (12)(34) 
		\lnAH
		+
		(132) + (134) + (142) + (143) + (1234) + (1243) 
	}
	{
		(243) + (234) + \gpu + (12)(34) 
		\lnAH
		+
		(132) + (134) + (142) + (143) + (124) + (123)
	}
	{
		\mpS{\gpAlt{4}} - (13)(24) - (14)(23)
		\qquad\qquad\qquad\qquad\qquad
		\modd \Gp{(34)}
	}
	respectively,
	which prove \refEq{3.1_Lem5ii_Eq5} and \refEq{3.1_Lem5ii_Eq6}.
Direct calculations show that
	\envPLine[3.1_Lem5iiPr_Eq4help]{\HLineCSCm[p]
	{
		(13)(24)\mpS{\Gp{(12)}}
	}
	{
		(13)(24) + (1423)
	}
	{
		(123)\mpS{\Gp{(12)}}
	}
	{
		(123) + (13)
	}
	{
		(124)\mpS{\Gp{(12)}}
	}
	{
		(124) + (14)
	}\HLineCFCm[*]
	{
		(234)\mpS{\Gp{(12)}}
	}
	{
		(234) + (1342)
	}
	{
		(243)\mpS{\Gp{(12)}}
	}
	{
		(243) + (1432)
	}
	}
	and so we see from \refEq{3.1_Lem5iiPr_Eq3} and \refEq{3.1_Lem5iiPr_Eq4help} that
	\envOTLineThCm[3.1_Lem5iiPr_Eq4]
	{
		\mpS{\pstWF} \mpS{\Gp{(12)}}
	}
	{
		\mpS{\Set{\gpu,(13)(24),(123),(124),(234),(243)}} \mpS{\Gp{(12)}}
		+
		\mpS{\pstWF[1]} \mpS{\Gp{(12)}}
	}
	{
		\mpS{\gpSym{4}}
	}
	which proves \refEq{3.1_Lem5ii_Eq7},
	and 
	completes the proof.
}

The following statement holds:
	the maps $\empWT{4}{3,1}$, $\empWT{4}{1,3}$, $\empWT{4}{2,2}$, $\empWT{4}{2,1,1}$, $\empWT{4}{1,2,1}$, and $\empWT{4}{1,1,2}$
	are invariant under the subgroups $\Gp{(12),(123)}$, $\Gp{(23),(234)}$, $\Gp{(12),(34)}$, $\Gp{(12)}$, $\Gp{(23)}$, and $\Gp{(34)}$,
	respectively.
In fact,
	this statement immediately follows from  \refEq{3.1_PL_DefFuncW} and the fact that $\mpWT{n}$ is invariant under $\gpSym{n}$,
	i.e.,
	\envM
	{
		\racFA{\mpWT{n}}{\sig}{\lbfL_n}
	}
	{
		\mpWTT{n}{\lbfL_n}
	}
	for any $\sig\in\gpSym{n}$.
Note that $\Gp{(12),(123)}$ and $\Gp{(23),(234)}$ are equivalent to the symmetric groups on $\Set{1,2,3}$ and $\Set{2,3,4}$,
	respectively.

We are now able to prove \refLem{3.1_Lem4}.
\envProof[\refLem{3.1_Lem4}]{
We can obtain \refEq{3.1_Lem4i_Eq1} by using \refEq{3.1_Lem5i_Eq1}
	because of the invariance of $\empWT{4}{2,2}$ under $\Gp{(12),(34)}$.
Similarly,
	we can obtain 
	the equations from \refEq{3.1_Lem4i_Eq2} through \refEq{3.1_Lem4ii_Eq7}
	by using 
	the congruence equations from \refEq{3.1_Lem5i_Eq2} through \refEq{3.1_Lem5ii_Eq7},
	respectively.
}

\subsection{Renormalization relations}\label{sectThreeTwo}
We define a characteristic function $\fcMChSoO{n}$ of the set $\setN^n$ by
	\envHLineDef
	{
		\fcMChSo{n}{\lbfL_n}
	}
	{
		\envCaseTPd{
			1			&	(l_1=\cdots=l_n=1)
		}{
			0			&	(otherwise)
		}
	}
Note that $\fcMChSoO{n}=(\fcMChSoO{1})^{\otimes n}$ and $\fcCVo[n]=\fcMChSO{n}+\fcMChSoO{n}$. 
For any real-valued functions $f_1,\ldots,f_j$ of $n$ variables,
	we define the product $f_1\cdots f_j$ of the functions 
	by using the multiplication in the real number field such that 
	\envHLineDefPd
	{
		\Fcr{f_1 \cdots f_j}{x_1,\ldots,x_n}
	}
	{
		\Fc{f_1}{x_1,\ldots,x_n}\cdots\Fc{f_j}{x_1,\ldots,x_n}
	}
For example,
	\envMThPd{
		\Fc{ \fcMChSoO{2} \cdot \fcMZO{1} \circ \mpWT{2} }{\lbfL_2}
	}
	{
		\fcMChSo{2}{\lbfL_2} \Fc{\fcMZO{1} \circ \mpWT{2}}{\lbfL_2}
	}
	{
		\fcMChSo{2}{l_1,l_2} \fcMZ{1}{l_1+l_2}
	}

The renormalization relations for depths less than $5$ are written in terms of real-valued functions, as follows. 
\begin{proposition}\label{3.2_Prop1}
We have 
	\envHLineCENmePd
	{\label{3.2_Prop1_EqFuncDp1}
		\fcMZHO{1}
	}
	{
		\fcMZSO{1}
	}
	{\label{3.2_Prop1_EqFuncDp2}\rule{0cm}{20pt}
		\fcMZHO{2}
	}
	{
		\fcMZSO{2} - \opF{1}{2} \fcMChSoO{2} \cdot \fcMZO{1} \circ \mpWT{2}
	}
	{\label{3.2_Prop1_EqFuncDp3}\rule{0cm}{20pt}
		\fcMZHO{3}
	}
	{
		\fcMZSO{3} 
		- 
		\opF{1}{2} (\fcMChSoO{2} \cdot \fcMZO{1} \circ \mpWT{2}) \otimes \fcMZSO{1} 
		+
		\opF{1}{3} \fcMChSoO{3} \cdot \fcMZO{1} \circ \mpWT{3}
	}
	{\label{3.2_Prop1_EqFuncDp4}\rule{0cm}{20pt}
		\fcMZHO{4}
	}
	{
		\fcMZSO{4}
		- 
		\opF{1}{2} (\fcMChSoO{2} \cdot \fcMZO{1} \circ \mpWT{2}) \otimes \fcMZSO{2} 
		+
		\opF{1}{3} (\fcMChSoO{3} \cdot \fcMZO{1} \circ \mpWT{3}) \otimes \fcMZSO{1} 
		+
		\opF{1}{16}\fcMChSoO{4} \cdot \fcMZO{1} \circ \mpWT{4}
	}
\end{proposition}

We require two lemmas to prove \refProp{3.2_Prop1}. 

\begin{lemma}\label{3.2_Lem1}
Let $\Fc{P}{T}=\SmT{j=0}{n}a_jT^j$ be a polynomial whose degree $n$ is less than $5$.
Then the constant term of $\pmpRH{\Fc{P}{T}} - \Fc{P}{T}$ is 
	\envHLine[3.2_Lem1_Eq]
	{\hspace{-15pt}
		\pmpRH{\Fc{P}{T}} \mVert_{T=0} - \Fc{P}{0} 
	}
	{
		\envCaseFPd[d]{
			0
			&
			(n<2)
		}{
			a_2\fcMZ{1}{2}
			&
			(n=2)
		}{
			a_2\fcMZ{1}{2} - 2a_3\fcMZ{1}{3}
			&
			(n=3)
		}{	
			a_2\fcMZ{1}{2} - 2a_3\fcMZ{1}{3} + \opF{27}{2}a_4\fcMZ{1}{4}  
			&
			(n=4)
		}
	}
\end{lemma}
\begin{lemma}\label{3.2_Lem2}
Let $n$ be an integer with $1\leq n\leq4$,
	and let $\lbfL_n=(l_1,\ldots,l_n)\in\setN^n$.
Then
	\envHLineCECmNmePt{\approx}
	{\label{3.2_Lem2_Eq1}
		\mpEVHi{\lbfL_1}{T}
	}
	{
		0
	}
	{\label{3.2_Lem2_Eq2}\rule{0pt}{20pt}
		\mpEVHi{\lbfL_2}{T}
	}
	{
		\opF{1}{2} \fcMChSo{2}{\lbfL_2} T^2
	}
	{\label{3.2_Lem2_Eq3}\rule{0pt}{20pt}
		\mpEVHi{\lbfL_3}{T}
	}
	{
		\opF{1}{2} \Fc{ \fcMChSoO{2} \otimes \fcMZHO{1} }{ \lbfL_3 } T^2
		+ 
		\opF{1}{6} \fcMChSo{3}{\lbfL_3}T^3
	}
	{\label{3.2_Lem2_Eq4}\rule{0pt}{20pt}
		\mpEVHi{\lbfL_4}{T}
	}
	{
		\opF{1}{2} \Fc{\fcMChSoO{2}\otimes\fcMZHO{2}}{\lbfL_4}T^2 + \opF{1}{6} \Fc{\fcMChSoO{3}\otimes\fcMZHO{1}}{\lbfL_4}T^3 + \opF{1}{24} \fcMChSo{4}{\lbfL_4}T^4
	}
	where $\approx$ means the congruence relation on $\setR[T]$ modulo $\setR T + \setR$,
	i.e.,
	$\Fc{P}{T}\approx\Fc{Q}{T}$ if and only if $\degg[n]{\Fc{P}{T}-\Fc{Q}{T}}<2$.
\end{lemma}

We will now prove \refProp{3.2_Prop1}.
We will then discuss the proofs of \refLem[s]{3.2_Lem1} and \ref{3.2_Lem2}.

\envProof[\refProp{3.2_Prop1}]{
We first introduce the following identity 
	for proving \refEq{3.2_Prop1_EqFuncDp1}, \refEq{3.2_Prop1_EqFuncDp2}, \refEq{3.2_Prop1_EqFuncDp3}, and \refEq{3.2_Prop1_EqFuncDp4}:
	for any index set $\lbfL_n=(l_1,\ldots,l_n)$ in $\setN^n$,
	\envHLinePd[3.2_Prop1Pr_Eq0]
	{
		\fcMChSo{n}{\lbfL_n} \fcMZ{1}{n}
	}
	{
		\Fc{ \fcMChSoO{n} \cdot \fcMZO{1} \circ \mpWT{n}}{\lbfL_n}
	}
This identity is easily obtained by the definition of $\fcMChSo[]{n}{}$
	and 
	the fact that
	\envMTh{
		\fcMZ{1}{n}
	}{
		\Fc{ \fcMZO{1}}{ \mpWTT{n}{\lbfL_n} }
	}{
		\Fc{ \fcMZO{1} \circ \mpWT{n} }{\lbfL_n}
	}
	if $l_1=\cdots=l_n=1$.

It follows from \refEq{3.2_Lem1_Eq} and \refEq{3.2_Lem2_Eq1} that
	\envHLinePd
	{
		\pmpRH{\mpEVHi{\lbfL_1}{T}} \mVert_{T=0} - \mpEVHi{\lbfL_1}{0}
	}
	{
		0
	}
Using \refEq{2_PL_DefRMZV} and \refEq{3.2_PL_EqFundRelRegIKZ} with $T=0$,
	we can restate this identity as
	\envHLineCm
	{
		\fcMZS{1}{\lbfL_1} - \fcMZH{1}{\lbfL_1}
	}
	{
		0
	}
	which proves \refEq{3.2_Prop1_EqFuncDp1}.
Similarly,
	we obtain from \refEq{3.2_Lem1_Eq} and \refEq{3.2_Lem2_Eq2} that
	\envHLineCm
	{
		\fcMZS{2}{\lbfL_2} - \fcMZH{2}{\lbfL_2}
	}
	{
		\opF{1}{2} \fcMChSo{2}{\lbfL_2} \fcMZ{1}{2}
	}
	which proves \refEq{3.2_Prop1_EqFuncDp2}
	since
	\envM
	{
		\fcMChSo{2}{\lbfL_2} \fcMZ{1}{2}
	}
	{
		\Fc{\fcMChSoO{2} \cdot \fcMZO{1} \circ \mpWT{2}}{\lbfL_2}
	}
	by \refEq{3.2_Prop1Pr_Eq0}.
We can obtain from \refEq{3.2_Lem1_Eq} and \refEq{3.2_Lem2_Eq3} that
	\envHLineCm
	{
		\fcMZS{3}{\lbfL_3} - \fcMZH{3}{\lbfL_3}
	}
	{
		\opF{1}{2} \Fc{ \fcMChSoO{2} \otimes \fcMZHO{1} }{ \lbfL_3 }  \fcMZ{1}{2} - \opF{1}{3}\fcMChSo{3}{\lbfL_3} \fcMZ{1}{3} 
	}
	which proves \refEq{3.2_Prop1_EqFuncDp3}
	since
	\envM{
		\fcMChSo{3}{\lbfL_3} \fcMZ{1}{3} 
	}{
		\Fc{\fcMChSoO{3} \cdot \fcMZO{1} \circ \mpWT{3}}{\lbfL_3}
	}
	and
	\envHLineF
	{
		\Fc{ \fcMChSoO{2} \otimes \fcMZHO{1} }{ \lbfL_3 }  \fcMZ{1}{2}
	}
	{
		\fcMChSo{2}{l_1,l_2} \fcMZ{1}{2} \fcMZH{1}{l_3}
	}
	{
		(\Fc{\fcMChSoO{2} \cdot \fcMZO{1} \circ \mpWT{2} }{l_1,l_2}) \fcMZS{1}{l_3}
	}
	{
		\Fc{ (\fcMChSoO{2} \cdot \fcMZO{1} \circ \mpWT{2}) \otimes \fcMZSO{1} }{\lbfL_3}
	}
	by \refEq{3.2_Prop1_EqFuncDp1} and \refEq{3.2_Prop1Pr_Eq0}.
We can obtain from \refEq{3.2_Lem1_Eq} and \refEq{3.2_Lem2_Eq4} that
	\envHLinePd[3.2_Prop1Pr_EqDp4_1]
	{
		\fcMZS{4}{\lbfL_4} - \fcMZH{4}{\lbfL_4}
	}
	{
		\opF{1}{2} \Fc{\fcMChSoO{2}\otimes\fcMZHO{2}}{\lbfL_4} \fcMZ{1}{2}
		-
		\opF{1}{3} \Fc{\fcMChSoO{3}\otimes\fcMZHO{1}}{\lbfL_4} \fcMZ{1}{3}
		+ 
		\opF{9}{16} \fcMChSo{4}{\lbfL_4} \fcMZ{1}{4}
	}
The first term on the right-hand side of \refEq{3.2_Prop1Pr_EqDp4_1} can be calculated as
	\envHLinePd[3.2_Prop1Pr_EqDp4_2]
	{
		\opF{1}{2} \Fc{\fcMChSoO{2}\otimes\fcMZHO{2}}{\lbfL_4} \fcMZ{1}{2}
	}
	{
		\opF{1}{2} \Fc{\fcMChSoO{2}\otimes\fcMZSO{2}}{\lbfL_4} \fcMZ{1}{2} - \opF{5}{8} \fcMChSo{4}{\lbfL_4}\fcMZ{1}{4}
	}
In fact,
	we see from \refEq{3.2_Prop1_EqFuncDp2} and \refEq{3.2_Prop1Pr_Eq0} that
	\envMCm
	{
		\fcMZH{2}{l_3,l_4}
	}
	{
		\fcMZS{2}{l_3,l_4} - \opF{1}{2} \fcMChSo{2}{l_3,l_4} \fcMZ{1}{2}
	}
	and so
	\envHLineFCm
	{
		\opF{1}{2} \Fc{\fcMChSoO{2}\otimes\fcMZHO{2}}{\lbfL_4} \fcMZ{1}{2}
	}
	{
		\opF{1}{2} \fcMChSo{2}{l_1,l_2} \fcMZH{2}{l_3,l_4} \fcMZ{1}{2}
	}
	{
		\opF{1}{2} \fcMChSo{2}{l_1,l_2} \fcMZS{2}{l_3,l_4} \fcMZ{1}{2} - \opF{1}{4} \fcMChSo{2}{l_1,l_2} \fcMChSo{2}{l_3,l_4} \fcMZ{1}{2}^2
	}
	{
		\opF{1}{2} \Fc{\fcMChSoO{2}\otimes\fcMZSO{2}}{\lbfL_4} \fcMZ{1}{2} - \opF{1}{4} \fcMChSo{4}{\lbfL_4} \fcMZ{1}{2}^2
	}
	where we note that 
	\envMPd
	{
		\fcMChSoO{4}
	}
	{
		(\fcMChSoO{2})^{\otimes2}
	}
This proves \refEq{3.2_Prop1Pr_EqDp4_2}
	because
	\envHLineCm[3.2_Prop1Pr_EqDp4_3]
	{
		\fcMZ{1}{2}^2
	}
	{
		\opF{5}{2} \fcMZ{1}{4}
	}
	which follows from Euler's results $\fcMZ{1}{2}=\opF[s]{\pi^2}{6}$ and $\fcMZ{1}{4}=\opF[s]{\pi^2}{90}$.
Then,
	noting \refEq{3.2_Prop1_EqFuncDp1}, \refEq{3.2_Prop1Pr_Eq0}, and \refEq{3.2_Prop1Pr_EqDp4_2},
	we can calculate \refEq{3.2_Prop1Pr_EqDp4_1} as
	\envLineThCm
	{
		\fcMZS{4}{\lbfL_4} - \fcMZH{4}{\lbfL_4}
	}
	{
		\opF{1}{2} \Fc{\fcMChSoO{2}\otimes\fcMZSO{2}}{\lbfL_4} \fcMZ{1}{2}
		-
		\opF{1}{3} \Fc{\fcMChSoO{3}\otimes\fcMZSO{1}}{\lbfL_4} \fcMZ{1}{3}
		-
		\opF{1}{16} \fcMChSo{4}{\lbfL_4} \fcMZ{1}{4}
	}
	{
		\opF{1}{2} \Fc{ (\fcMChSoO{2} \cdot \fcMZO{1} \circ \mpWT{2}) \otimes \fcMZSO{2} }{\lbfL_4}
		-
		\opF{1}{3} \Fc{ (\fcMChSoO{3} \cdot \fcMZO{1} \circ \mpWT{3}) \otimes \fcMZSO{1} }{\lbfL_4}
		-
		\opF{1}{16} \Fc{ \fcMChSoO{4} \cdot \fcMZO{1} \circ \mpWT{4} }{\lbfL_4}
	}
	which proves \refEq{3.2_Prop1_EqFuncDp4}.
}

We will now show \refLem[s]{3.2_Lem1} and \ref{3.2_Lem2} for the completeness of the proof of \refProp{3.2_Prop1}.

\envProof[\refLem{3.2_Lem1}]{
Let $\sbL[]{}$ denote the Landau symbol.
By definition,
	\envOTLine
	{
		\pfcA{u}
	}
	{
		\SmT{k=0}{\infty} \gam_k u^k
	}
	{
		\exx[G]{\SmT{m=2}{\infty} \opF{\mo^m\fcMZ{1}{m}}{m} u^m }
	}
	near $u=0$.
Thus,
	\envHLineThCm
	{
		\pfcA{u}
	}
	{
		1 
		+ 
		\bkR[g]{ \opF{\fcMZ{1}{2}}{2}u^2 - \opF{\fcMZ{1}{3}}{3}u^3 + \opF{\fcMZ{1}{4}}{4}u^4 + \sbL{u^5} }
		+
		\opF{1}{2}\bkR[g]{ \opF{\fcMZ{1}{2}}{2}u^2 + \sbL{u^3} }^2
		+
		\cdots
	}
	{
		1 + \opF{\fcMZ{1}{2}}{2}u^2 - \opF{\fcMZ{1}{3}}{3}u^3 + \bkR[a]{  \opF{\fcMZ{1}{4}}{4} +  \opF{\fcMZ{1}{2}^2}{8} } u^4 + \sbL{u^5} 
	}
	and so 
	\envPLine{\text{
	\envMCm{
		\gam_0
	}{
		1
	}\quad
	\envMCm{
		\gam_1
	}{
		0
	}\quad
	\envMCm{
		\gam_2
	}{
		\opF{\fcMZ{1}{2}}{2}
	}\quad
	\envMCm{
		\gam_3
	}{
		- \opF{\fcMZ{1}{3}}{3}
	}\quad
	and
	\quad\envMThCm{
		\gam_4
	}{
		\opF{2\fcMZ{1}{4}+\fcMZ{1}{2}^2}{8}
	}{
		\opF{9\fcMZ{1}{4}}{16}
	}
	}}
	where we have used \refEq{3.2_Prop1Pr_EqDp4_3} for the last equality.
Therefore,
	we  see from \refEq{2_PL_DefMapRhoDetail} that
	\envPLine{\HLineCSCm[p]
	{
		\pmpRH{1}
	}
	{
		1
	}
	{
		\pmpRH{T}
	}
	{
		T
	}
	{
		\pmpRH{T^2}
	}
	{
		T^2+\fcMZ{1}{2}
	}\HLineCFCm[*]
	{
		\pmpRH{T^3}
	}
	{
		T^3+3\fcMZ{1}{2}T-2\fcMZ{1}{3}
	}
	{
		\pmpRH{T^4}
	}
	{
		T^4+6\fcMZ{1}{2}T^2-8\fcMZ{1}{3}T+ \opF{27}{2}\fcMZ{1}{4} 
	}
	}
	which proves \refEq{3.2_Lem1_Eq} since 
	\envMPd{
		\pmpRH{\Fc{P}{T}} \mVert_{T=0} - \Fc{P}{0}
	}
	{
		\SmT{j=2}{n}a_j \pmpRH{T^j} \mVert_{T=0}
	}
}

In order to prove \refLem{3.2_Lem2},
	we require the following detailed results on polynomials $\mpEVHi{\lbfL_n}{T}$ and values $\fcMZH{n}{\lbfL_n} $ and $\fcMZS{n}{\lbfL_n}$ for smaller depths. 

\begin{lemma}\label{3.2_Lem3}
\mbox{}\\\bfTx{(i)}
Let $\lbfL_n=(l_1,\ldots,l_n)\in\setN^n$.
Then
	\envHLineCSNmePd
	{\label{3.2_Lem3iEq1}
		\mpEVHi{\lbfL_1}{T}
	}
	{
		\fcMChSo{1}{\lbfL_1}T + \fcMZH{1}{\lbfL_1} 
	}
	{\label{3.2_Lem3iEq2}\rule{0pt}{20pt}
		\mpEVHi{\lbfL_2}{T}
	}
	{
		\opF{1}{2} \fcMChSo{2}{\lbfL_2}T^2 + \Fc{ \fcMChSoO{1} \otimes \fcMZHO{1} }{\lbfL_2}T  + \fcMZH{2}{\lbfL_2} 
	}
	{\label{3.2_Lem3iEq3}\rule{0pt}{20pt}
		\mpEVHi{\lbfL_3}{T}
	}
	{
		\opF{1}{6}\fcMChSo{3}{\lbfL_3}T^3  + \opF{1}{2} \Fc{ \fcMChSoO{2} \otimes \fcMZHO{1} }{\lbfL_3}T^2 + \Fc{\fcMChSoO{1}\otimes\fcMZHO{2}}{\lbfL_3}T + \fcMZH{3}{\lbfL_3} 
	}	
\mbox{}\\\bfTx{(ii)}
Suppose that $l_2>1$,
	and let $\letB\in\Set{*,\sh}$.
Then 
	\envHLineCSNme
	{\label{3.2_Lem3iiEq1}
		\fcMZB{1}{1} 
	}
	{
		0
	,}
	{\label{3.2_Lem3iiEq2}
		\fcMZB{2}{1,l_2} 
	}
	{
		-\bkR{ \fcMZ{2}{l_2,1} + \fcMZ{1}{l_2+1} }
	,}
	{\label{3.2_Lem3iiEq3}\rule{0pt}{30pt}
		\fcMZB{2}{1,1}
	}
	{
		\envCaseTPd[d]{
			-\opF{1}{2} \fcMZ{1}{2}
			&
			(\letB=*)	
		}{
			0
			&
			(\letB=\sh)	
		}
	}
\end{lemma}

\envProof{
We recall that $\mpEVH[]{}:\alHNh[y](\simeq\alHRh)\to\setR[T]$ 
	is a homomorphism 
	characterized by the properties that it extends the evaluation map $\mpEVN[]{}:\alHN\to\setR$ and sends $y=z_1$ to $T$,
	and that
	$\mpEVHi{\lbfL_n}{T}$ denotes $\mpEVH{z_{l_1}\cdots z_{l_n}}$.
We note that
	$\mpEVHi{\lbfL_n}{T}=\mpEVNi{\lbfL_n}{T}=\fcMZ{n}{\lbfL_n}=\fcMZH{n}{\lbfL_n}$ if $l_1>1$,
	or
	\envHLinePd[3.2_Lem3zPr_Eq]
	{
		\mpEVHi{\lbfL_n}{T}
	}
	{
		\fcMZH{n}{\lbfL_n}
		\qquad
		(l_1>1)
	}

We will first prove \refEq{3.2_Lem3iEq1} and \refEq{3.2_Lem3iiEq1}.
It follows from the definition of $\mpEVH[]{}$ that
	\envMThPd
	{
		\mpEVHi{1}{T}
	}
	{
		\mpEVH{z_1}
	}
	{
		T
	}
Since
	\envM
	{
		\mpEVHi{1}{0}
	}
	{
		\fcMZH{1}{1}
	}
	by definition,
	we have 
	\envHLineCm[3.2_Lem3aPr_Eq2]
	{
		\fcMZH{1}{1} 
	}
	{
		0
	}
	and
	\envHLinePd[3.2_Lem3aPr_Eq0]
	{
		\mpEVHi{1}{T}
	}
	{
		T + \fcMZH{1}{1}
	}
We easily see from the definition of $\fcMChSo[]{1}{}$ that
	$\fcMChSo{1}{l_1} T$ is $T$ if $l_1=1$, 
	and $0$ otherwise.	
Thus,
	we obtain from \refEq{3.2_Lem3zPr_Eq} for $n=1$ and \refEq{3.2_Lem3aPr_Eq0} that 
	\envHLineCm[3.2_Lem3aPr_Eq1]
	{
		\mpEVHi{l_1}{T}
	}
	{
		\fcMChSo{1}{l_1} T + \fcMZH{1}{l_1} 
		\qquad
		(l_1\in\setN)
	}
	which proves \refEq{3.2_Lem3iEq1}.
It follows from \refEq{3.2_Lem1_Eq} for $n<2$ and \refEq{3.2_Lem3aPr_Eq0} that
	\envMCm
	{
		\pmpRH{\mpEVHi{1}{T}} \mVert_{T=0} - \mpEVHi{1}{0}
	}
	{
		0
	}
	and so
	we see from \refEq{3.2_PL_EqFundRelRegIKZ} that
	\envMCm
	{
		\fcMZS{1}{1} 
	}
	{
		\fcMZH{1}{1} 
	}
	which together with \refEq{3.2_Lem3aPr_Eq2} proves \refEq{3.2_Lem3iiEq1}.

We will next prove \refEq{3.2_Lem3iEq2}, \refEq{3.2_Lem3iiEq2}, and \refEq{3.2_Lem3iiEq3}.
By \refEq{3.1_Lem1_Eq} with $l_1=1$ and $l_1=l_2=1$,
	we obtain
	\envHLineCFLaaCm 
	{
		z_1z_{l_2}
	}
	{
		z_{l_2}*z_1-z_{l_2}z_1- z_{l_2+1}
	}
	{
		z_1z_1
	}
	{
		\opF{1}{2}z_1*z_1-\opF{1}{2}z_2
	}
	respectively. 
Since
	\envMCm
	{
		\mpEVHi{1,l_2}{0}
	}
	{
		\fcMZH{2}{1,l_2}
	}
	applying  $\mpEVH[]{}$ to both sides of these equations yields
	\envHLine[3.2_Lem3bPr_Eq2]
	{
		\fcMZH{2}{1,l_2} 
	}
	{
		\envCaseTCm[d]{
			-(\fcMZ{2}{l_2,1}+\fcMZ{1}{l_2+1})
			&
			(l_2>1)
		}{
			-\opF{1}{2}\fcMZ{1}{2}
			&
			(l_2=1)
		}
	}
	and
	\envHLine[3.2_Lem3bPr_Eq0]
	{
		\mpEVHi{1,l_2}{T}
	}
	{
		\envCaseTPd[d]{
			\fcMZ{1}{l_2}T + \fcMZH{2}{1,l_2} 
			&
			(l_2>1)
		}{
			\opF{1}{2}T^2 + \fcMZH{2}{1,1} 
			&
			(l_2=1)
		}
	}
We see that
	$\opF[s]{\fcMChSo{2}{l_1,l_2}T^2}{2} $ is $\opF[s]{T^2}{2}$ if $l_1=l_2=1$,
	and $0$ otherwise,
	and that
	$\fcMChSo{1}{l_1}\fcMZH{1}{l_2}T$ is $\fcMZ{1}{l_2}T$ if $l_1=1$ and $l_2>1$,
	and $0$ otherwise,
	where we have used \refEq{3.2_Lem3aPr_Eq2} for the latter statement.
Thus,
	we obtain from \refEq{3.2_Lem3zPr_Eq} for $n=2$ and \refEq{3.2_Lem3bPr_Eq0} that	
	\envHLineCm[3.2_Lem3bPr_Eq1]
	{
		\mpEVHi{\lbfL_2}{T}
	}
	{
		\opF{\fcMChSo{2}{l_1,l_2}}{2} T^2 + \fcMChSo{1}{l_1}\fcMZH{1}{l_2}  T + \fcMZH{2}{l_1,l_2} 
		\qquad
		(\lbfL_2\in\setN^2)
	}
	which proves \refEq{3.2_Lem3iEq2}.
It follows from \refEq{3.2_Lem1_Eq} for $n\leq2$ and \refEq{3.2_Lem3bPr_Eq0}
	that
	$\pmpRH{\mpEVHi{1,l_2}{T}} \mVert_{T=0} - \mpEVHi{1,l_2}{0}$ is $0$ if $l_2>1$,
	and
	$\opF[s]{\fcMZ{1}{2}}{2}$ if $l_2=1$.
Therefore
	we see from \refEq{3.2_PL_EqFundRelRegIKZ} that
	\envHLine
	{
		\fcMZS{1}{1,l_2} 
	}
	{
		\fcMZH{1}{1,l_2} 
		+
		\envCaseTCm[d]{
			0
			&
			(l_2>1)
		}{
			\opF{1}{2}\fcMZ{1}{2}
			&
			(l_2=1)
		}
	}
	which together with \refEq{3.2_Lem3bPr_Eq2} proves \refEq{3.2_Lem3iiEq2} and \refEq{3.2_Lem3iiEq3}.
		
We now prove \refEq{3.2_Lem3iEq3}.
Replacing $(l_1,l_2,l_3)$ of \refEq{3.1_Lem2_Eq1} with $(l_2,l_3,1)$, 
	we obtain
	\envHLinePd[3.2_Lem3cPr_EqHarRel1]
	{
		z_1z_{l_2}z_{l_3}
	}
	{
		z_{l_2}z_{l_3}*z_1 - z_{l_2}z_{l_3}z_1 - z_{l_2}z_{1}z_{l_3} - z_{l_2+1}z_{l_3} - z_{l_2}z_{l_3+1}
	}
Substituting $l_2=1$ and $l_2=l_3=1$ in \refEq{3.2_Lem3cPr_EqHarRel1} yields
	\envHLineCFCmNme
	{\label{3.2_Lem3cPr_EqHarRel2}
		z_1z_1z_{l_3}
	}
	{
		\opF{1}{2} z_1z_{l_3}*z_1 - \opF{1}{2} (z_1z_{l_3}z_1  + z_{2}z_{l_3} + z_1z_{l_3+1})
	}
	{\label{3.2_Lem3cPr_EqHarRel3}
		z_1z_1z_1
	}
	{
		\opF{1}{3}z_1z_1*z_1 - \opF{1}{3} (z_2z_1+z_1z_2)
	}
	respectively.
Let $\sim$ mean the congruence relation on $\setR[T]$ modulo $\setR$,
	i.e.,
	$\Fc{P}{T}\sim\Fc{Q}{T}$ if and only if $\Fc{P}{T}-\Fc{Q}{T}\in\setR$.
Note that
	\envMPt{\sim}
	{
		\mpEVHi{\lbfL_n}{T}
	}
	{
		0
	}
	if $l_1>1$.
Applying  $\mpEVH[]{}$ to both sides of \refEq{3.2_Lem3cPr_EqHarRel1},
	we obtain
	\envHLinePdPt[3.2_Lem3cPr_Eq0help1]{\sim}
	{
		\mpEVHi{1,l_2,l_3}{T}
	}
	{
		\fcMZ{2}{l_2,l_3}T
		\qquad
		(l_2>1)
	}
Because of \refEq{3.2_Lem3iEq2} and \refEq{3.2_Lem3cPr_Eq0help1},
	applying  $\mpEVH[]{}$ to both sides of \refEq{3.2_Lem3cPr_EqHarRel2} yields
	\envHLineFPt{\sim}
	{
		\mpEVHi{1,1,l_3}{T}
	}
	{
		\opF{1}{2} \mpEVHi{1,l_3}{T} T - \opF{1}{2} (\mpEVHi{1,l_3,1}{T} + \mpEVHi{1,l_3+1}{T})
	}
	{
		\opF{1}{2} \bkR{ \fcMZH{1}{l_3} T + \fcMZH{2}{1,l_3} } T 
		- 
		\opF{1}{2} (\fcMZ{2}{l_3,1}T + \fcMZ{1}{l_3+1}  T  )
	}
	{
		\opF{1}{2} \fcMZH{1}{l_3} T^2  + \opF{1}{2} \bkR{ \fcMZH{2}{1,l_3} - \fcMZ{2}{l_3,1} - \fcMZ{1}{l_3+1} }T
	}
	for $l_3>1$,
	which together with \refEq{3.2_Lem3iiEq2} gives 
	\envHLinePdPt[3.2_Lem3cPr_Eq0help2]{\sim}
	{
		\mpEVHi{1,1,l_3}{T}
	}
	{
		\opF{1}{2} \fcMZH{1}{l_3} T^2  + \fcMZH{2}{1,l_3} T
		\qquad
		(l_3>1)
	}
Applying  $\mpEVH[]{}$ to both sides of \refEq{3.2_Lem3cPr_EqHarRel3} also yields
	\envHLineFCmPt{\sim}
	{
		\mpEVHi{1,1,1}{T}
	}
	{
		\opF{1}{3}\mpEVHi{1,1}{T} T - \opF{1}{3}\mpEVHi{1,2}{T}
	}
	{
		\opF{1}{3} \bkR[a]{ \opF{1}{2}T^2 + \fcMZH{2}{1,1} } T - \opF{1}{3} \fcMZH{1}{2}T
	}
	{
		\opF{1}{6} T^3  + \opF{1}{3} \bkR{ \fcMZH{2}{1,1} - \fcMZ{1}{2} }T
	}
	which together with \refEq{3.2_Lem3iiEq3} gives 
	\envHLinePdPt[3.2_Lem3cPr_Eq0help3]{\sim}
	{
		\mpEVHi{1,1,1}{T}
	}
	{
		\opF{1}{6} T^3  + \fcMZH{2}{1,1}T 
	}
Similarly to \refEq{3.2_Lem3aPr_Eq1} and \refEq{3.2_Lem3bPr_Eq1},
	we can obtain from \refEq{3.2_Lem3zPr_Eq} for $n=3$,
	\refEq{3.2_Lem3cPr_Eq0help1}, \refEq{3.2_Lem3cPr_Eq0help2}, and \refEq{3.2_Lem3cPr_Eq0help3} that
	\envHLineCmPt{\sim}
	{
		\mpEVHi{\lbfL_3}{T}
	}
	{
		\opF{\fcMChSo{3}{l_1,l_2,l_3}}{6} T^3 
		+ 
		\opF{\fcMChSo{2}{l_1,l_2}}{2} \fcMZH{1}{l_3}T^2 
		+ 
		\fcMChSo{1}{l_1}\fcMZH{2}{l_2,l_3} T
		\qquad
		(\lbfL_3\in\setN^3)
	}
	which proves \refEq{3.2_Lem3iEq3},
	since
	\envMPd
	{
		\mpEVHi{\lbfL_3}{0}
	}
	{
		\fcMZH{3}{\lbfL_3}
	}
}

We are now able to prove \refLem{3.2_Lem2}.
\envProof[\refLem{3.2_Lem2}]{
Identities \refEq{3.2_Lem2_Eq1}, \refEq{3.2_Lem2_Eq2}, and \refEq{3.2_Lem2_Eq3} immediately follow from \refEq{3.2_Lem3iEq1}, \refEq{3.2_Lem3iEq2}, and \refEq{3.2_Lem3iEq3},
	respectively.
	
We now prove \refEq{3.2_Lem2_Eq4}.
Replacing $(l_1,l_2,l_3,l_4)$ of \refEq{3.1_Lem3_Eq1} with $(l_2,l_3,l_4,1)$, 	
	we obtain
	\envHLinePd[3.2_Lem2Pr_EqHarRel1]
	{
		z_1z_{l_2}z_{l_3}z_{l_4}
	}
	{
		z_{l_2}z_{l_3}z_{l_4}*z_1 - z_{l_2}z_{l_3}z_{l_4}z_1 - z_{l_2}z_{l_3}z_1z_{l_4} - z_{l_2}z_1z_{l_3}z_{l_4}
		\lnAH
		-
		z_{l_2+1}z_{l_3}z_{l_4} - z_{l_2}z_{l_3+1}z_{l_4} - z_{l_2}z_{l_3}z_{l_4+1}
	}
Substituting $l_2=1$, $l_2=l_3=1$, and $l_1=l_2=l_3=1$ in \refEq{3.2_Lem2Pr_EqHarRel1} yields
	\envHLineCSCmNme
	{\label{3.2_Lem2Pr_EqHarRel2}
		z_1z_1z_{l_3}z_{l_4}
	}
	{
		\opF{1}{2} z_1z_{l_3}z_{l_4}*z_1 
		\lnAH
		- 
		\opF{1}{2}( z_1z_{l_3}z_{l_4}z_1 + z_1z_{l_3}z_1z_{l_4} + z_1z_{l_3+1}z_{l_4} + z_1z_{l_3}z_{l_4+1} + z_2z_{l_3}z_{l_4} )
	}
	{\label{3.2_Lem2Pr_EqHarRel3}
		z_1z_1z_1z_{l_4}
	}
	{
		\opF{1}{3} z_1z_1z_{l_4}*z_1 
		- 
		\opF{1}{3}( z_1z_1z_{l_4}z_1 + z_1z_2z_{l_4} + z_1z_1z_{l_4+1} + z_2z_1z_{l_4} )
	}
	{\label{3.2_Lem2Pr_EqHarRel4}
		z_1z_1z_1z_1
	}
	{
		\opF{1}{4} z_1z_1z_1*z_1 
		- 
		\opF{1}{4}( z_1z_2z_1 + z_1z_1z_2 + z_2z_1z_1 )
	}
	respectively.
Applying  $\mpEVH[]{}$ to both sides of \refEq{3.2_Lem2Pr_EqHarRel1},
	we obtain 
	\envHLinePdPt[3.2_Lem2Pr_Eq0help1]{\approx}
	{
		\mpEVHi{1,l_2,l_3,l_4}{T}
	}
	{
		0
		\qquad
		(l_2>1)
	}
Because of \refEq{3.2_Lem3iEq3} and \refEq{3.2_Lem2Pr_Eq0help1},
	applying  $\mpEVH[]{}$ to both sides of \refEq{3.2_Lem2Pr_EqHarRel2} yields
	\envHLineFPt{\approx}
	{
		\mpEVHi{1,1,l_3,l_4}{T}
	}
	{
		\opF{1}{2}\mpEVHi{1,l_3,l_4}{T}T
	}
	{
		\opF{1}{2}\bkR{ \fcMZH{2}{l_3,l_4} T + \fcMZH{3}{1,l_3,l_4}  }T
	}
	{
		\opF{1}{2} \fcMZH{2}{l_3,l_4} T^2
	}
	for $l_3>1$,
	which is summarized as
	\envHLinePdPt[3.2_Lem2Pr_Eq0help2]{\approx}
	{
		\mpEVHi{1,1,l_3,l_4}{T}
	}
	{
		\opF{1}{2}\fcMZH{2}{l_3,l_4} T^2
		\qquad
		(l_3>1)
	}
Because of \refEq{3.2_Lem3iEq3} and \refEq{3.2_Lem2Pr_Eq0help2},
	applying $\mpEVH[]{}$ to both sides of  \refEq{3.2_Lem2Pr_EqHarRel3} also yields
	\envHLineFPt{\approx}
	{
		\mpEVHi{1,1,1,l_4}{T}
	}
	{
		\opF{1}{3} \mpEVHi{1,1,l_4}{T}T  - \opF{1}{3} \bkR{ \mpEVHi{1,1,l_4,1}{T} + \mpEVHi{1,1,l_4+1}{T} }
	}
	{
		\opF{1}{3}\bkR[a]{ \opF{1}{2} \fcMZH{1}{l_4}T^2 + \fcMZH{2}{1,l_4} T + \fcMZH{3}{1,1,l_4}  }T - \opF{1}{6} \bkR{ \fcMZH{2}{l_4,1} + \fcMZH{1}{l_4+1} }T^2 
	}
	{
		\opF{1}{6} \fcMZH{1}{l_4}T^3 + \bkR[a]{ \opF{1}{3} \fcMZH{2}{1,l_4} - \opF{1}{6} \bkR{ \fcMZH{2}{l_4,1} + \fcMZH{1}{l_4+1} } }T^2 
	}
	for $l_4>1$,
	which together with \refEq{3.2_Lem3iiEq2} gives
	\envHLinePdPt[3.2_Lem2Pr_Eq0help3]{\approx}
	{
		\mpEVHi{1,1,1,l_4}{T}
	}
	{
		\opF{1}{6} \fcMZH{1}{l_4}T^3 + \opF{1}{2} \fcMZH{2}{1,l_4} T^2
		\qquad
		(l_4>1)
	}
Moreover,
	because of \refEq{3.2_Lem3iEq3},
	applying $\mpEVH[]{}$ to both sides of \refEq{3.2_Lem2Pr_EqHarRel4} yields
	\envHLineFCmPt{\approx}
	{
		\mpEVHi{1,1,1,1}{T}
	}
	{
		\opF{1}{4} \mpEVHi{1,1,1}{T}T - \opF{1}{4} \mpEVHi{1,1,2}{T} 
	}
	{
		\opF{1}{4} \bkR[a]{ \opF{1}{6} T^3+ \fcMZH{2}{1,1} T + \fcMZH{3}{1,1,1} } T - \opF{1}{8} \fcMZH{1}{2} T^2
	}
	{
		\opF{1}{24}T^4 + \bkR[a]{ \opF{1}{4} \fcMZH{2}{1,1} - \opF{1}{8} \fcMZH{1}{2} } T^2
	}
	which together with \refEq{3.2_Lem3iiEq3} gives
	\envHLinePdPt[3.2_Lem2Pr_Eq0help4]{\approx}
	{
		\mpEVHi{1,1,1,1}{T}
	}
	{
		\opF{1}{24}T^4 + \opF{1}{2} \fcMZH{2}{1,1} T^2 
	}
Then,
	we can deduce from 
	\refEq{3.2_Lem3zPr_Eq} for $n=4$, \refEq{3.2_Lem2Pr_Eq0help1}, \refEq{3.2_Lem2Pr_Eq0help2}, \refEq{3.2_Lem2Pr_Eq0help3}, and \refEq{3.2_Lem2Pr_Eq0help4} that
	\envHLineCmPt{\approx}
	{
		\mpEVHi{\lbfL_4}{T}
	}
	{
		\opF{\fcMChSo{4}{l_1,l_2,l_3,l_4}}{24} T^4
		+
		\opF{\fcMChSo{3}{l_1,l_2,l_3}}{6} \fcMZH{1}{l_4}T^3
		+
		\opF{\fcMChSo{2}{l_1,l_2}}{2}\fcMZH{2}{l_3,l_4}  T^2 
		\quad
		(\lbfL_4\in\setN^4)
	}
	which proves \refEq{3.2_Lem2_Eq4}.
}

\section{Proofs}\label{sectFour}
\subsection{Proof of \refThm{1_Thm1}}\label{sectFourOne}
We first prove \refEq{1_Thm1_EqDpt2} and \refEq{1_Thm1_EqDpt3} in \refThm{1_Thm1}.
\envProof[\refEq{1_Thm1_EqDpt2}]{
By \refEq{3.1_Prop1_Eq},
	we easily obtain 
	\envHLineCm[4.1_1_Thm1PrDp2_EqHar]
	{
		\racF{\fcMZHO{2}}{ \mpS{\gpC{2}} }
	}
	{
		{\fcMZHO{1}}^{\otimes2} - \fcMZO{1} \circ \mpWT{2}
	}
	which proves \refEq{1_Thm1_EqDpt2} for $\letB=*$.

We can deduce the following identities from \refEq{3.2_Prop1_EqFuncDp1} and \refEq{3.2_Prop1_EqFuncDp2}: 
	\envHLineCFLaaCm[4.1_1_Thm1PrDp2_EqShaHelp]
	{
		{\fcMZHO{1}}^{\otimes2}
	}
	{
		{\fcMZSO{1}}^{\otimes2}
	}
	{
		\racF{\fcMZHO{2}}{ \mpS{\gpC{2}}}
	}
	{
		\racF{\fcMZSO{2}}{ \mpS{\gpC{2}}}
		-
		\fcMChSoO{2} \cdot \fcMZO{1} \circ \mpWT{2} 
	}	
	where we have used in the second identity the property that $\fcMChSoO{2} \cdot \fcMZO{1} \circ \mpWT{2} $ is invariant under $\gpSym{2}$,
	or
	\envMPd{
		\racF{ \fcMChSoO{2} \cdot \fcMZO{1} \circ \mpWT{2} }{ \mpS{\gpC{2}} }
	}{
		2\fcMChSoO{2} \cdot \fcMZO{1} \circ \mpWT{2}
	}
Since $\fcMChSO{2}+\fcMChSoO{2}=\fcCVo[2]$,
	substituting \refEq{4.1_1_Thm1PrDp2_EqShaHelp} into \refEq{4.1_1_Thm1PrDp2_EqHar} yields
	\envHLineCm[4.1_1_Thm1PrDp2_EqSha]
	{
		\racF{\fcMZSO{2}}{ \mpS{\gpC{2}} }
	}
	{
		{\fcMZSO{1}}^{\otimes2} - \fcMChSO{2} \cdot \fcMZSO{1} \circ \mpWT{2}
	}
	which proves \refEq{1_Thm1_EqDpt2} for $\letB=\sh$.
}
\envProof[\refEq{1_Thm1_EqDpt3}]{
Since
	$\gpC{3}=\Set{\gpu, (123), (132)}$ and $\pstUT=\Set{\gpu,(23),(123)}$,
	direct calculations give the following equations in $\setZ[\gpSym{3}]$:
	\envHLineCFLqqPd
	{
		(123)\mpS{\gpC{3}} 
	}
	{
		\mpS{\gpC{3}} 
	}
	{
		\mpS{\pstUT} \mpS{\gpC{3}}
	}
	{
		\mpS{\gpSym{3}} + \mpS{\gpC{3}}
	}
We thus see from \refEq{3.1_Prop2_Eq1} that
	\envHLinePd[4.1_1_Thm1PrDp3_EqHarHelp]
	{
		\racF{\fcMZHO{2}\otimes\fcMZHO{1}}{ \mpS{\gpC{3}} }
	}
	{
		\racF[n]{\fcMZHO{3}}{ \mpS{\gpSym{3}} + \mpS{\gpC{3}} } 
		+ 
		\racF{ \fcMZHO{2} \circ (\empWT{3}{2,1} + \empWT{3}{1,2}) }{ \mpS{\gpC{3}} } 
	}
Subtracting \refEq{4.1_1_Thm1PrDp3_EqHarHelp} from \refEq{3.1_Prop2_Eq2},
	we obtain
	\envMPd
	{
		{\fcMZHO{1}}^{\otimes3} - \racF{\fcMZHO{2}\otimes\fcMZHO{1}}{ \mpS{\gpC{3}} }
	}
	{
		- \racF{\fcMZHO{3}}{ \mpS{\gpC{3}} } + \fcMZO{1} \circ \mpWT{3}
	}
This identity is equivalent to
	\envHLineCm[4.1_1_Thm1PrDp3_EqHar]
	{
		\racF{\fcMZHO{3}}{ \mpS{\gpC{3}} } 
	}
	{
		- {\fcMZHO{1}}^{\otimes3} + \racF{\fcMZHO{2}\otimes\fcMZHO{1}}{ \mpS{\gpC{3}} } + \fcMZO{1} \circ \mpWT{3}
	}
	which proves \refEq{1_Thm1_EqDpt3} for $\letB=*$.

We can deduce the following identities from \refEq{3.2_Prop1_EqFuncDp1}, \refEq{3.2_Prop1_EqFuncDp2}, and \refEq{3.2_Prop1_EqFuncDp3}:
	\envHLineCSCmNme
	{\label{4.1_1_Thm1PrDp3_EqShaHelp1}
		{\fcMZHO{1}}^{\otimes3}
	}
	{
		{\fcMZSO{1}}^{\otimes3}
	}
	{\label{4.1_1_Thm1PrDp3_EqShaHelp2}
		\racF{\fcMZHO{2}\otimes\fcMZHO{1}}{ \mpS{\gpC{3}} }
	}
	{
		\racF{\fcMZSO{2}\otimes\fcMZSO{1}}{ \mpS{\gpC{3}} }
		-
		\opF{1}{2} \racF{ (\fcMChSoO{2} \cdot \fcMZO{1} \circ \mpWT{2}) \otimes \fcMZSO{1} }{ \mpS{\gpC{3}} }
	}
	{\label{4.1_1_Thm1PrDp3_EqShaHelp3}
		\racF{\fcMZHO{3}}{ \mpS{\gpC{3}} }
	}
	{
		\racF{ \fcMZSO{3} }{ \mpS{\gpC{3}} }
		- 
		\opF{1}{2} \racF{ (\fcMChSoO{2} \cdot \fcMZO{1} \circ \mpWT{2}) \otimes \fcMZSO{1} }{ \mpS{\gpC{3}} }
		+
		\fcMChSoO{3} \cdot \fcMZO{1} \circ \mpWT{3} 
	}
	where we have used in the third identity the property that
	$\fcMChSoO{3} \cdot (\fcMZO{1} \circ \mpWT{3})$ is invariant under $\gpSym{3}$.
Since $\fcMChSO{3}+\fcMChSoO{3}=\fcCVo[3]$,
	substituting \refEq{4.1_1_Thm1PrDp3_EqShaHelp1}, \refEq{4.1_1_Thm1PrDp3_EqShaHelp2}, and \refEq{4.1_1_Thm1PrDp3_EqShaHelp3} 
	into \refEq{4.1_1_Thm1PrDp3_EqHar} yields
	\envHLineCm[4.1_1_Thm1PrDp3_EqSha]
	{
		\racF{\fcMZSO{3}}{ \mpS{\gpC{3}} } 
	}
	{
		- {\fcMZSO{1}}^{\otimes3} + \racF{\fcMZSO{2}\otimes\fcMZSO{1}}{ \mpS{\gpC{3}} } + \fcMChSO{3} \cdot \fcMZO{1} \circ \mpWT{3}
	}
	which proves \refEq{1_Thm1_EqDpt3} for $\letB=\sh$,
	and completes the proof.
}

We now prepare two lemmas before proving \refEq{1_Thm1_EqDpt4}, 
	because the proof of \refEq{1_Thm1_EqDpt4} is more complicated than those of \refEq{1_Thm1_EqDpt2} and \refEq{1_Thm1_EqDpt3}.
The identities of \refLem{4.1_Lem1} \respTx{\refLem{4.1_Lem2}} correspond to \refEq{4.1_1_Thm1PrDp3_EqHarHelp} 
	\respTx{\refEq{4.1_1_Thm1PrDp3_EqShaHelp1}, \refEq{4.1_1_Thm1PrDp3_EqShaHelp2}, and \refEq{4.1_1_Thm1PrDp3_EqShaHelp3} } 
	in the proof of \refEq{1_Thm1_EqDpt3}.
	
\begin{lemma}\label{4.1_Lem1}
We have  
	\envHLineCSNmePd
	{\label{4.1_Lem1_EqHar1}\hspace{-20pt}
		\racF{ \fcMZHO{3}\otimes\fcMZHO{1} }{ \mpS{\gpC{4}} }
	}
	{
		\racF[n]{ \fcMZHO{4} }{ 2\mpS{\gpC{4}} + (12)\mpS{\gpC{4}} + (34)\mpS{\gpC{4}} }
		\lnAH
		+
		\racF[n]{ \fcMZHO{3} \circ (\empWT{4}{2,1,1}+\empWT{4}{1,2,1}+\empWT{4}{1,1,2}) }{ \mpS{\gpAlt{4}} - (13)\mpS{\gpC{4}} - (23)\mpS{\gpC{4}} }
	}
	{\label{4.1_Lem1_EqHar2}\hspace{-20pt}\rule{0cm}{23pt}
		\racF{ {\fcMZHO{2}}^{\otimes2} }{ \mpS{\gpCs{4}} }
	}
	{
		\racF[n]{ \fcMZHO{4} }{ \mpS{\gpC{4}} + (14)\mpS{\gpC{4}} + (23)\mpS{\gpC{4}} }
		\lnAH
		+
		\racF{ \fcMZHO{3} \circ (\empWT{4}{2,1,1}+\empWT{4}{1,2,1}+\empWT{4}{1,1,2})  }{ (23)\mpS{\gpC{4}} }
		+
		\racF{ \fcMZHO{2} \circ \empWT{4}{2,2}  }{ (23)\mpS{\gpCs{4}} }
	}
	{\label{4.1_Lem1_EqHar3}\hspace{-20pt}\rule{0cm}{23pt}
		\racF{ \fcMZHO{2}\otimes{\fcMZHO{1}}^{\otimes2} }{ \mpS{\gpC{4}} }
	}
	{
		\racF[n]{ \fcMZHO{4} }{ 2\mpS{\gpSym{4}} + \mpS{\gpC{4}} - (13)\mpS{\gpC{4}} }
		\lnAH
		+
		\racF[n]{ \fcMZHO{3} \circ (\empWT{4}{2,1,1}+\empWT{4}{1,2,1}+\empWT{4}{1,1,2})  }{ 2\mpS{\gpAlt{4}} - (13)\mpS{\gpC{4}} }
		\lnAH[]
		+
		\fcMZHO{2} \circ \bkR[b]{
			\racF[n]{  \empWT{4}{2,2}  }{ \mpS{\gpC{4}} + 2(23)\mpS{\gpCs{4}} }
			+
			\racF{ (\empWT{4}{3,1}+\empWT{4}{1,3})  }{ \mpS{\gpC{4}} }
		}
		\nonumber
	}
\end{lemma}
\begin{lemma}\label{4.1_Lem2}
We have
	\envPLineNme{\HLineCSCmNme[p]
	{\label{4.1_Lem1_EqShaHelp1}\hspace{-15pt}
		{\fcMZHO{1}}^{\otimes4}
	}
	{
		{\fcMZSO{1}}^{\otimes4}
	}
	{\label{4.1_Lem1_EqShaHelp2}\hspace{-15pt}\rule{0cm}{23pt}
		\racF{ \fcMZHO{2}\otimes{\fcMZHO{1}}^{\otimes2} }{ \mpS{\gpC{4}} }
	}
	{
		\racF{ \fcMZSO{2}\otimes{\fcMZSO{1}}^{\otimes2} }{ \mpS{\gpC{4}} }
		-
		\opF{1}{2} \racF{ (\fcMChSoO{2}\cdot \fcMZO{1} \circ \mpWT{2})\otimes{\fcMZSO{1}}^{\otimes2} }{ \mpS{\gpC{4}} }
	}
	{\label{4.1_Lem1_EqShaHelp3}\hspace{-15pt}\rule{0cm}{23pt}
		\racF{ {\fcMZHO{2}}^{\otimes2} }{ \mpS{\gpCs{4}} }
	}
	{
		\racF{ {\fcMZSO{2}}^{\otimes2} }{ \mpS{\gpCs{4}} }
		- 
		\opF{1}{2} \racF{ (\fcMChSoO{2} \cdot \fcMZO{1} \circ \mpWT{2}) \otimes \fcMZSO{2} }{ \mpS{\gpC{4}} }
		+
		\opF{5}{4} \fcMChSoO{4} \cdot \fcMZO{1} \circ \mpWT{4}
	}\HLineCFNmePd
	{\label{4.1_Lem1_EqShaHelp4}\hspace{-15pt}\rule{0cm}{23pt}
		\racF{ \fcMZHO{3}\otimes\fcMZHO{1} }{ \mpS{\gpC{4}} }
	}
	{
		\racF{ \fcMZSO{3}\otimes\fcMZSO{1} }{ \mpS{\gpC{4}} }		- 
		\opF{1}{2} \racF{ (\fcMChSoO{2} \cdot \fcMZO{1} \circ \mpWT{2}) \otimes {\fcMZSO{1}}^{\otimes2} }{ \mpS{\gpC{4}} }
		\lnAH
		+
		\opF{1}{3} \racF{ (\fcMChSoO{3} \cdot \fcMZO{1} \circ \mpWT{3}) \otimes\fcMZSO{1} }{ \mpS{\gpC{4}} }
	}
	{\label{4.1_Lem1_EqShaHelp5}\hspace{-15pt}\rule{0cm}{23pt}
		\racF{ \fcMZHO{4} }{ \mpS{\gpC{4}} } 
	}
	{
		\racF{ \fcMZSO{4} }{ \mpS{\gpC{4}} }
		- 
		\opF{1}{2} \racF{ (\fcMChSoO{2} \cdot \fcMZO{1} \circ \mpWT{2}) \otimes \fcMZSO{2} }{ \mpS{\gpC{4}} }
		\lnAH
		+
		\opF{1}{3} \racF{ (\fcMChSoO{3} \cdot \fcMZO{1} \circ \mpWT{3}) \otimes \fcMZSO{1} }{ \mpS{\gpC{4}} } 
		+
		\opF{1}{4} \fcMChSoO{4} \cdot \fcMZO{1} \circ \mpWT{4} 
	}
	}
\end{lemma}

We now prove \refEq{1_Thm1_EqDpt4}.
We will then discuss proofs of \refLem[s]{4.1_Lem1} and \ref{4.1_Lem2}.

\envProof[identity \refEq{1_Thm1_EqDpt4}]{
Direct calculations show that
	\envPLine[4.1_1_Thm1Pr_EqLeftCstGpC1]{\HLineCECm[p]
	{
		\mpS{\gpC{4}}
	}
	{
		\gpu + (13)(24) + (1234) + (1432)
	}
	{
		(12)\mpS{\gpC{4}}
	}
	{
		(12) + (143) + (234) + (1324)
	}
	{
		(13)\mpS{\gpC{4}}
	}
	{
		(13) + (24) + (12)(34) + (14)(23)
	}
	{
		(14)\mpS{\gpC{4}}
	}
	{
		(14) + (123) + (243) + (1342)
	}\HLineCFCm
	{
		(23)\mpS{\gpC{4}}
	}
	{
		(23) + (134) + (142) + (1243)
	}
	{
		(34)\mpS{\gpC{4}}
	}
	{
		(34) + (124) + (132) + (1423)
		\nonumber
	}
	}
	from which we see that
	\envHLineCm[4.1_1_Thm1PrDp4_EqLeftCstGpC2]
	{
		\mpS{\gpSym{4}}
	}
	{
		\bkR{ \gpu + (12) + (13) + (14) + (23) + (34) } \mpS{\gpC{4}} 
	}
	i.e.,
	\envMO{
		\Set{ \gpC{4}, (12)\gpC{4}, (13)\gpC{4}, (14)\gpC{4}, (23)\gpC{4}, (34)\gpC{4} }
	}
	gives a left $\gpC{4}$-coset decomposition of $\gpSym{4}$.
By \refEq{4.1_1_Thm1PrDp4_EqLeftCstGpC2},
	the sum of \refEq{4.1_Lem1_EqHar1} and \refEq{4.1_Lem1_EqHar2} yields
	\envHLinePd
	{
		\racF{ \fcMZHO{3}\otimes\fcMZHO{1} }{ \mpS{\gpC{4}}} + \racF{ {\fcMZHO{2}}^{\otimes2} }{ \mpS{\gpCs{4}}}
	}
	{
		\racF[n]{ \fcMZHO{4} }{ \mpS{\gpSym{4}} + 2\mpS{\gpC{4}} - (13)\mpS{\gpC{4}} }
		\lnAH
		+
		\racF[n]{ \fcMZHO{3} \circ (\empWT{4}{2,1,1}+\empWT{4}{1,2,1}+\empWT{4}{1,1,2})  }{ \mpS{\gpAlt{4}} - (13)\mpS{\gpC{4}} }
		\lnAH
		+
		\racF{ \fcMZHO{2} \circ \empWT{4}{2,2}  }{ (23)\mpS{\gpCs{4}} }
	}
Subtracting \refEq{4.1_Lem1_EqHar3} from this identity,
	we obtain
	\envLinePd[4.1_1_Thm1PrDp4_EqHarHelp]
	{
		\racF{ \fcMZHO{3}\otimes\fcMZHO{1} }{ \mpS{\gpC{4}} } + \racF{ {\fcMZHO{2}}^{\otimes2} }{ \mpS{\gpCs{4}}}
		-
		\racF{ \fcMZHO{2}\otimes{\fcMZHO{1}}^{\otimes2} }{ \mpS{\gpC{4}} }
	}
	{
		- 
		\racF[n]{ \fcMZHO{4} }{ \mpS{\gpSym{4}} - \mpS{\gpC{4}} }
		\lnAH
		-
		\racF{ \fcMZHO{3} \circ (\empWT{4}{2,1,1}+\empWT{4}{1,2,1}+\empWT{4}{1,1,2})  }{ \mpS{\gpAlt{4}}  }
		\lnAH[]
		-
		\fcMZHO{2} \circ \bkR[b]{
			\racF[n]{  \empWT{4}{2,2}  }{ \mpS{\gpC{4}} + (23)\mpS{\gpCs{4}} }
			+
			\racF{ (\empWT{4}{3,1}+\empWT{4}{1,3})  }{ \mpS{\gpC{4}} }
		}
		\nonumber
	}
We see from \refEq{3.1_PL_DefGpDpt4X} and the equivalence classes modulo $\Gp{(12), (34)}$ in \refTab{3.1_PL_Table1}
	that 
	\envMFCmPt{\equiv}
	{
		\mpS{\pstXF}
	}
	{
		(14) + (23) + \mpS{\gpC{4}}
	}
	{
		(23) + (134) + \mpS{\gpC{4}}
	}
	{
		(23)\mpS{\gpCs{4}} + \mpS{\gpC{4}}
		\mod \Gp{(12), (34)}
	}
	and so 
	\envHLinePd
	{
		\racF{ \empWT{4}{2,2} }{ \mpS{\pstXF} }
	}
	{
		\racF[n]{  \empWT{4}{2,2}  }{ \mpS{\gpC{4}} + (23)\mpS{\gpCs{4}} }
	}
Thus the sum of \refEq{3.1_Prop3_Eq4} and \refEq{4.1_1_Thm1PrDp4_EqHarHelp} yields
	\envMPd
	{
		\racF{ \fcMZHO{3}\otimes\fcMZHO{1} }{ \mpS{\gpC{4}} } + \racF{ {\fcMZHO{2}}^{\otimes2} }{ \mpS{\gpCs{4}}}
		-
		\racF{ \fcMZHO{2}\otimes{\fcMZHO{1}}^{\otimes2} }{ \mpS{\gpC{4}} } + {\fcMZHO{1}}^{\otimes4}
	}
	{
		\racF{ \fcMZHO{4} }{ \mpS{\gpC{4}} }  + \fcMZO{1} \circ \mpWT{4}
	}
This identity is equivalent to
	\envHLineCm[4.1_1_Thm1PrDp4_EqHar]
	{
		\racF{ \fcMZHO{4} }{ \mpS{\gpC{4}} } 
	}
	{
		{\fcMZHO{1}}^{\otimes4}
		-
		\racF{ \fcMZHO{2}\otimes{\fcMZHO{1}}^{\otimes2} }{ \mpS{\gpC{4}} }
		+
		\racF{ {\fcMZHO{2}}^{\otimes2} }{ \mpS{\gpCs{4}} }
		+
		\racF{ \fcMZHO{3}\otimes\fcMZHO{1} }{ \mpS{\gpC{4}} }
		- 
		\fcMZO{1} \circ \mpWT{4}
	}
	which proves \refEq{1_Thm1_EqDpt4} for $\letB=*$. 

Combining \refEq{4.1_Lem1_EqShaHelp1}, \refEq{4.1_Lem1_EqShaHelp2}, \refEq{4.1_Lem1_EqShaHelp3}, and \refEq{4.1_Lem1_EqShaHelp4}
	(or considering  
	\refEq{4.1_Lem1_EqShaHelp1} $-$ \refEq{4.1_Lem1_EqShaHelp2} $+$ \refEq{4.1_Lem1_EqShaHelp3} $+$ \refEq{4.1_Lem1_EqShaHelp4},
	roughly speaking),
	we can restate the right-hand side of \refEq{4.1_1_Thm1PrDp4_EqHar} as
	\envHLinePd[4.1_1_Thm1PrDp4_EqShaHelp1]
	{
		\text{(RHS of \refEq{4.1_1_Thm1PrDp4_EqHar})}
	}
	{
		{\fcMZSO{1}}^{\otimes4}
		-
		\racF{ \fcMZSO{2}\otimes{\fcMZSO{1}}^{\otimes2} }{ \mpS{\gpC{4}} }
		+
		\racF{ {\fcMZSO{2}}^{\otimes2} }{ \mpS{\gpCs{4}} }
		+
		\racF{ \fcMZSO{3}\otimes\fcMZSO{1} }{ \mpS{\gpC{4}} }
		\lnAH
		- 
		\opF{1}{2} \racF{ (\fcMChSoO{2} \cdot \fcMZO{1} \circ \mpWT{2}) \otimes \fcMZSO{2} }{ \mpS{\gpC{4}} }
		+
		\opF{1}{3} \racF{ (\fcMChSoO{3} \cdot \fcMZO{1} \circ \mpWT{3}) \otimes \fcMZSO{1} }{ \mpS{\gpC{4}} } 
		\lnAH[]
		+
		(\opF{5}{4} \fcMChSoO{4} - \fcCVo[4] ) \cdot (\fcMZO{1} \circ \mpWT{4}) 
		\nonumber
	}
Equating \refEq{4.1_Lem1_EqShaHelp5}, \refEq{4.1_1_Thm1PrDp4_EqHar}, and \refEq{4.1_1_Thm1PrDp4_EqShaHelp1},
	we obtain
	\envHLineCm
	{
		\racF{ \fcMZSO{4} }{ \mpS{\gpC{4}} } 
	}
	{
		{\fcMZSO{1}}^{\otimes4}
		-
		\racF{ \fcMZSO{2}\otimes{\fcMZSO{1}}^{\otimes2} }{ \mpS{\gpC{4}} }
		+
		\racF{ {\fcMZSO{2}}^{\otimes2} }{ \mpS{\gpCs{4}} }
		+
		\racF{ \fcMZSO{3}\otimes\fcMZSO{1} }{ \mpS{\gpC{4}} }
		+
		(\fcMChSoO{4} - \fcCVo[4] ) \cdot (\fcMZO{1} \circ \mpWT{4}) 
	}
	which together with $\fcMChSO{4}+\fcMChSoO{4}=\fcCVo[4]$ proves \refEq{1_Thm1_EqDpt4} for $\letB=\sh$,
	and completes the proof.
}

We will show \refLem[s]{4.1_Lem1} and \ref{4.1_Lem2} for the completeness of the proof of \refEq{1_Thm1_EqDpt4}.
We first prove \refLem{4.1_Lem2}.

\envProof[\refLem{4.1_Lem2}]{
Identity \refEq{4.1_Lem1_EqShaHelp1} immediately follows from \refEq{3.2_Prop1_EqFuncDp1}.
We see from  \refEq{3.2_Prop1_EqFuncDp1} and \refEq{3.2_Prop1_EqFuncDp2} that
	\envHLinePd
	{
		\fcMZHO{2}\otimes{\fcMZHO{1}}^{\otimes2} 
	}
	{
		\fcMZSO{2}\otimes{\fcMZSO{1}}^{\otimes2} - \opF{1}{2} (\fcMChSoO{2}\cdot \fcMZO{1} \circ \mpWT{2})\otimes{\fcMZSO{1}}^{\otimes2}
	}
Applying $\mpS{\gpC{4}}$ to both sides of this equation,
	we obtain \refEq{4.1_Lem1_EqShaHelp2}.
In a similar way,
	we obtain \refEq{4.1_Lem1_EqShaHelp4} from \refEq{3.2_Prop1_EqFuncDp1} and \refEq{3.2_Prop1_EqFuncDp3}.
We also obtain  \refEq{4.1_Lem1_EqShaHelp5} by applying $\mpS{\gpC{4}}$ to both sides of \refEq{3.2_Prop1_EqFuncDp4},
	since $\fcMChSoO{4} \cdot (\fcMZO{1} \circ \mpWT{4}) $ is invariant under $\gpSym{4}$.
	
We prove \refEq{4.1_Lem1_EqShaHelp3}.
We easily see that 
	$\racF{ f\otimes g }{ (13)(24) }= g\otimes f$ for any functions $f$ and $g$ of two variables,
	and so we obtain from \refEq{3.2_Prop1_EqFuncDp2} that
	\envHLineThPd[4.1_Lem1Proof_ShaHelp3_Eq1]
	{\hspace{-15pt}
		{\fcMZHO{2}}^{\otimes2} 
	}
	{
		\bkR[B]{ \fcMZSO{2} - \opF{1}{2} \fcMChSoO{2} \cdot \fcMZO{1} \circ \mpWT{2} } 
		\otimes 
		\bkR[B]{ \fcMZSO{2} - \opF{1}{2} \fcMChSoO{2} \cdot \fcMZO{1} \circ \mpWT{2} }
	}
	{	
		{\fcMZSO{2}}^{\otimes2}  
		- 
		\opF{1}{2} \racF[n]{ (\fcMChSoO{2} \cdot \fcMZO{1} \circ \mpWT{2}) \otimes \fcMZSO{2} }{ \gpu + (13)(24) }
		+
		\opF{1}{4} (\fcMChSoO{2} \cdot \fcMZO{1} \circ \mpWT{2})^{\otimes2}
	}
We see from \refEq{3.2_Prop1Pr_Eq0}, \refEq{3.2_Prop1Pr_EqDp4_3}, and $\fcMChSoO{2}{}^{\otimes2}=\fcMChSoO{4}$ that
	\envOFLineCm
	{
		(\fcMChSoO{2} \cdot \fcMZO{1} \circ \mpWT{2})^{\otimes2}
	}
	{
		\fcMZ{1}{2}^2\fcMChSoO{2}{}^{\otimes2} 
	}
	{
		\opF{5}{2} \fcMZ{1}{4} \fcMChSoO{4}
	}
	{
		\opF{5}{2} \fcMChSoO{4} \cdot \fcMZO{1} \circ \mpWT{4}
	}
	by which we can restate \refEq{4.1_Lem1Proof_ShaHelp3_Eq1} as
	\envHLinePd[4.1_Lem1Proof_ShaHelp3_Eq2]
	{
		{\fcMZHO{2}}^{\otimes2} 
	}
	{	
		{\fcMZSO{2}}^{\otimes2}
		- 
		\opF{1}{2} \racF{ (\fcMChSoO{2} \cdot \fcMZO{1} \circ \mpWT{2}) \otimes \fcMZSO{2} }{ \mpS{\Gp{(13)(24)}} }
		+
		\opF{5}{8} \fcMChSoO{4} \cdot \fcMZO{1} \circ \mpWT{4}
	}
Since $\gpCs{4}=\Set{\gpu,(1234)}\subset\gpC{4}=\Set{\gpu,(1234),(13)(24),(1432)}$, 
	\envHLinePd[4.1_1_Thm1PrDp4_EqHar22help1]
	{
		\mpS{\Gp{(13)(24)}} \mpS{\gpCs{4}} 
	}
	{
		\mpS{\gpC{4}}
	}
Applying $\mpS{\gpCs{4}}$ to both sides of \refEq{4.1_Lem1Proof_ShaHelp3_Eq2},
	we thus obtain \refEq{4.1_Lem1_EqShaHelp3},
	which 
	completes the proof.
}

We now prove \refLem{4.1_Lem1}.

\envProof[\refLem{4.1_Lem1}]{
Let $\sig\in\Set{(12),(23),(34)}$.
By the equivalence classes modulo $\Gp{\sig}$ in \refTab{3.1_PL_Table1} and straightforward calculations,
	\refEq{3.1_Lem5iiPr_Eq3help2} yields 
	\envHLineCmPt[4.1_Lem1Proof_EqModAltSymHelp]{\equiv}
	{
		2\mpS{\gpAlt{4}}
	}
	{
		\mpS{\gpSym{4}}
		\qquad
		\modd\Gp{\sig}
	}
	and \refEq{4.1_1_Thm1Pr_EqLeftCstGpC1} yields
	\envPLine[4.1_Lem1Proof_EqModLeftCosetHelp]
	{\begin{array}{rclrclrclc}
		\rule{0pt}{15pt}	\mpS{\gpC{4}}&\equiv&(12)\mpS{\gpC{4}},\quad 	
			&(13)\mpS{\gpC{4}}&\equiv&(34)\mpS{\gpC{4}},\quad 		
			&(14)\mpS{\gpC{4}}&\equiv&(23)\mpS{\gpC{4}}	
			&\qquad\modd\Gp{(12)},\\
		\rule{0pt}{15pt}	\mpS{\gpC{4}}&\equiv&(23)\mpS{\gpC{4}},\quad 	
			&(12)\mpS{\gpC{4}}&\equiv&(34)\mpS{\gpC{4}},\quad 		
			&(13)\mpS{\gpC{4}}&\equiv&(14)\mpS{\gpC{4}}
			&\qquad\modd\Gp{(23)},\\
		\rule{0pt}{15pt}	\mpS{\gpC{4}}&\equiv&(34)\mpS{\gpC{4}},\quad 	
			&(12)\mpS{\gpC{4}}&\equiv&(13)\mpS{\gpC{4}},\quad 		
			&(14)\mpS{\gpC{4}}&\equiv&(23)\mpS{\gpC{4}}
			&\qquad\modd\Gp{(34)}.
	\end{array}}
Thus,
	we deduce from \refEq{4.1_1_Thm1PrDp4_EqLeftCstGpC2} that 
	\envHLineCmPt[4.1_Lem1Proof_EqModLeftCoset]{\equiv}
	{
		\mpS{\gpAlt{4}}
	}
	{
		\alp\mpS{\gpC{4}} + \beta\mpS{\gpC{4}} + \gam\mpS{\gpC{4}}
		\qquad
		\modd\Gp{\sig}
	}
	where $(\alp,\beta,\gam)$ is a $3$-tuple of $\Set{ \gpu, (12), (13), (14), (23), (34) }$ such that
	\envPLine[4.1_Lem1Proof_EqModLeftCosetCond]
	{
		\envCaseThPd{
			\alp\in\Set{\gpu,(12)},\;\; \beta\in\Set{(13),(34)},\;\; \gam\in\Set{(14),(23)}
			\quad&
			(\sig=(12))
		}{
			\alp\in\Set{\gpu,(23)},\;\; \beta\in\Set{(12),(34)},\;\; \gam\in\Set{(13),(14)}
			&
			(\sig=(23))
		}{
			\alp\in\Set{\gpu,(34)},\;\; \beta\in\Set{(12),(13)},\;\; \gam\in\Set{(14),(23)}
			&
			(\sig=(34))
		}
	}

We now prove \refEq{4.1_Lem1_EqHar1}.
Since either $g\gpC{4}=h\gpC{4}$ or $g\gpC{4}\cap h\gpC{4}=\phi$ for any $g, h\in\gpSym{4}$,
	we can see from the first and second equations of \refEq{4.1_1_Thm1Pr_EqLeftCstGpC1} that
	\envHLineCFLaaCm
	{
		(1234)\mpS{\gpC{4}}
	}
	{
		\mpS{\gpC{4}}
	}
	{
		(234)\mpS{\gpC{4}}
	}
	{
		(12)\mpS{\gpC{4}}
	}
	respectively.
By \refEq{3.1_PL_DefGpDpt4U} and \refEq{3.1_Lem5iPr_Eq00}, 
	we obtain 
	\envOTLinePd
	{
		\mpS{\pstUF} \mpS{\gpC{4}}
	}
	{
		\mpS{\gpC{4}}+(34)\mpS{\gpC{4}}+(234)\mpS{\gpC{4}}+(1234)\mpS{\gpC{4}}
	}
	{
		2\mpS{\gpC{4}} + (12)\mpS{\gpC{4}} + (34)\mpS{\gpC{4}}
	}
Thus,
	applying  $\mpS{\gpC{4}}$ to both sides of \refEq{3.1_Prop3_Eq1} yields
	\envHLinePd[4.1_1_Thm1PrDp4_EqHar31help1]
	{
		\racF{ \fcMZHO{3}\otimes\fcMZHO{1} }{ \mpS{\gpC{4}} }
	}
	{
		\racF[n]{\fcMZHO{4}}{ 2\mpS{\gpC{4}} + (12)\mpS{\gpC{4}} + (34)\mpS{\gpC{4}} } 
		\lnAH
		+
		\fcMZHO{3} \circ \bkR[b]{ \racF{  (\empWT{4}{2,1,1}+ \empWT{4}{1,2,1}) }{ (12)\mpS{\gpC{4}} } + \racF{ \empWT{4}{1,1,2} }{ \mpS{\gpC{4}} } }
	}
We know from \refEq{4.1_Lem1Proof_EqModLeftCoset} and \refEq{4.1_Lem1Proof_EqModLeftCosetCond} that
	\envHLinePt{\equiv}
	{
		\mpS{\gpAlt{4}} 
	}
	{
		(13)\mpS{\gpC{4}}  
		+
		(23)\mpS{\gpC{4}}
		+
		\envCaseTCm{
			(12)\mpS{\gpC{4}}
			\quad&
			\modd \Gp{(12)}  \;\;\mathrm{or}\;\;  \Gp{(23)}
		}{
			\mpS{\gpC{4}}
			\quad&
			\modd \Gp{(34)}
		}
	}
	and so 
	\envHLinePt{\equiv}
	{
		\racF[n]{ \empWT{4}{i,j,k} }{ \mpS{\gpAlt{4}} - (13)\mpS{\gpC{4}} - (23)\mpS{\gpC{4}} }
	}
	{
		\envCaseTCm{
			\racF{ \empWT{4}{i,j,k} }{ (12)\mpS{\gpC{4}} }
			\quad&
			((i,j,k)\in J)
		}{
			\racF{ \empWT{4}{i,j,k} }{ \mpS{\gpC{4}} }
			\quad&
			((i,j,k)=(1,1,2))
		}
	}
	where $J$ means the set $\Set{(2,1,1),(1,2,1)}$.
Therefore we have
	\envLinePd
	{
		\fcMZHO{3} \circ \bkR[b]{ \racF{  (\empWT{4}{2,1,1}+ \empWT{4}{1,2,1}) }{ (12)\mpS{\gpC{4}} } + \racF{ \empWT{4}{1,1,2} }{ \mpS{\gpC{4}} } }
	}
	{
		\racF[n]{ \fcMZHO{3} \circ (\empWT{4}{2,1,1}+\empWT{4}{1,2,1}+\empWT{4}{1,1,2}) }{ \mpS{\gpAlt{4}} - (13)\mpS{\gpC{4}} - (23)\mpS{\gpC{4}} }
	}
Substituting this into \refEq{4.1_1_Thm1PrDp4_EqHar31help1}  proves \refEq{4.1_Lem1_EqHar1}.

We can easily see that
	\envHLineCFLaaCm
	{
		\mpS{\pstVF[0]}
	}
	{
		(23) \mpS{\Gp{(13)(24)}}
	}
	{
		\mpS{\pstVF}
	}
	{
		( \gpu+(123)+(23) )\mpS{\Gp{(13)(24)}}
	}
	which together with \refEq{4.1_1_Thm1PrDp4_EqHar22help1} give
	\envHLineCFLaaCm
	{
		\mpS{\pstVF[0]} \mpS{\gpCs{4}}
	}
	{
		(23) \mpS{\gpC{4}}
	}
	{
		\mpS{\pstVF} \mpS{\gpCs{4}}
	}
	{
		\mpS{\gpC{4}} + (14)\mpS{\gpC{4}} + (23)\mpS{\gpC{4}}  
	}
	respectively,
	where we note that $(123)\mpS{\gpC{4}}=(14)\mpS{\gpC{4}}$ by the fourth equation of \refEq{4.1_1_Thm1Pr_EqLeftCstGpC1}.
Thus,
	applying  $\mpS{\gpCs{4}}$ to both sides of \refEq{3.1_Prop3_Eq2}, 
	we obtain \refEq{4.1_Lem1_EqHar2}.

Lastly, we prove \refEq{4.1_Lem1_EqHar3}.	
We can obtain the following identity by applying $\mpS{\gpC{4}}$ to both sides of \refEq{3.1_Prop3_Eq3}:
	\envHLinePd[4.1_Lem1Pr_EqqHar3help]
	{\hspace{-20pt}
		\racF{ \fcMZHO{2}\otimes{\fcMZHO{1}}^{\otimes2} }{ \mpS{\gpC{4}} }
	}
	{
		\racF[n]{\fcMZHO{4}}{ 2 \mpS{\gpSym{4}} + \mpS{\gpC{4}} - (13)\mpS{\gpC{4}} } 
		\lnAH
		+
		\racF[n]{ \fcMZHO{3} \circ (\empWT{4}{2,1,1} + \empWT{4}{1,2,1} + \empWT{4}{1,1,2}) }{ 2\mpS{\gpAlt{4}} + (23)\mpS{\gpC{4}} - (14)\mpS{\gpC{4}} }
		\lnAH
		-
		\fcMZHO{3} \circ \bkR[b]{ \racF{ \empWT{4}{2,1,1} }{ (34)\mpS{\gpC{4}} } + \racF{ \empWT{4}{1,2,1} }{ \mpS{\gpC{4}} } + \racF{ \empWT{4}{1,1,2} }{ (12)\mpS{\gpC{4}} } }
		\lnAH[]
		+
		\fcMZHO{2} \circ \bkR[b]{
			\racF[n]{  \empWT{4}{2,2}  }{ \mpS{\gpC{4}} + 2(23)\mpS{\gpCs{4}} }
			+
			\racF{ (\empWT{4}{3,1}+\empWT{4}{1,3})  }{ \mpS{\gpC{4}} }
		}
		\nonumber
	}
(We will prove \refEq{4.1_Lem1Pr_EqqHar3help} in \refLem{4.1_Lem3} below because the proof is not short.)
We can also obtain by \refEq{4.1_Lem1Proof_EqModLeftCosetHelp} 
	\envHLinePt{\equiv}
	{
		(23)\mpS{\gpC{4}} + (13)\mpS{\gpC{4}} 
	}
	{
		(14)\mpS{\gpC{4}}
		+
		\envCaseThPd{
			(34)\mpS{\gpC{4}}		 	&	\modd	\Gp{(12)}
		}{
			\mpS{\gpC{4}}				&	\modd	\Gp{(23)}
		}{
			(12)\mpS{\gpC{4}}		 	&	\modd	\Gp{(34)}
		}
	}
Thus,
	$\racF{  (\empWT{4}{2,1,1}+\empWT{4}{1,2,1}+\empWT{4}{1,1,2})  }{ (13)\mpS{\gpC{4}} }$ can be expressed as
	\envLineCm
	{
		\racF{  (\empWT{4}{2,1,1}+\empWT{4}{1,2,1}+\empWT{4}{1,1,2})  }{ (13)\mpS{\gpC{4}} }
	}
	{
		\racF[n]{  (\empWT{4}{2,1,1}+\empWT{4}{1,2,1}+\empWT{4}{1,1,2})  }{ - (23)\mpS{\gpC{4}} + (14)\mpS{\gpC{4}} }
		\lnAH
		+
		\bkR{ \racF{ \empWT{4}{2,1,1} }{ (34)\mpS{\gpC{4}} } + \racF{ \empWT{4}{1,2,1} }{ \mpS{\gpC{4}} } + \racF{ \empWT{4}{1,1,2} }{ (12)\mpS{\gpC{4}} } }
	}
	which gives
	\envLinePd
	{
		\racF[n]{ \fcMZHO{3} \circ (\empWT{4}{2,1,1}+\empWT{4}{1,2,1}+\empWT{4}{1,1,2})  }{ 2\mpS{\gpAlt{4}} - (13)\mpS{\gpC{4}} }
	}
	{
		\racF[n]{ \fcMZHO{3} \circ (\empWT{4}{2,1,1} + \empWT{4}{1,2,1} + \empWT{4}{1,1,2}) }{ 2\mpS{\gpAlt{4}} + (23)\mpS{\gpC{4}} - (14)\mpS{\gpC{4}} }
		\lnAH
		-
		\fcMZHO{3} \circ \bkR[b]{ \racF{ \empWT{4}{2,1,1} }{ (34)\mpS{\gpC{4}} } + \racF{ \empWT{4}{1,2,1} }{ \mpS{\gpC{4}} } + \racF{ \empWT{4}{1,1,2} }{ (12)\mpS{\gpC{4}} } }
	}
Combining \refEq{4.1_Lem1Pr_EqqHar3help} and this equation  proves \refEq{4.1_Lem1_EqHar3}.
}

\begin{lemma}\label{4.1_Lem3}\mbox{}\\
{\bf(i)}
Let $\sig\in\Set{(12),(23),(34)}$.
Then the following congruence equations hold:
	\envHLineCSNmePdPt{\equiv}
	{\label{4.1_Lem3i_Eq1}
		\mpS{\pstWF[0]} \mpS{\gpC{4}}  
	}
	{
		\mpS{\gpC{4}} + 2(23)\mpS{\gpCs{4}}
		\qquad\hspace{51pt}	
		\modd	\Gp{(12), (34)}
	}
	{\label{4.1_Lem3i_Eq2}
		\mpS{\pstWF[1]} \mpS{\gpC{4}}  
	}
	{
		2\mpS{\gpAlt{4}} + (23)\mpS{\gpC{4}} - (14)\mpS{\gpC{4}}
		\qquad	
		\modd	\Gp{\sig}
	}
	{\label{4.1_Lem3i_Eq3}
		\mpS{\pstWF} \mpS{\gpC{4}} 
	}
	{
		2 \mpS{\gpSym{4}} + \mpS{\gpC{4}} - (13)\mpS{\gpC{4}}
		\qquad\hspace{21pt}
		\modd	\Gp{\gpu}
	}
{\bf(ii)}
Identity \refEq{4.1_Lem1Pr_EqqHar3help} holds.
	
\end{lemma}
\envProof{
We first prove the assertion (i). 
We see from \refEq{3.1_PL_DefGpDpt4W} and \refEq{3.1_Lem5iPr_Eq00} that
	\envMCm{
		\mpS{\pstWF[0]} \mpS{\gpC{4}}  
	}{
		(23)\mpS{\gpC{4}} + (24)\mpS{\gpC{4}}
	}
	and from the third equation of \refEq{4.1_1_Thm1Pr_EqLeftCstGpC1}, we see that
	\envMPd
	{
		(24)\mpS{\gpC{4}}
	}
	{
		(13)\mpS{\gpC{4}} 
	}
We thus obtain 
	\envHLinePd[4.1_Lem3iPr_Eqq1]
	{
		\mpS{\pstWF[0]} \mpS{\gpC{4}}  
	}
	{
		(13)\mpS{\gpC{4}} + (23)\mpS{\gpC{4}}
	}
Equation \refEq{4.1_Lem3iPr_Eqq1} proves \refEq{4.1_Lem3i_Eq1},
	since
	\envHLineCFCmPt{\equiv}
	{
		(13)\mpS{\gpC{4}} 
	}
	{
		(1432) + (1234) + \gpu + (13)(24)
	\lnP{\equiv}
		\mpS{\gpC{4}}
		\qquad	
		\modd	\Gp{(12), (34)}
	}
	{
		(23)\mpS{\gpC{4}}
	}
	{
		2( (23) + (134) )
	\lnP{\equiv}
		2(23)\mpS{\gpCs{4}}
		\qquad\hspace{48pt}	
		\modd	\Gp{(12), (34)}
	}
	which can be seen from the equivalence classes modulo $\Gp{(12), (34)}$ in \refTab{3.1_PL_Table1}.
By virtue of \refEq{4.1_1_Thm1PrDp4_EqLeftCstGpC2},
	calculations similar to \refEq{4.1_Lem3iPr_Eqq1} show that
	\envHLineFiCm[4.1_Lem3iPr_Eqq2]
	{
		\mpS{\pstWF[1]} \mpS{\gpC{4}}  
	}
	{
		\mpS{\Set{(34),(1234),(1243),(1324)}} \mpS{\gpC{4}} + \mpS{\pstWF[0]} \mpS{\gpC{4}} 
	}
	{
		\bkR{ (34)\mpS{\gpC{4}} + \mpS{\gpC{4}} + (23)\mpS{\gpC{4}} + (12)\mpS{\gpC{4}}  }
		+
		\bkR{ (13)\mpS{\gpC{4}} + (23)\mpS{\gpC{4}} }
	}
	{
		\mpS{\gpC{4}} + (12)\mpS{\gpC{4}} + (13)\mpS{\gpC{4}} + 2(23)\mpS{\gpC{4}} + (34)\mpS{\gpC{4}} 
	}
	{
		\mpS{\gpSym{4}} + (23)\mpS{\gpC{4}} - (14)\mpS{\gpC{4}}
	}
	and 
	\envHLineFiPd[4.1_Lem3iPr_Eqq3]
	{
		\mpS{\pstWF} \mpS{\gpC{4}} 
	}
	{
		\mpS{\Set{\gpu[4],(13)(24),(123),(124),(234),(243)}} \mpS{\gpC{4}} + \mpS{\pstWF[1]} \mpS{\gpC{4}} 
	}
	{
		\bkR{ \mpS{\gpC{4}}+\mpS{\gpC{4}}+(14)\mpS{\gpC{4}}+(34)\mpS{\gpC{4}}+(12)\mpS{\gpC{4}}+(14)\mpS{\gpC{4}} } 
		\lnAH
		+ 
		\bkR{ \mpS{\gpSym{4}} + (23)\mpS{\gpC{4}} - (14)\mpS{\gpC{4}} }
	}
	{
		\mpS{\gpSym{4}} 
		+
		2\mpS{\gpC{4}} + (12)\mpS{\gpC{4}} + (14)\mpS{\gpC{4}} + (23)\mpS{\gpC{4}} + (34)\mpS{\gpC{4}} 
	}
	{
		2 \mpS{\gpSym{4}} + \mpS{\gpC{4}} - (13)\mpS{\gpC{4}}
	}
Then we obtain \refEq{4.1_Lem3i_Eq2} by \refEq{4.1_Lem1Proof_EqModAltSymHelp} and \refEq{4.1_Lem3iPr_Eqq2},
	and 
	\refEq{4.1_Lem3i_Eq3} by \refEq{4.1_Lem3iPr_Eqq3}.

We now prove the assertion (ii), or \refEq{4.1_Lem1Pr_EqqHar3help}.
We can deduce from \refEq{4.1_Lem3i_Eq1}, \refEq{4.1_Lem3i_Eq2}, and \refEq{4.1_Lem3i_Eq3} that
	\envHLineCSCmNme
	{\label{4.1_Lem3iiPr_Eqq1}
		\fcMZHO{2} \circ \racF{ \empWT{4}{2,2} }{ \mpS{\pstWF[0]}\mpS{\gpC{4}} } 
	}
	{
		\fcMZHO{2} \circ \racF[n]{ \empWT{4}{2,2} }{ \mpS{\gpC{4}} + 2(23)\mpS{\gpCs{4}}} 
	}
	{\label{4.1_Lem3iiPr_Eqq2}
		\racF{ \fcMZHO{3} \circ \empWT{4}{i,j,k} }{ \mpS{\pstWF[1]} \mpS{\gpC{4}}   }
	}
	{
		\racF[n]{ \fcMZHO{3} \circ \empWT{4}{i,j,k} }{ 2\mpS{\gpAlt{4}} + (23)\mpS{\gpC{4}} - (14)\mpS{\gpC{4}} }
	}
	{\label{4.1_Lem3iiPr_Eqq3}
		\racF{\fcMZHO{4}}{ \mpS{\pstWF} \mpS{\gpC{4}}  } 
	}
	{
		\racF[n]{\fcMZHO{4}}{ 2 \mpS{\gpSym{4}} + \mpS{\gpC{4}} - (13)\mpS{\gpC{4}} } 
	}
	respectively,
	where $(i,j,k)\in\Set{(2,1,1),(1,2,1),(1,1,2)}$.
Applying $\mpS{\gpC{4}}$ to both sides of \refEq{3.1_Prop3_Eq3} 
	and 
	substituting \refEq{4.1_Lem3iiPr_Eqq1}, \refEq{4.1_Lem3iiPr_Eqq2}, and \refEq{4.1_Lem3iiPr_Eqq3} into it,
	we obtain
	\envLinePd[4.1_Lem3iiPr_EqMain]
	{\hspace{-20pt}
		\racF{ \fcMZHO{2}\otimes{\fcMZHO{1}}^{\otimes2} }{ \mpS{\gpC{4}} }
	}
	{
		\racF[n]{\fcMZHO{4}}{ 2 \mpS{\gpSym{4}} + \mpS{\gpC{4}} - (13)\mpS{\gpC{4}} } 
		\lnAH
		+
		\racF[n]{ \fcMZHO{3} \circ (\empWT{4}{2,1,1} + \empWT{4}{1,2,1} + \empWT{4}{1,1,2}) }{ 2\mpS{\gpAlt{4}} + (23)\mpS{\gpC{4}} - (14)\mpS{\gpC{4}} }
		\lnAH
		-
		\fcMZHO{3} \circ \bkR[b]{ \racF{ \empWT{4}{2,1,1} }{ (34)\mpS{\gpC{4}} } + \racF{ \empWT{4}{1,2,1} }{ (1234)\mpS{\gpC{4}} } + \racF{ \empWT{4}{1,1,2} }{ (1324)\mpS{\gpC{4}} } }
		\lnAH[]
		+
		\fcMZHO{2} \circ \bkR[b]{
			\racF[n]{ \empWT{4}{2,2} }{ \mpS{\gpC{4}} + 2(23)\mpS{\gpCs{4}} } 
			+
			\racF{ \empWT{4}{3,1} }{ (24)\mpS{\gpC{4}} }
			+
			\racF{ \empWT{4}{1,3} }{ \mpS{\gpC{4}} }
		}
		\nonumber
	} 
We see from the third equation of \refEq{4.1_1_Thm1Pr_EqLeftCstGpC1} and the equivalence classes modulo $\Gp{(12),(123)}$ in \refTab{3.1_PL_Table1} that
	\envMFCmPt{\equiv}{
		(24)\mpS{\gpC{4}} 
	}{
		(13) + (24) + (12)(34) + (14)(23)
	}{
		\gpu + (13)(24) + (1234) + (1432)
	}{
		\mpS{\gpC{4}}
		\mod 
		\Gp{(12),(123)}
	}
	and so
	\envHLinePd
	{
		\fcMZHO{2} \circ \racF{ \empWT{4}{3,1} }{ (24)\mpS{\gpC{4}} }
	}
	{
		\fcMZHO{2} \circ \racF{ \empWT{4}{3,1} }{ \mpS{\gpC{4}} }
	}
This equation together with \refEq{4.1_Lem3iiPr_EqMain} proves \refEq{4.1_Lem1Pr_EqqHar3help},
	since
	\envM{
		(1234)\mpS{\gpC{4}}
	}{
		\mpS{\gpC{4}}
	}
	and
	\envM{
		(1324)\mpS{\gpC{4}}
	}{
		(12)\mpS{\gpC{4}}
	}
	by the first and second equations of \refEq{4.1_1_Thm1Pr_EqLeftCstGpC1},
	respectively.
}

\subsection{Proof of \refCor{1_Cor1}}\label{sectFourTwo}
Let $\stP{n}$ be the set of partitions of $\Set{1,\ldots,n}$,
	and let $\stPh{n}{m}$ be its subset consisting of partitions $\Pi=\bkB{P_1,\ldots,P_m}$
	such that the number of the parts is $m$.
Note that $\stP{n}=\UnT[t]{m=1}{n}\stPh{n}{m}$.
We identify a partition $\Pi=\Set{ \Set{n_1^{(1)},\ldots,n_{a_1}^{(1)}}, \ldots, \Set{n_1^{(m)},\ldots,n_{a_m}^{(m)}} }$ with 
	\envPLinePd
	{
		n_1^{(1)}\ldots n_{a_1}^{(1)} |\ldots| n_1^{(m)}\ldots n_{a_m}^{(m)}
	}	
For example,
	$\Set{\Set{1,2,3}}=123$, $\Set{\Set{1,2},\Set{3}}=12|3$, and $\Set{\Set{1},\Set{2}, \Set{3}}=1|2|3$.
Let $n$ and $n'$ be positive integers with $n<n'$.
For convenience, we embedded $\gpSym{n}$ into $\gpSym{n'}$ in the following way:
	a permutation $\SMatT{1&\ldots&n}{j_1&\ldots&j_n}$ of $\gpSym{n}$ is identified 
	with the permutation $\SMatT{1&\ldots&n&n+1&\ldots&n'}{j_1&\ldots&j_n&n+1&\ldots&n'}$ of $\gpSym{n'}$,
	which fixes integers between $n+1$ and $n'$. 

We require the following lemmas to prove \refEq{1_Cor_EqSymRMZV} of \refCor{1_Cor1} for $n=2,3$, and $4$.
We assume that $\lbfL_n=(l_1,\ldots,l_n)\in\setN^n$ and $\letB\in\Set{*,\sh}$ in the lemmas.

\begin{lemma}[Case of depth $2$]\label{4.2_Lem1} 
	\envHLineCFNmePd
	{\label{4.2_Lem1_Eq1}
		{\fcMZBO{1}}^{\otimes2}(\lbfL_2)
	}
	{
		\Sm{\Pi\in\stPh{2}{2}} \fcZPB{\lbfL_2}{\Pi}
	}
	{\label{4.2_Lem1_Eq2}
		\Fcr{ \fcMChBO{2} \cdot \fcMZ[]{1}{} \circ \mpWT{2} }{\lbfL_2}
	}
	{
		\Sm{\Pi\in\stPh{2}{1}} \fcZPB{\lbfL_2}{\Pi} 
	}
\end{lemma}

\begin{lemma}[Case of depth $3$]\label{4.2_Lem2} 
	\envHLineCSNmePd
	{\label{4.2_Lem2_Eq1}
		\racFA{ {\fcMZBO{1}}^{\otimes3} }{ \mpS{\gpSym{2}} }{ \lbfL_3 }
	}
	{
		2 \Sm{\Pi\in\stPh{3}{3}} \fcZPB{\lbfL_3}{\Pi}
	}
	{\label{4.2_Lem2_Eq2}
		\racFA{ \fcMZBO{2}\otimes\fcMZBO{1} }{ \mpS{\gpC{3}}\mpS{\gpSym{2}} }{ \lbfL_3 }
	}
	{
		3\Sm{\Pi\in\stPh{3}{3}} \fcZPB{\lbfL_3}{\Pi}
		-
		\Sm{\Pi\in\stPh{3}{2}} \fcZPB{\lbfL_3}{\Pi}
	}
	{\label{4.2_Lem2_Eq3}
		\racFA{ \fcMChBO{3}\cdot\fcMZO{1}\circ\mpWT{3} }{ \mpS{\gpSym{2}} }{ \lbfL_3 }
	}
	{
		2 \Sm{\Pi\in\stPh{3}{1}} \fcZPB{\lbfL_3}{\Pi} 
	}
\end{lemma}
\begin{lemma}[Case of depth $4$]\label{4.2_Lem3} 
	\envPLineNme{\HLineCSCmNme[p]
	{\label{4.2_Lem3_Eq1}\hspace{-20pt}
		\racFA{ {\fcMZBO{1}}^{\otimes4} }{ \mpS{\gpSym{3}} }{ \lbfL_4 }
	}
	{
		6 \Sm{\Pi\in\stPh{4}{4}} \fcZPB{\lbfL_4}{\Pi}
	}
	{\label{4.2_Lem3_Eq2}\hspace{-20pt}
		\racFA{ \fcMZBO{2}\otimes{\fcMZBO{1}}^{\otimes2} }{ \mpS{\gpC{4}} \mpS{\gpSym{3}} }{ \lbfL_4 }
	}
	{
		12 \Sm{\Pi\in\stPh{4}{4}} \fcZPB{\lbfL_4}{\Pi} - 2 \Sm{\Pi\in\stPh{4}{3}} \fcZPB{\lbfL_4}{\Pi}
	}
	{\label{4.2_Lem3_Eq3}\hspace{-20pt}
		\racFA{ {\fcMZBO{2}}^{\otimes2} }{ \mpS{\gpCs{4}} \mpS{\gpSym{3}} }{ \lbfL_4 }
	}
	{
		3 \Sm{\Pi\in\stPh{4}{4}} \fcZPB{\lbfL_4}{\Pi}
		-
		\Sm{\Pi\in\stPh{4}{3}} \fcZPB{\lbfL_4}{\Pi}
		+
		\Sm{\Pi\in\stPhTh{4}{2}{2,2}} \fcZPB{\lbfL_4}{\Pi}
	}\HLineCFCmNme
	{\label{4.2_Lem3_Eq4}\hspace{-20pt}
		\racFA{ \fcMZBO{3}\otimes\fcMZBO{1} }{ \mpS{\gpC{4}} \mpS{\gpSym{3}} }{ \lbfL_4 }
	}
	{
		4 \Sm{\Pi\in\stPh{4}{4}} \fcZPB{\lbfL_4}{\Pi}
		-
		2 \Sm{\Pi\in\stPh{4}{3}} \fcZPB{\lbfL_4}{\Pi}
		+
		2 \Sm{\Pi\in\stPhTh{4}{2}{3,1}} \fcZPB{\lbfL_4}{\Pi}
	}
	{\label{4.2_Lem3_Eq5}\hspace{-20pt}
		\racFA{ \fcMChBO{4}\cdot\fcMZO{1}\circ\mpWT{4} }{ \mpS{\gpSym{3}} }{ \lbfL_4 }
	}
	{
		6 \Sm{\Pi\in\stPh{4}{1}} \fcZPB{\lbfL_4}{\Pi}
	}
	}
	where
	$\stPhTh{4}{2}{2,2}$ and $\stPhTh{4}{2}{3,1}$ in \refEq{4.2_Lem3_Eq3} and \refEq{4.2_Lem3_Eq4} are subsets in $\stPh{4}{2}$ defined by
	\envHLineCFLaaCmDef
	{
		\stPhTh{4}{2}{2,2}
	}
	{
		\Set{12|34, 13|24, 14|23}	
	}
	{
		\stPhTh{4}{2}{3,1}
	}
	{
		\Set{123|4, 124|3, 134|2, 234|1}	
	}
	respectively.
Note that
	\envMPd{
		\stPh{4}{2}
	}{
		\stPhTh{4}{2}{2,2} \cup \stPhTh{4}{2}{3,1}
	}
\end{lemma}

We will give proofs of \refEq{1_Cor_EqSymRMZV} for $n=2,3$, and $4$ before discussing those of the lemmas.
 
\envProof[\refEq{1_Cor_EqSymRMZV} for $n=2$]{
Substituting \refEq{4.2_Lem1_Eq1} and \refEq{4.2_Lem1_Eq2} into the right-hand side of \refEq{1_Thm1_EqDpt2} yields
	\envHLineCm
	{
		\racFA{\fcMZBO{2}}{ \mpS{\gpSym{2}} }{\lbfL_2}
	}
	{
		\Sm{\Pi\in\stPh{2}{2}} \fcZPB{\lbfL_2}{\Pi} - \Sm{\Pi\in\stPh{2}{1}} \fcZPB{\lbfL_2}{\Pi} 
	}
	since
	\envM
	{
		\gpC{2}
	}
	{
		\gpSym{2}
	}
	and
	\envMPd
	{
		\fcMChB{2}{\lbfL_2}\fcMZ{1}{L_2}
	}
	{
		\Fcr{ \fcMChBO{2} \cdot \fcMZ[]{1}{} \circ \mpWT{2} }{\lbfL_2}
	}
We have by definition (see \refEq{1_PL_DefNmPartition})
	\envHLine
	{
		\pfcCP{2}{\Pi}
	}
	{
		\envCaseTCm{
			1
			&
			(\Pi\in\stPh{2}{2})
		}{
			-1
			&
			(\Pi\in\stPh{2}{1})
		}
	}
	and thus we obtain by $\stP{2}=\UnT[t]{m=1}{2}\stPh{2}{m}$
	\envHLineCm
	{
		\racFA{\fcMZBO{2}}{ \mpS{\gpSym{2}} }{\lbfL_2}
	}
	{
		\Sm{\Pi\in\stP{2}} \pfcCP{2}{\Pi} \fcZPB{\lbfL_2}{\Pi}
	}
	which proves \refEq{1_Cor_EqSymRMZV} for $n=2$.
}

\envProof[\refEq{1_Cor_EqSymRMZV} for $n=3$]{
Applying $\mpS{\gpSym{2}}$ to both sides of \refEq{1_Thm1_EqDpt3}, 
	we obtain
	\envHLineCm
	{
		\racFA{\fcMZBO{3}}{ \mpS{\gpSym{3}} }{\lbfL_3}
	}
	{
		-
		\racFA{ {\fcMZBO{1}}^{\otimes3} }{ \mpS{\gpSym{2}} }{ \lbfL_3 }
		+
		\racFA{ \fcMZBO{2}\otimes\fcMZBO{1} }{ \mpS{\gpC{3}} \mpS{\gpSym{2}} }{ \lbfL_3 }
		+
		\racFA{ \fcMChBO{3}\cdot\fcMZO{1}\circ\mpWT{3} }{ \mpS{\gpSym{2}} }{ \lbfL_3 }
	}
	where we have used $\mpS{\gpSym{3}}=\mpS{\gpC{3}}\mpS{\gpSym{2}}$ on the left-hand side of this equation.
Substituting \refEq{4.2_Lem2_Eq1}, \refEq{4.2_Lem2_Eq2}, and \refEq{4.2_Lem2_Eq3} into the right-hand side of this equation yields
	\envHLineCm
	{
		\racFA{\fcMZBO{3}}{ \mpS{\gpSym{3}} }{\lbfL_3}
	}
	{
		\Sm{\Pi\in\stPh{3}{3}} \fcZPB{\lbfL_3}{\Pi}
		-
		\Sm{\Pi\in\stPh{3}{2}} \fcZPB{\lbfL_3}{\Pi}
		+
		2 \Sm{\Pi\in\stPh{3}{1}} \fcZPB{\lbfL_3}{\Pi} 
	}
	which proves \refEq{1_Cor_EqSymRMZV} for $n=3$,
	since 
	\envHLine
	{
		\pfcCP{3}{\Pi}
	}
	{
		\envCaseThCm{
			1
			&
			(\Pi\in\stPh{3}{3})
		}{
			-1
			&
			(\Pi\in\stPh{3}{2})
		}{
			2
			&
			(\Pi\in\stPh{3}{1})
		}
	}
	and 
	\envMPd{
		\stP{3}
	}
	{
		\UnT[t]{m=1}{3}\stPh{3}{m}
	}
}

\envProof[\refEq{1_Cor_EqSymRMZV} for $n=4$]{
We see can from \refEq{4.1_1_Thm1Pr_EqLeftCstGpC1} and \refEq{4.1_1_Thm1PrDp4_EqLeftCstGpC2} that
	\envMPd
	{
		\mpS{\gpSym{4}}
	}
	{
		\mpS{\gpSym{3}}\mpS{\gpC{4}}
	}
Taking the inverses of both sides of this equation,
	we obtain
	\envMPd
	{
		\mpS{\gpSym{4}}
	}
	{
		\mpS{\gpC{4}}\mpS{\gpSym{3}}
	}
Thus,
	applying $\mpS{\gpSym{3}}$ to both sides of \refEq{1_Thm1_EqDpt4}
	and 
	combining the identities in \refLem{4.2_Lem3}
	(or considering  
	\refEq{4.2_Lem3_Eq1} $-$ \refEq{4.2_Lem3_Eq2} $+$ \refEq{4.2_Lem3_Eq3} $+$ \refEq{4.2_Lem3_Eq4} $-$ \refEq{4.2_Lem3_Eq5}),
	we can obtain
	\envHLineCm
	{
		\racFA{\fcMZBO{4}}{ \mpS{\gpSym{4}} }{\lbfL_4}
	}
	{
		\Sm{\Pi\in\stPh{4}{4}} \fcZPB{\lbfL_4}{\Pi}
		-
		\Sm{\Pi\in\stPh{4}{3}} \fcZPB{\lbfL_4}{\Pi}
		+
		\Sm{\Pi\in\stPhTh{4}{2}{2,2}} \fcZPB{\lbfL_4}{\Pi}
		+
		2 \Sm{\Pi\in\stPhTh{4}{2}{3,1}} \fcZPB{\lbfL_4}{\Pi}
		\lnAH
		-
		6 \Sm{\Pi\in\stPh{4}{1}} \fcZPB{\lbfL_4}{\Pi}
	}
	which proves \refEq{1_Cor_EqSymRMZV} for $n=4$,
	since
	\envHLine
	{
		\pfcCP{4}{\Pi}
	}
	{
		\envCaseFiCm{
			1
			&
			(\Pi\in\stPh{4}{4})
		}{
			-1
			&
			(\Pi\in\stPh{4}{3})
		}{
			1
			&
			(\Pi\in\stPhTh{4}{2}{2,2})
		}{
			2
			&
			(\Pi\in\stPhTh{4}{2}{3,1})
		}{
			-6
			&
			(\Pi\in\stPh{4}{1})
		}
	}
	\envMCm{
		\stPh{4}{2}
	}
	{
		\stPhTh{4}{2}{2,2}\cup\stPhTh{4}{2}{3,1}
	}
	and 
	\envMPd{
		\stP{4}
	}
	{
		\UnT[t]{m=1}{4}\stPh{4}{m}
	}
}

Recall that $\fcMZ{1}{1}=\fcMZB{1}{1}=0$ and $\fcMZ{1}{k}=\fcMZB{1}{k}$.
By the definition of $\fcMChB[]{1}{}$,
	\envHLine
	{
		\fcMChB{1}{k}\fcMZ{1}{k}
	}
	{
		\fcMZB{1}{k}
	}
	for any positive integer $k$;
	this will be used repeatedly below.

We now give proofs of \refLem[s]{4.2_Lem1}, \ref{4.2_Lem2}, and \ref{4.2_Lem3}.

\envProof[\refLem{4.2_Lem1}]{
We have 
	\envM
	{
		\stPh{2}{2}
	}
	{
		\Set{1|2}
	}
	and
	\envM
	{
		\stPh{2}{1}
	}
	{
		\Set{12}
	}
	by definition.
Thus, 
	\envHLineCFCm
	{
		\Sm{\Pi\in\stPh{2}{2}} \fcZPB{\lbfL_2}{\Pi}
	}
	{
		\fcMChB{1}{l_1}\fcMZ{1}{l_1} \fcMChB{1}{l_2}\fcMZ{1}{l_2} 
	\lnP{=}
		{\fcMZBO{1}}^{\otimes2}(\lbfL_2)
	}
	{
		\Sm{\Pi\in\stPh{2}{1}} \fcZPB{\lbfL_2}{\Pi}
	}
	{
		\fcMChB{2}{l_1,l_2}\fcMZ{1}{l_1+l_2}
	\lnP{=}
		\Fcr{ \fcMChBO{2} \cdot \fcMZ[]{1}{} \circ \mpWT{2} }{\lbfL_2}
	}
	which prove \refEq{4.2_Lem1_Eq1} and \refEq{4.2_Lem1_Eq2},
	respectively.
}

\envProof[\refLem{4.2_Lem2}]{
We have 
	\envMCm
	{
		\stPh{3}{3}
	}
	{
		\Set{1|2|3}
	}
	\envMCm
	{
		\stPh{3}{2}
	}
	{
		\Set{12|3,13|2,23|1}
	}
	and
	\envM
	{
		\stPh{3}{1}
	}
	{
		\Set{123}
	}
	by definition.
In particular,
	$\stPh{3}{2}$ is expressed as $\Un[t]{\sig\in\gpC{3}} \Set{\sig^{-1}(1)\sig^{-1}(2)|\sig^{-1}(3)}$,
	and so
	\envHLinePd
	{
		\Sm{\Pi\in\stPh{3}{2}} \fcZPB{\lbfL_3}{\Pi}
	}
	{
		\Sm{\sig\in\gpC{3}} \fcMChB{2}{l_{\sig^{-1}(1)},l_{\sig^{-1}(2)}}\fcMZ{1}{l_{\sig^{-1}(1)}+l_{\sig^{-1}(2)}} \fcMChB{1}{l_{\sig^{-1}(3)}}\fcMZ{1}{l_{\sig^{-1}(3)}}
	}
Thus, 
	\envHLineCSNmePd
	{\label{4.2_Lem2Proof_EqRight1}
		\Sm{\Pi\in\stPh{3}{3}} \fcZPB{\lbfL_3}{\Pi}
	}
	{
		{\fcMZBO{1}}^{\otimes3}(\lbfL_3)
	}
	{\label{4.2_Lem2Proof_EqRight2}
		\Sm{\Pi\in\stPh{3}{2}} \fcZPB{\lbfL_3}{\Pi}
	}
	{
		\racFA{ (\fcMChBO{2}\cdot\fcMZBO{1}\circ\mpWT{2})\otimes\fcMZBO{1} }{ \mpS{\gpC{3}} }{ \lbfL_3 }
	}
	{\label{4.2_Lem2Proof_EqRight3}
		\Sm{\Pi\in\stPh{3}{1}} \fcZPB{\lbfL_3}{\Pi}
	}
	{
		\Fcr{ \fcMChBO{3}\cdot\fcMZO{1}\circ\mpWT{3} }{\lbfL_3} 
	}	

Since ${\fcMZBO{1}}^{\otimes3}$ is invariant under $\gpSym{3}$,
	we have
	\envMCm
	{
		\racFA{ {\fcMZBO{1}}^{\otimes3} }{ \mpS{\gpSym{2}} }{\lbfL_3}
	}
	{
		2 \Fc{ {\fcMZBO{1}}^{\otimes3} }{\lbfL_3}
	}
	which together with \refEq{4.2_Lem2Proof_EqRight1} proves \refEq{4.2_Lem2_Eq1}.
Similarly,
	we have
	\envMCm
	{
		\racFA{ \fcMChBO{3}\cdot\fcMZO{1}\circ\mpWT{3} }{ \mpS{\gpSym{2}} }{\lbfL_3}
	}
	{
		2 \Fc{ \fcMChBO{3} \cdot \fcMZO{1} \circ \mpWT{3} }{\lbfL_3}
	}
	which together with \refEq{4.2_Lem2Proof_EqRight3} proves \refEq{4.2_Lem2_Eq3}.	
We know from \refEq{1_Thm1_EqDpt2} that 
	\envHLinePd[4.2_Lem2Proof_EqPL1]
	{
		\racF{\fcMZBO{2}}{ \mpS{\gpSym{2}} } 
	}
	{
		{\fcMZBO{1}}^{\otimes2} - \fcMChBO{2} \cdot \fcMZO{1} \circ \mpWT{2}
	}
Since $\gpC{3}\gpSym{2}=\gpSym{2}\gpC{3}$,
	\envHLineFCm
	{
		\racF{ \fcMZBO{2}\otimes\fcMZBO{1} }{ \mpS{\gpC{3}}\mpS{\gpSym{2}} }
	}
	{
		\bkR{ \racF{ \fcMZBO{2}\otimes\fcMZBO{1} }{ \mpS{\gpSym{2}} } } | \mpS{\gpC{3}} 
	}
	{
		\racFr{ {\fcMZBO{1}}^{\otimes3} - (\fcMChBO{2} \cdot \fcMZO{1} \circ \mpWT{2}) \otimes\fcMZBO{1} }{ \mpS{\gpC{3}}  }
	}
	{
		3 {\fcMZBO{1}}^{\otimes3} 
		-
		\racF{ (\fcMChBO{2} \cdot \fcMZO{1} \circ \mpWT{2}) \otimes\fcMZBO{1} }{ \mpS{\gpC{3}} }
	}
	which together with \refEq{4.2_Lem2Proof_EqRight1} and \refEq{4.2_Lem2Proof_EqRight2} proves \refEq{4.2_Lem2_Eq2},
	and completes the proof.
}
	
\envProof[\refLem{4.2_Lem3}]{
Let $\gpAlts{4}$ be the subset of $\gpAlt{4}$ given by
	\envOTLineThPd
	{
		\gpAlts{4}
	}
	{
		\Set{\gpu, (13)(24), (123), (132), (142), (234)}
	}  
	{
		\Gp{(13)(24)}\gpC{3}
	}
Note that 
	\envMPd
	{
		\mpS{\gpAlt{4}}
	}
	{
		\mpS{\Gp{(12)(34)}}\mpS{\gpAlts{4}}
	}
We can see from the definitions of $\stPh{4}{m}$ and $\stPhTh{4}{2}{i,j}$ that
	\envPLine{\HLineCECm[p]
	{
		\stPh{4}{4}
	}
	{
		\Set{1|2|3|4}
	}
	{
		\stPh{4}{3}
	}
	{
		\Un{\sig\in\gpAlts{4}} \Set{\sig^{-1}(1)\sig^{-1}(2)|\sig^{-1}(3)|\sig^{-1}(4)}
	}
	{
		\stPhTh{4}{2}{2,2}
	}
	{
		\Un{\sig\in\gpC{3}} \Set{\sig^{-1}(1)\sig^{-1}(2)|\sig^{-1}(3)\sig^{-1}(4)}
	}
	{
		\stPhTh{4}{2}{3,1}
	}
	{
		\Un{\sig\in\gpC{4}} \Set{\sig^{-1}(1)\sig^{-1}(2)\sig^{-1}(3)|\sig^{-1}(4)}
	}\HLineCTPd
	{
		\stPh{4}{1}
	}
	{
		\Set{1234}
	}
	}
Thus, 
	\envPLineNme{\HLineCECmNme[p]
	{\label{4.2_Lem3Proof_EqRight1}
		\Sm{\Pi\in\stPh{4}{4}} \fcZPB{\lbfL_4}{\Pi}
	}
	{
		{\fcMZBO{1}}^{\otimes4}(\lbfL_4)
	}
	{\label{4.2_Lem3Proof_EqRight2}
		\Sm{\Pi\in\stPh{4}{3}} \fcZPB{\lbfL_4}{\Pi}
	}
	{
		\racFA{ (\fcMChBO{2}\cdot\fcMZO{1}\circ\mpWT{2})\otimes{\fcMZBO{1}}^{\otimes2} }{ \mpS{\gpAlts{4}} }{ \lbfL_4 }
	}
	{\label{4.2_Lem3Proof_EqRight3.1}
		\Sm{\Pi\in\stPhTh{4}{2}{2,2}} \fcZPB{\lbfL_4}{\Pi}
	}
	{
		\racFA{ (\fcMChBO{2}\cdot\fcMZO{1}\circ\mpWT{2})^{\otimes2} }{ \mpS{\gpC{3}} }{ \lbfL_4 }
	}
	{\label{4.2_Lem3Proof_EqRight3.2}
		\Sm{\Pi\in\stPhTh{4}{2}{3,1}} \fcZPB{\lbfL_4}{\Pi}
	}
	{
		\racFA{ (\fcMChBO{3}\cdot\fcMZO{1}\circ\mpWT{3})\otimes\fcMZBO{1} }{ \mpS{\gpC{4}} }{ \lbfL_4 }
	}\HLineCTNmePd
	{\label{4.2_Lem3Proof_EqRight4}
		\Sm{\Pi\in\stPh{4}{1}} \fcZPB{\lbfL_4}{\Pi}
	}
	{
		\Fcr{ \fcMChBO{4} \cdot \fcMZO{1} \circ \mpWT{4} }{ \lbfL_4 }
	}
	}

Since ${\fcMZBO{1}}^{\otimes4}$ and $\fcMChBO{4}\cdot\fcMZO{1}\circ\mpWT{4}$ are invariant under $\gpSym{4}$,
	we have
	\envM
	{
		\racFA{ {\fcMZBO{1}}^{\otimes4} }{ \mpS{\gpSym{3}} }{\lbfL_4}
	}
	{
		6 \Fc{ {\fcMZBO{1}}^{\otimes4} }{\lbfL_4}
	}
	and
	\envMCm
	{
		\racFA{ \fcMChBO{4}\cdot\fcMZO{1}\circ\mpWT{4} }{ \mpS{\gpSym{3}} }{\lbfL_4}
	}
	{
		6 \Fc{ \fcMChBO{4}\cdot\fcMZO{1}\circ\mpWT{4} }{\lbfL_4}
	}	
	which together with \refEq{4.2_Lem3Proof_EqRight1} and \refEq{4.2_Lem3Proof_EqRight4} prove \refEq{4.2_Lem3_Eq1} and \refEq{4.2_Lem3_Eq5},
	respectively.

We now prove \refEq{4.2_Lem3_Eq2}.	
We can easily see that
	\envMThCm{
		\mpS{\gpSym{4}}
	}
	{
		\mpS{\gpC{4}} \mpS{\gpSym{3}}
	}
	{
		\mpS{\gpSym{2}}\mpS{\gpAlt{4}}
	}
	and so,
	by \refEq{4.2_Lem2Proof_EqPL1}, 
	\envHLineFPd[4.2_Lem3Pr_Eqq2.1]
	{
		\racF{ \fcMZBO{2}\otimes{\fcMZBO{1}}^{\otimes2} }{ \mpS{\gpC{4}}\mpS{\gpSym{3}} }
	}
	{
		\bkR{ \racF{ \fcMZBO{2}\otimes{\fcMZBO{1}}^{\otimes2} }{ \mpS{\gpSym{2}} } } | \mpS{\gpAlt{4}}
	}
	{
		\racFr{ {\fcMZBO{1}}^{\otimes4} - (\fcMChBO{2} \cdot \fcMZO{1} \circ \mpWT{2}) \otimes{\fcMZBO{1}}^{\otimes2} }{ \mpS{\gpAlt{4}} }
	}
	{
		12 {\fcMZBO{1}}^{\otimes4} - \racF{ \bkR{ \fcMChBO{2} \cdot \fcMZO{1} \circ \mpWT{2}} \otimes{\fcMZBO{1}}^{\otimes2} }{ \mpS{\gpAlt{4}} }
	}
Since $\mpS{\gpAlt{4}}=\mpS{\Gp{(12)(34)}}\mpS{\gpAlts{4}}$ and
	$\bkR{ \fcMChBO{2} \cdot \fcMZO{1} \circ \mpWT{2}} \otimes{\fcMZBO{1}}^{\otimes2}$ is invariant under $\Gp{(12)(34)}$,
	\envHLinePd[4.2_Lem3Pr_Eqq2.2]
	{
		\racF{ \bkR{ \fcMChBO{2} \cdot \fcMZO{1} \circ \mpWT{2}} \otimes{\fcMZBO{1}}^{\otimes2} }{ \mpS{\gpAlt{4}} }
	}
	{
		2 \racF{ \bkR{ \fcMChBO{2} \cdot \fcMZO{1} \circ \mpWT{2}} \otimes{\fcMZBO{1}}^{\otimes2} }{ \mpS{\gpAlts{4}} }
	}
Combining \refEq{4.2_Lem3Pr_Eqq2.1} and \refEq{4.2_Lem3Pr_Eqq2.2},
	we obtain
	\envHLineCm
	{
		\racF{ \fcMZBO{2}\otimes{\fcMZBO{1}}^{\otimes2} }{ \mpS{\gpC{4}}\mpS{\gpSym{3}} }
	}
	{
		12 {\fcMZBO{1}}^{\otimes4} - 2 \racF{ \bkR{ \fcMChBO{2} \cdot \fcMZO{1} \circ \mpWT{2}} \otimes{\fcMZBO{1}}^{\otimes2} }{ \mpS{\gpAlts{4}} }
	}
	which together with \refEq{4.2_Lem3Proof_EqRight1} and \refEq{4.2_Lem3Proof_EqRight2} proves \refEq{4.2_Lem3_Eq2}.

We now prove \refEq{4.2_Lem3_Eq3}.
For this,
	we require the identity 
	\envHLineCm[4.2_Lem3Pr_Eq3.1]
	{
		\racF{ {\fcMZBO{2}}^{\otimes2} }{ \mpS{\gpCs{4}}\mpS{\gpSym{3}} }
	}
	{
		\bkR{ \racF{ \fcMZBO{2} }{ \mpS{\Gp{(12)}} } \otimes \racF{ \fcMZBO{2} }{ \mpS{\Gp{(34)}} } } | \mpS{\gpC{3}} 
	}
	which can be verified as follows.
A direct calculation shows that
	\envHLineThPd
	{
		\mpS{\Gp{(12), (34)}} \mpS{\gpC{3}}
 	}
	{
		\bkR{\gpu + (12) + (34) + (12)(34)} \bkR{\gpu + (123) + (132)}
	}
	{
		\gpu + (12) + (13) + (23) + (34) + (12)(34) 
		\lnAH
		+ 
		(123) + (132) + (143) + (243) + (1243) + (1432)
	}
From this equation and the equivalence classes modulo $\Gp{(13)(24)}$ in \refTab{3.1_PL_Table1},
	we see that
	\envHLinePdPt{\equiv}
	{
		\mpS{\gpSym{4}}
	}
	{
		2 \mpS{\Gp{(12), (34)}} \mpS{\gpC{3}}
		\quad
		\modd \Gp{(13)(24)}
	}
We also see from 
	\envM{
		\mpS{\gpSym{4}}
	}
	{
		\mpS{\gpC{4}}\mpS{\gpSym{3}}
	}
	and
	\refEq{4.1_1_Thm1PrDp4_EqHar22help1} 
	that
	\envMCmPt{\equiv}
	{
		\mpS{\gpSym{4}}
	}
	{
		2\mpS{\gpCs{4}}\mpS{\gpSym{3}}
		\;\modd \Gp{(13)(24)}
	}
	and so
	\envHLineCmPt{\equiv}
	{
		\mpS{\gpCs{4}}\mpS{\gpSym{3}}
	}
	{
		\mpS{\Gp{(12), (34)}} \mpS{\gpC{3}}
		\quad
		\modd \Gp{(13)(24)}
	}
	which verifies \refEq{4.2_Lem3Pr_Eq3.1},
	since
	${\fcMZBO{2}}^{\otimes2}$ is invariant under $\Gp{(13)(24)}$.
Then,	
	by \refEq{4.2_Lem2Proof_EqPL1},
	the right-hand side of \refEq{4.2_Lem3Pr_Eq3.1} is calculated as
	\envLineFPd
	{
		\text{(RHS of \refEq{4.2_Lem3Pr_Eq3.1})}
	}
	{
		\bkR{{\fcMZBO{1}}^{\otimes2} - \fcMChBO{2} \cdot \fcMZO{1} \circ \mpWT{2}}
		\otimes
		\bkR{{\fcMZBO{1}}^{\otimes2} - \fcMChBO{2} \cdot \fcMZO{1} \circ \mpWT{2}}
		| \mpS{\gpC{3}}
	}
	{
		\racF{\bkB{
			{\fcMZBO{1}}^{\otimes4} 
			- 
			(\fcMChBO{2} \cdot \fcMZO{1} \circ \mpWT{2}) \otimes {\fcMZBO{1}}^{\otimes2}
			-
			{\fcMZBO{1}}^{\otimes2} \otimes (\fcMChBO{2} \cdot \fcMZO{1} \circ \mpWT{2})
			+
			(\fcMChBO{2} \cdot \fcMZO{1} \circ \mpWT{2})^{\otimes2}
		}}{ \mpS{\gpC{3}} }
		\nonumber
	}
	{
		3 {\fcMZBO{1}}^{\otimes4} 
		-
		\racF{ (\fcMChBO{2} \cdot \fcMZO{1} \circ \mpWT{2}) \otimes {\fcMZBO{1}}^{\otimes2} }{ \mpS{\Gp{(13)(24)}} \mpS{\gpC{3}} }
		+
		\racF{ (\fcMChBO{2} \cdot \fcMZO{1} \circ \mpWT{2})^{\otimes2} }{ \mpS{\gpC{3}} }
	}
Since $\mpS{\Gp{(13)(24)}}\mpS{\gpC{3}}=\mpS{\gpAlts{4}}$,
	\refEq{4.2_Lem3Pr_Eq3.1} can be restated as
	\envHLineCm
	{
		\racF{ {\fcMZBO{2}}^{\otimes2} }{ \mpS{\gpCs{4}}\mpS{\gpSym{3}} }
	}
	{
		3 {\fcMZBO{1}}^{\otimes4} 
		-
		\racF{ (\fcMChBO{2} \cdot \fcMZO{1} \circ \mpWT{2}) \otimes {\fcMZBO{1}}^{\otimes2} }{ \mpS{\gpAlts{4}} }
		+
		\racF{ (\fcMChBO{2} \cdot \fcMZO{1} \circ \mpWT{2})^{\otimes2} }{ \mpS{\gpC{3}} }
	}
	which together with \refEq{4.2_Lem3Proof_EqRight1}, \refEq{4.2_Lem3Proof_EqRight2}, and \refEq{4.2_Lem3Proof_EqRight3.1} proves \refEq{4.2_Lem3_Eq3}.

We lastly prove \refEq{4.2_Lem3_Eq4} in a similar way to \refEq{4.2_Lem3_Eq3}.
We require the identity
	\envHLineCm[4.2_Lem3Pr_Eq4]
	{\hspace{-15pt}
		\racF{ \fcMZBO{3}\otimes\fcMZBO{1} }{ \mpS{\gpC{4}}\mpS{\gpSym{3}} }
	}
	{
		4 {\fcMZBO{1}}^{\otimes4} 
		- 
		2 \racF{ \fcMChBO{2}\cdot(\fcMZBO{2}\circ\mpWT{2})\otimes{\fcMZBO{1}}^{\otimes2} }{ \mpS{\gpAlts{4}} } 
		+
		2 \racF{ (\fcMChBO{3}\cdot\fcMZO{1}\circ\mpWT{3})\otimes\fcMZBO{1} }{ \mpS{\gpC{4}} } 
	}
	which can be verified as follows.
Identity \refEq{1_Cor_EqSymRMZV} for $n=3$ can be restated as
	\envHLineCm
	{
		\racF{\fcMZBO{3}}{ \mpS{\gpSym{3}} }
	}
	{
		{\fcMZBO{1}}^{\otimes3} 
		-
		\racF{ ( \fcMChBO{2}\cdot\fcMZBO{1}\circ\mpWT{2} ) \otimes\fcMZBO{1} }{ \mpS{\gpC{3}} }
		+
		2 \fcMChBO{3}\cdot\fcMZO{1}\circ\mpWT{3}
	}	
	because of \refEq{4.2_Lem2Proof_EqRight1}, \refEq{4.2_Lem2Proof_EqRight2}, and \refEq{4.2_Lem2Proof_EqRight3}.
A direct calculation shows that
	\envHLineCm
	{
		\mpS{\gpC{3}}\mpS{\gpC{4}}
	}
	{
		\gpu + (14) + (34) + (13)(24) + (123) + (124) + (132) + (243)
		\lnAH
		+
		(1234) + (1342) + (1423) + (1432)
	}
	and so 
	we see from the equivalence classes modulo $\Gp{(12),(34)}$ in \refTab{3.1_PL_Table1} that
	\envHLinePdPt{\equiv}
	{
		\mpS{\gpC{3}}\mpS{\gpC{4}}
	}
	{
		2 \mpS{\gpAlts{4}}
		\qquad\modd 
		\Gp{(12),(34)}
	} 
Since
	\envMCm{
		\mpS{\gpC{4}}\mpS{\gpSym{3}}
	}{
		\mpS{\gpSym{3}}\mpS{\gpC{4}}
	}
	we thus have
	\envHLineFCm
	{
		\racF{ \fcMZBO{3}\otimes\fcMZBO{1} }{ \mpS{\gpC{4}}\mpS{\gpSym{3}} }
	}
	{
		\racF{ \racFrr{ \fcMZBO{3} }{\mpS{\gpSym{3}}}\otimes\fcMZBO{1}  }{ \mpS{\gpC{4}} }
	}
	{
		\racF{ {\fcMZBO{1}}^{\otimes4} }{ \mpS{\gpC{4}} }
		-
		\racF{ (\fcMChBO{2}\cdot\fcMZBO{1}\circ\mpWT{2} ) \otimes{\fcMZBO{1}}^{\otimes2} }{ \mpS{\gpC{3}}\mpS{\gpC{4}} }
		+
		2 \racF{ (\fcMChBO{3}\cdot\fcMZO{1}\circ\mpWT{3})\otimes\fcMZBO{1} }{ \mpS{\gpC{4}} }
	}
	{
		4 {\fcMZBO{1}}^{\otimes4} 
		- 
		2 \racF{ \fcMChBO{2}\cdot(\fcMZBO{2}\circ\mpWT{2})\otimes{\fcMZBO{1}}^{\otimes2} }{ \mpS{\gpAlts{4}} } 
		+
		2 \racF{ (\fcMChBO{3}\cdot\fcMZO{1}\circ\mpWT{3})\otimes\fcMZBO{1} }{ \mpS{\gpC{4}} } 
	}
	which verifies \refEq{4.2_Lem3Pr_Eq4}.
Then,
	combining \refEq{4.2_Lem3Proof_EqRight1}, \refEq{4.2_Lem3Proof_EqRight2}, \refEq{4.2_Lem3Proof_EqRight3.2}, and \refEq{4.2_Lem3Pr_Eq4},
	we obtain \refEq{4.2_Lem3_Eq4},
	which completes the proof.
}

\begin{table}[!t]
\begin{center}{\tiny
	\caption{Examples of \refEq{1_Thm1_EqDpt3} (or \refEq{5_PL_EqDepth3}) }\label{5_PL_TableExThmDp3}
\mbox{}\\
	\begin{tabular}{|c|lr|}  \hline
		\parbox[c][15pt][c]{0pt}{}
		{\bf Index set}&{\bf Linear relation}&\\\hline \parbox[c][15pt][c]{0pt}{}
		(1,1,1)&$3\fcMZB{3}{1,1,1} = \fcMChB{3}{1,1,1} \fcMZ{1}{3}$.
			&(d3-1)   		\\\hline \parbox[c][15pt][c]{0pt}{}
		(1,1,2)&$\fcMZB{3}{1,1,2} + \fcMZB{3}{1,2,1} + \fcMZ{3}{2,1,1} = \fcMZB{2}{1,1}\fcMZ{1}{2} + \fcMZ{1}{4}$.
			&(d3-2)   		\\\hline \parbox[c][15pt][c]{0pt}{}
		(1,1,3)&$\fcMZB{3}{1,1,3} + \fcMZB{3}{1,3,1} + \fcMZ{3}{3,1,1} = \fcMZB{2}{1,1}\fcMZ{1}{3} + \fcMZ{1}{5}$,
			&(d3-3)   	\\ \parbox[c][15pt][c]{0pt}{}
		(1,2,2)&$\fcMZB{3}{1,2,2} + \fcMZ{3}{2,2,1} + \fcMZ{3}{2,1,2} = - \fcMZ{1}{2}\fcMZ{1}{3} + \fcMZ{1}{5}$.
			&(d3-4)		\\\hline \parbox[c][15pt][c]{0pt}{}
		(1,1,4)&$\fcMZB{3}{1,1,4} + \fcMZB{3}{1,4,1} + \fcMZ{3}{4,1,1} = \fcMZB{2}{1,1}\fcMZ{1}{4} + \fcMZ{1}{6}$,
			&(d3-5)		\\ \parbox[c][15pt][c]{0pt}{}
		(1,2,3)&$\fcMZB{3}{1,2,3} + \fcMZ{3}{2,3,1} + \fcMZ{3}{3,1,2} = \fcMZB{2}{1,2}\fcMZ{1}{3} + \fcMZ{2}{3,1}\fcMZ{1}{2} + \fcMZ{1}{6}$,
			&(d3-6)     	\\ \parbox[c][15pt][c]{0pt}{}
		(1,3,2)&$\fcMZB{3}{1,3,2} + \fcMZ{3}{3,2,1} + \fcMZ{3}{2,1,3} = \fcMZB{2}{1,3}\fcMZ{1}{2} + \fcMZ{2}{2,1}\fcMZ{1}{3} + \fcMZ{1}{6}$,
			&(d3-7)    	\\ \parbox[c][15pt][c]{0pt}{}
		(2,2,2)&$3\fcMZ{3}{2,2,2} = - \fcMZ{1}{2}^3 + 3\fcMZ{2}{2,2}\fcMZ{1}{2} + \fcMZ{1}{6}$.
			&(d3-8)     	\\\hline
	\end{tabular} 
}\end{center}
\begin{center}{\tiny
	\caption{Examples of \refEq{1_Cor_EqSymRMZV} for $n=3$ }\label{5_PL_TableExCorDp3}
\mbox{}\\
	\begin{tabular}{|c|lr|}  \hline
		\parbox[c][15pt][c]{0pt}{}
		{\bf Index set}&{\bf Linear relation}&\\\hline \parbox[c][15pt][c]{0pt}{}
		(1,1,1)&$6\fcMZB{3}{1,1,1} = 2\fcMChB{3}{1,1,1} \fcMZ{1}{3}$.
			&(d3'-1)   		\\\hline \parbox[c][15pt][c]{0pt}{}
		(1,1,2)&$2(\fcMZB{3}{1,1,2} + \fcMZB{3}{1,2,1} + \fcMZ{3}{2,1,1}) = - \fcMChB{2}{1,1}\fcMZ{1}{2}^2 + 2\fcMZ{1}{4}$.
			&(d3'-2)   		\\\hline \parbox[c][15pt][c]{0pt}{}
		(1,1,3)&$2(\fcMZB{3}{1,1,3} + \fcMZB{3}{1,3,1} + \fcMZ{3}{3,1,1}) = - \fcMChB{2}{1,1}\fcMZ{1}{2}\fcMZ{1}{3} + 2\fcMZ{1}{5}$,
			&(d3'-3)   	\\ \parbox[c][15pt][c]{0pt}{}
		(1,2,2)&$2(\fcMZB{3}{1,2,2} + \fcMZ{3}{2,2,1} + \fcMZ{3}{2,1,2}) = - 2\fcMZ{1}{2}\fcMZ{1}{3} + 2\fcMZ{1}{5}$.
			&(d3'-4)		\\\hline \parbox[c][15pt][c]{0pt}{}
		(1,1,4)&$2(\fcMZB{3}{1,1,4} + \fcMZB{3}{1,4,1} + \fcMZ{3}{4,1,1}) = - \fcMChB{2}{1,1}\fcMZ{1}{2}\fcMZ{1}{4} + 2\fcMZ{1}{6}$,
			&(d3'-5)		\\ \parbox[c][15pt][c]{0pt}{}
		(1,2,3)&$\fcMZB{3}{1,2,3} + \fcMZB{3}{1,3,2} + \fcMZ{3}{2,1,3} +\fcMZ{3}{2,3,1} + \fcMZ{3}{3,1,2} + \fcMZ{3}{3,2,1} $
			&\\ \quad \parbox[c][10pt][c]{0pt}{}
			&$= -(\fcMZ{1}{2}\fcMZ{1}{4} + \fcMZ{1}{3}^2 ) + 2\fcMZ{1}{6}$,
			&(d3'-6)     	\\ \parbox[c][20pt][c]{0pt}{}
		(2,2,2)&$6\fcMZ{3}{2,2,2} = \fcMZ{1}{2}^3 - 3\fcMZ{1}{2}\fcMZ{1}{4} + 2\fcMZ{1}{6}$.
			&(d3'-7)     	\\\hline
	\end{tabular} 
}\end{center}
\begin{center}{\tiny
	\caption{Examples of \refEq{1_Thm1_EqDpt4} (or \refEq{5_PL_EqDepth4}) }\label{5_PL_TableExThmDp4}
\mbox{}\\
	\begin{tabular}{|c|lc|}  \hline
		\parbox[c][15pt][c]{0pt}{}
		{\bf Index set}&{\bf Linear relation}&\\\hline \parbox[c][15pt][c]{0pt}{}
		(1,1,1,1)&$4\fcMZB{4}{1,1,1,1} =  2\fcMZB{2}{1,1}^2 -  \fcMChB{4}{1,1,1,1} \fcMZ{1}{4} $.
			&(d4-1)   		\\\hline \parbox[c][15pt][c]{0pt}{}
		(1,1,1,2)&$\fcMZB{4}{1,1,1,2}+\fcMZB{4}{1,1,2,1}+\fcMZB{4}{1,2,1,1}+\fcMZ{4}{2,1,1,1} =-\fcMZB{2}{1,1}\fcMZ{1}{3}+\fcMZB{3}{1,1,1}\fcMZ{1}{2}-\fcMZ{1}{5}$.
			&(d4-2)   	\\\hline \parbox[c][15pt][c]{0pt}{}
		(1,1,1,3)&$\fcMZB{4}{1,1,1,3}+\fcMZB{4}{1,1,3,1}+\fcMZB{4}{1,3,1,1}+\fcMZ{4}{3,1,1,1} =-\fcMZB{2}{1,1}\fcMZ{1}{4}+\fcMZB{3}{1,1,1}\fcMZ{1}{3}-\fcMZ{1}{6}$,
			&(d4-3)   	\\  \parbox[c][15pt][c]{0pt}{} 
		(1,1,2,2)&$\fcMZB{4}{1,1,2,2}+\fcMZ{4}{1,2,2,1}+\fcMZB{4}{2,2,1,1}+\fcMZB{4}{2,1,1,2}$ 
			&\\ \quad \parbox[c][10pt][c]{0pt}{}
			&$=-\fcMZB{2}{1,1}\fcMZ{1}{2}^2 + \fcMZB{2}{1,1}\fcMZ{2}{2,2} + \fcMZB{2}{1,2}\fcMZ{2}{2,1} + \bkR[n]{\fcMZB{3}{1,1,2}+\fcMZ{3}{2,1,1}}\fcMZ{1}{2} - \fcMZ{1}{6}$,
			&(d4-4)   	\\  \parbox[c][20pt][c]{0pt}{}
		(1,2,1,2)&$2(\fcMZB{4}{1,2,1,2}+\fcMZ{4}{2,1,2,1}) = \fcMZB{2}{1,2}^2 + \fcMZ{2}{2,1}^2 + 2\fcMZB{3}{1,2,1}\fcMZ{1}{2} - \fcMZ{1}{6}$.
			&(d4-5)   	\\\hline
	\end{tabular} 
}\end{center}
\begin{center}{\tiny
	\caption{Examples of \refEq{1_Cor_EqSymRMZV} for $n=4$ }\label{5_PL_TableExCorDp4}
\mbox{}\\
	\begin{tabular}{|c|lc|}  \hline
		\parbox[c][15pt][c]{0pt}{}
		{\bf Index set}&{\bf Linear relation}&\\\hline \parbox[c][15pt][c]{0pt}{}
		(1,1,1,1)&$24\fcMZB{4}{1,1,1,1} =   3\fcMChB{2}{1,1}\fcMZ{1}{2}^2  - 6\fcMChB{4}{1,1,1,1} \fcMZ{1}{4} $.
			&(d4'-1)   		\\\hline \parbox[c][15pt][c]{0pt}{}
		(1,1,1,2)&$6(\fcMZB{4}{1,1,1,2}+\fcMZB{4}{1,1,2,1}+\fcMZB{4}{1,2,1,1}+\fcMZ{4}{2,1,1,1}) $
			&\\ \quad \parbox[c][15pt][c]{0pt}{}
			&$=3\fcMChB{2}{1,1}\fcMZ{1}{2}\fcMZ{1}{3} + 2\fcMChB{3}{1,1,1}\fcMZ{1}{2}\fcMZ{1}{3} - 6\fcMZ{1}{5}$.
			&(d4'-2)   	\\\hline \parbox[c][15pt][c]{0pt}{}
		(1,1,1,3)&$6(\fcMZB{4}{1,1,1,3}+\fcMZB{4}{1,1,3,1}+\fcMZB{4}{1,3,1,1}+\fcMZ{4}{3,1,1,1})$ 
			&\\ \quad \parbox[c][15pt][c]{0pt}{}
			&$=3\fcMChB{2}{1,1}\fcMZ{1}{2}\fcMZ{1}{4} + 2\fcMChB{3}{1,1,1}\fcMZ{1}{3}^2 - 6\fcMZ{1}{6}$,
			&(d4'-3)   	\\  \parbox[c][15pt][c]{0pt}{} 
		(1,1,2,2)&$4(\fcMZB{4}{1,1,2,2}+\fcMZ{4}{1,2,1,2}+\fcMZ{4}{1,2,2,1}+\fcMZB{4}{2,1,1,2}+\fcMZB{4}{2,1,2,1}+\fcMZB{4}{2,2,1,1})$ 
			&\\ \quad \parbox[c][15pt][c]{0pt}{}
			&$=-\fcMChB{2}{1,1}\fcMZ{1}{2}^3 +(\fcMChB{2}{1,1}+4)\fcMZ{1}{2}\fcMZ{1}{4} + 2\fcMZ{1}{3}^2 -6\fcMZ{1}{6}$.
			&(d4'-4)   	\\   \hline
	\end{tabular} 
}\end{center}
\end{table}

\section{Examples}\label{sectFive}
We list examples of \refEq{1_Thm1_EqDpt3} and \refEq{1_Thm1_EqDpt4} in \refTab{5_PL_TableExThmDp3} and \refTab{5_PL_TableExThmDp4},
	respectively.
We also list examples of \refEq{1_Cor_EqSymRMZV} for $n=3$ and $n=4$ in \refTab{5_PL_TableExCorDp3} and \refTab{5_PL_TableExCorDp4}, respectively, for comparison.
The examples treat the case of weights less than $7$.
We omit examples of \refEq{1_Thm1_EqDpt2} and \refEq{1_Cor_EqSymRMZV} for $n=2$
	because they are essentially the harmonic relations.
The following straightforward expressions of \refEq{1_Thm1_EqDpt3} and \refEq{1_Thm1_EqDpt4} 
	are convenient for calculating the examples in \refTab{5_PL_TableExThmDp3} and \refTab{5_PL_TableExThmDp4}:
	\envLineCm[5_PL_EqDepth3]
	{
		\fcMZB{3}{l_1,l_2,l_3} + \fcMZB{3}{l_2,l_3,l_1} + \fcMZB{3}{l_3,l_1,l_2}
	}
	{
		-
		\fcMZB{1}{l_1}\fcMZB{1}{l_2}\fcMZB{1}{l_3}
		\lnAH
		+
		\fcMZB{2}{l_1,l_2}\fcMZB{1}{l_3} + \fcMZB{2}{l_2,l_3}\fcMZB{1}{l_1} + \fcMZB{2}{l_3,l_1}\fcMZB{1}{l_2}
		\lnAH[]
		+
		\fcMChB{3}{l_1,l_2,l_3} \fcMZ{1}{l_1+l_2+l_3}
		\nonumber
	}
	and
	\envLinePd[5_PL_EqDepth4]
	{
		\fcMZB{4}{l_1,l_2,l_3,l_4} + \fcMZB{4}{l_2,l_3,l_4,l_1} + \fcMZB{4}{l_3,l_4,l_1,l_2} + \fcMZB{4}{l_4,l_1,l_2,l_3}
	}
	{
		\fcMZB{1}{l_1}\fcMZB{1}{l_2}\fcMZB{1}{l_3}\fcMZB{1}{l_4}
		\lnAH
		-
		\fcMZB{2}{l_1,l_2}\fcMZB{1}{l_3}\fcMZB{1}{l_4} - \fcMZB{2}{l_2,l_3}\fcMZB{1}{l_4}\fcMZB{1}{l_1}
		\lnAH
		-
		\fcMZB{2}{l_3,l_4}\fcMZB{1}{l_1}\fcMZB{1}{l_2} - \fcMZB{2}{l_4,l_1}\fcMZB{1}{l_2}\fcMZB{1}{l_3} 
		\lnAH[]
		+
		\fcMZB{2}{l_1,l_2}\fcMZB{2}{l_3,l_4} + \fcMZB{2}{l_2,l_3}\fcMZB{2}{l_4,l_1}
		\lnAH
		+
		\fcMZB{3}{l_1,l_2,l_3}\fcMZB{1}{l_4} + \fcMZB{3}{l_2,l_3,l_4}\fcMZB{1}{l_1} + \fcMZB{3}{l_3,l_4,l_1}\fcMZB{1}{l_2} + \fcMZB{3}{l_4,l_1,l_2}\fcMZB{1}{l_3}
		\lnAH
		-
		\fcMChB{4}{l_1,l_2,l_3,l_4} \fcMZ{1}{l_1+l_2+l_3+l_4}
		\nonumber
	}
We note that
	we have used
	\envM
	{
		\fcMZB{1}{1}
	}
	{
		0
	}
	for all equations in the tables,
	and
	\envM
	{
		\fcMZB{2}{1,k}+\fcMZB{2}{k,1}
	}
	{
		-\fcMZ{1}{k+1}
		\;
		(k>1)
	} 
	for \refEqA{d3-4}, \refEqA{d4-2}, and \refEqA{d4-3}.

Lastly,
	we derive the following from (d3-1) and (d4-1) as applications of examples:
	\envHLineCSNmePd
	{\label{5_PL_EqSha}
		\fcMZS{3}{1,1,1}	\lnP{=}	\fcMZS{4}{1,1,1,1}
	}
	{
		0
	}
	{\label{5_PL_Eq111Har}
		\fcMZH{3}{1,1,1}
	}
	{
		\opF{1}{3}\fcMZ{1}{3}
	}
	{\label{5_PL_Eq1111Har}
		\fcMZH{4}{1,1,1,1}
	}
	{
		\opF{1}{16}\fcMZ{1}{4}
	}
We can easily obtain \refEq{5_PL_EqSha} from (d3-1) and (d4-1), 
	since 	
	\envM
	{
		\fcMChS{n}{1,\ldots,1}
	}
	{
		0
	}
	and
	\envMPd
	{
		\fcMZS{2}{1,1}
	}
	{
		0
	}
We can also obtain \refEq{5_PL_Eq111Har} from (d3-1) because  
	\envMPd
	{
		\fcMChH{n}{1,\ldots,1}
	}
	{
		1
	}
We see from \refEq{3.2_Lem3iiEq3} and (d4-1) that
	\envMThCm
	{
		4 \fcMZH{4}{1,1,1,1}
	}
	{
		2\fcMZH{2}{1,1}^2 - \fcMZ{1}{4}
	}
	{
		\opF{1}{2}\fcMZ{1}{2}^2 - \fcMZ{1}{4}
	}
	which together with \refEq{3.2_Prop1Pr_EqDp4_3} proves \refEq{5_PL_Eq1111Har}.




\end{document}